\documentclass[12pt,british,reqno,openany]{amsbook}

\usepackage{babel}
\usepackage{amssymb,amsthm}
\usepackage[latin1]{inputenc}
\usepackage[T1]{fontenc}
\usepackage{textcomp}
\usepackage{eucal}
\usepackage{enumerate}
\usepackage{calc}

\usepackage[matrix,arrow]{xy}

\makeatletter
\AtBeginDocument{%
  \renewcommand\mainmatter{\cleardoublepage\pagenumbering{arabic}}%
  \renewcommand\backmatter{\if@openright\cleardoublepage\else\clearpage\fi}%
  \renewcommand\chapter{\if@openright\cleardoublepage\else\clearpage\fi
    \thispagestyle{plain}\global\@topnum\z@
    \@afterindenttrue \secdef\@chapter\@schapter}%
  }

\def\maketitle{\par
  \@topnum\z@ 
  \begingroup
  \@maketitle
  \endgroup
  \def\do##1{\let##1\relax}%
  \do\maketitle \do\@maketitle \do\title \do\@xtitle \do\@title
  \do\author \do\@xauthor \do\address \do\@xaddress
  \do\email \do\@xemail \do\curraddr \do\@xcurraddr
  \do\dedicatory \do\@dedicatory
}
\def\tableofcontents{%
  \@starttoc{toc}\contentsname
  \@setabstract
}
\def\@setabstract{%
  \ifdim\dp\abstractbox>\z@
    {\ifdim\ht\abstractbox<38pc \@openrightfalse \fi
     \def\@tocwriteb##1##2##3{}%
     \chapter*{Abstract}
     \@maketopfootnotes
     \noindent\unvbox\abstractbox
    }%
  \fi
  \c@footnote\z@
}
\def\@maketopfootnotes{%
  {\let\@makefnmark\relax  \let\@thefnmark\relax
   \ifx\@empty\relax\else \@footnotetext{\relax}\fi
   \ifx\@empty\@subjclass\else \@footnotetext{\@setsubjclass}\fi
   \ifx\@empty\@keywords\else \@footnotetext{\@setkeywords}\fi
   \ifx\@empty\thankses\else \@footnotetext{%
     \def\par{\let\par\@par}\@setthanks}\fi
  }%
  \def\do##1{\let##1\relax}%
  \do\keywords \do\@keywords \do\subjclass \do\@subjclass
}
\def\@setsubjclass{{\let\\=\space
  {\itshape\subjclassname.}\enspace\@subjclass\@addpunct.}}
\def\@setkeywords{%
  {\itshape \keywordsname.}\enspace \@keywords\@addpunct.}
\def\@maketitle{%
  \cleardoublepage \thispagestyle{empty}%
  \begingroup \topskip\z@skip
  \null\vfil
  \begingroup
  \LARGE\bfseries \centering
  \openup\medskipamount
  \@title\par\vspace{24pt}%
  \def\and{\par\medskip}\centering
  \mdseries\authors\par\bigskip
  \endgroup
  \vfil
  \ifx\@empty\addresses\else
  Author address:
  \@setaddresses\fi
  \bigskip
  \vfil
  \begin{center}
    \ifx\@empty\@translators\else\vfil\@settranslators\fi
  \end{center}
  \vfil\vfil
  \endgroup
}
\renewenvironment{abstract}{%
  \ifx\maketitle\relax
    \ClassWarning{\@classname}{Abstract should precede
      \protect\maketitle\space in AMS documentclasses; reported}%
  \fi
  \global\setbox\abstractbox=\vtop\bgroup
      \normalsize
      \indent
}{%
    \egroup
}

\renewcommand{\p@enumii}{}

\def\@enum@{\list{\csname label\@enumctr\endcsname}%
           {\usecounter{\@enumctr}\def\makelabel##1{
\normalfont\ignorespaces\emph{{##1}~}}
\setlength{\labelsep}{0pt}
\settowidth{\labelwidth}{0pt}
\setlength{\leftmargin}{\labelwidth}
\addtolength{\leftmargin}{\labelsep}
\setlength{\itemindent}{-4pt}
 \setlength{\leftmargini}{\labelwidth}
\addtolength{\leftmargini}{-4pt}
\setlength{\leftmarginii}{\leftmargini}
\setlength{\leftmarginiii}{\leftmargini}
}}

\makeatother

\renewcommand{\epsilon}{\ensuremath{\varepsilon}}
\renewcommand{\phi}{\ensuremath{\varphi}}
\newcommand{\vide}{\ensuremath{\varnothing}}
\renewcommand{\to}{\ensuremath{\rightarrow}}
\newcommand{\comm}{\ensuremath{\rightleftharpoons}}

\newcommand{\R}{\ensuremath{\mathbb R}}

\newcommand{\N}{\ensuremath{\mathbb N}}
\newcommand{\Z}{\ensuremath{\mathbb Z}}

\newcommand{\dt}{\ensuremath{\mathbb D}^{2}}
\newcommand{\St}[1][2]{\ensuremath{\mathbb S}^{#1}}
\newcommand{\FF}{\ensuremath{\mathbb F}}
\newcommand{\F}[1][n]{\ensuremath{\FF_{{#1}}}}
\newcommand{\sn}[1][n]{\ensuremath{S_{{#1}}}}
\newcommand{\an}[1][n]{\ensuremath{A_{{#1}}}}
\DeclareMathOperator{\id}{\text{Id}}

\newcommand{\pgcd}[2]{\ensuremath{\operatorname{\text{gcd}}(#1,#2)}}
\newcommand{\ppcm}[2]{\ensuremath{\operatorname{\text{lcm}}(#1,#2)}}
\newcommand{\tr}[1]{\ensuremath{\operatorname{\text{Tr}}{#1}}}
\renewcommand{\ker}[1]{\ensuremath{\operatorname{\text{Ker}}\left({#1}\right)}}
\newcommand{\im}[1]{\ensuremath{\operatorname{\text{Im}}\left({#1}\right)}}
\newcommand{\coker}[1]{\ensuremath{\operatorname{\text{Coker}}\left({#1}\right)}}
\newcommand{\aut}[1]{\ensuremath{\operatorname{\text{Aut}}({#1})}}
\newcommand{\Int}[1]{\ensuremath{\operatorname{\text{Int}}({#1})}}

\newcommand{\pnm}[1][n]{\ensuremath{P_{{#1}}(M)}}
\newcommand{\gpab}[1][G]{\ensuremath{{#1}\textsuperscript{Ab}}}
\newcommand{\bnab}[1]{\ensuremath{\left(B_{#1}(\St)\right)\textsuperscript{Ab}}}
\newcommand{\bmmn}[2]{\ensuremath{B_{#1}(\St\setminus\brak{x_1,\ldots,x_{#2}})}}

\makeatletter
\def\@map#1#2[#3]{\mbox{$#1 \colon\thinspace #2 \to #3$}}
\def\map#1#2{\@ifnextchar [{\@map{#1}{#2}}{\@map{#1}{#2}[#2]}}
\DeclareRobustCommand*\textsubscript[1]{\@textsubscript{\selectfont#1}}
\def\@textsubscript#1{{\m@th\ensuremath{_{\mbox{\fontsize\sf@size\z@#1}}}}}
\makeatother

\newcommand{\lhra}{\mathrel{\lhook\joinrel\to}}

\DeclareRobustCommand*{\up}[1]{\textsuperscript{#1}}

\newcommand{\ft}[1][n]{\ensuremath{\Delta_{#1}}}

\newcommand{\brak}[1]{\ensuremath{\left\{ #1 \right\}}}
\newcommand{\ang}[1]{\ensuremath{\langle #1\rangle}}
\newcommand{\set}[2]{\ensuremath{\brak{#1 \,\mid\, #2}}}
\newcommand{\setang}[2]{\ensuremath{\ang{#1 \,\mid\, #2}}}

\newcommand{\si}[2][{}]{\ensuremath{\sigma_{#2}^{#1}}}
\newcommand{\sii}[2][1]{\ensuremath{\sigma_{#2}^{-{#1}}}}

\theoremstyle{plain}
\newtheorem{thm}{Theorem}
\newtheorem{lem}[thm]{Lemma}
\newtheorem{prop}[thm]{Proposition}
\newtheorem{cor}[thm]{Corollary}

\theoremstyle{remark}
\newtheorem{rem}[thm]{Remark}
\newtheorem{rems}[thm]{Remarks}

\newtheoremstyle{citing}
  {}
  {}
  {\itshape}
  {\parindent}
  {\scshape}
  {.}
  {.5em}
  {\thmnote{#3}}

\theoremstyle{citing}
\newtheorem*{varthm}{}

\newcommand{\reth}[1]{Theorem~\protect\ref{th:#1}}
\newcommand{\relem}[1]{Lemma~\protect\ref{lem:#1}}
\newcommand{\repr}[1]{Proposition~\protect\ref{prop:#1}}
\newcommand{\reco}[1]{Corollary~\protect\ref{cor:#1}}
\newcommand{\resec}[1]{Section~\protect\ref{sec:#1}}
\newcommand{\rechap}[1]{Chapter~\protect\ref{chap:#1}}
\newcommand{\rerem}[1]{Remark~\protect\ref{rem:#1}}
\newcommand{\rerems}[1]{Remarks~\protect\ref{rem:#1}}
\newcommand{\req}[2][{}]{equation~(\protect\ref{eq:#2}\textsubscript{${#1}$})}
\newcommand{\reqref}[2][{}]{(\protect\ref{eq:#2}\textsubscript{${#1}$})}
\newcommand{\unreqref}[2][{}]{\protect\ref{eq:#2}\textsubscript{${#1}$}}

\begin{document}

\frontmatter
\title{The lower central and derived series of the braid groups
$B_n(\St)$ and $\bmmn{m}{n}$}

\author{Daciberg~Lima~Gon\c{c}alves}
\address{Departamento de Matem\'atica - IME-USP,\newline
Caixa Postal~\textup{66281}~-~Ag.~Cidade de S\~ao Paulo,\newline CEP:~\textup{05311-970} - S\~ao Paulo - SP - Brazil.}
\email{dlgoncal@ime.usp.br}

\author{John~Guaschi\vspace{0.25cm}\linebreak
{\textmd{\large 17th March 2006}}}
\address{Laboratoire de Math\'ematiques Emile Picard,\newline 
UMR CNRS~\textup{5580}, UFR-MIG, Universit\'e Toulouse~III,\newline \textup{31062}~Toulouse Cedex~\textup{9}, France.}
\email{guaschi@picard.ups-tlse.fr}

\date{17\up{th}~March~2006}




\begin{abstract}
Our aim in this paper is to determine the lower central and derived
series for the braid groups of the sphere and of the
finitely-punctured sphere. We are motivated in part by the study of
the generalised Fadell-Neuwirth short exact sequence~\cite{GG2,GG4},
but the problem is of interest in its own right.

The braid groups of the $2$-sphere $\St$ were studied by Fadell, Van
Buskirk and Gillette during the 1960's, and are of particular interest
due to the fact that they have torsion elements (which were
characterised by Murasugi). We first prove that for all $n\in\N$, the
lower central series of the $n$-string braid group $B_n(\St)$ is
constant from the commutator subgroup onwards. We obtain a
presentation of $\Gamma_2(B_n(\St))$, from which we observe that
$\Gamma_2(B_4(\St))$ is a semi-direct product of the quaternion group
of order~$8$ by a free group of rank $2$. As for the derived series of
$B_n(\St)$, we show that for all $n\geq 5$, it is constant from the
derived subgroup onwards. The group $B_n(\St)$ being finite and
soluble for $n\leq 3$, the critical case is $n=4$ for which the
derived subgroup is the semi-direct product obtained above. By proving
a general result concerning the structure of the derived subgroup of a
semi-direct product, we are able to determine completely the derived
series of $B_4(\St)$ which from $(B_4(\St))^{(4)}$ onwards coincides
with that of the free group of rank $2$, as well as its successive
derived series quotients.

For $n \geq 1$, the class of $m$-string braid groups $\bmmn{m}{n}$ of
the $n$-punctured sphere includes the usual Artin braid groups $B_m$
(for $n=1$), those of the annulus, which are Artin groups of type $B$
(for $n=2$), and affine Artin groups of type $\widetilde{C}$ (for
$n=3$). Motivated by the study of almost periodic solutions of
algebraic equations with almost periodic coefficients, Gorin and Lin
determined the commutator subgroup of the Artin braid groups. We
extend their results, and show that the lower central series of $B_m$
is completely determined for all $m\in\N$, and that the derived series
is determined for all $m\neq 4$. In the exceptional case $m=4$, we
determine some higher elements of the derived series and its
quotients. 

When $n\geq 2$, we prove that the lower central series (respectively
derived series) of $\bmmn{m}{n}$ is constant from the commutator
subgroup onwards for all $m\geq 3$ (respectively $m\geq 5$). The case
$m=1$ is that of the free group of rank $n-1$. The case $n=2$ is of
particular interest notably when $m=2$ also. In this case, the
commutator subgroup is a free group of infinite rank.  We then go on
to show that $B_2(\St\setminus\brak{x_1,x_2})$ admits various
interpretations, as the Baumslag-Solitar group
$\operatorname{BS}(2,2)$, or as a one relator group with non-trivial
centre for example. We conclude from this latter fact that
$B_2(\St\setminus\brak{x_1,x_2})$ is residually nilpotent, and that
from the commutator subgroup onwards, its lower central series
coincides with that of the free product $\Z_2\ast \Z$. Further, its
lower central series quotients $\Gamma_i/\Gamma_{i+1}$ are direct sums
of copies of $\Z_2$, the number of summands being determined
explicitly. In the case $m\geq 3$ and $n=2$, we obtain a presentation
of the derived subgroup, from which we deduce its Abelianisation.
Finally, in the case $n=3$, we obtain partial results for the derived
series, and we prove that the lower central series quotients
$\Gamma_i/\Gamma_{i+1}$ are $2$-elementary finitely-generated groups. 
\end{abstract}

\maketitle

\tableofcontents

\chapter*{Preface}

\section{Generalities and definitions}

Let $n\in\N$. The braid groups of the plane $\mathbb{E}^2$, denoted by $B_n$, and known as \emph{Artin braid groups}, were introduced by E.~Artin in~1925~\cite{A1}, and further studied in~\cite{A2,A3,Ch}. Artin showed that $B_n$ admits the following well-known presentation: $B_n$ is generated by elements $\sigma_1,\ldots, \sigma_{n-1}$, subject to the classical \emph{Artin relations}:
\begin{equation}\label{eq:presartin}
\left.
\begin{aligned}
&\text{$\si{i}\si{j}=\si{j}\si{i}$ if $\lvert i-j\rvert\geq 2$ and $1\leq i,j\leq n-1$}\\
&\text{$\si{i}\si{i+1}\si{i}=\si{i+1}\si{i}\si{i+1}$ for all $1\leq i\leq n-2$.}
\end{aligned}
\right\}
\end{equation}
A natural generalisation to braid groups of arbitrary topological
spaces was made at the beginning of the 1960's by Fox (using the
notion of configuration space)~\cite{FoN}. In that paper, Fox and
Neuwirth proved some basic results about the braid groups of arbitrary
manifolds. In particular, if $M^r$ is a connected manifold of
dimension $r\geq 3$ then there is no braid theory (as formulated in
this paper). The braid groups of compact, connected
surfaces have been widely studied; (finite) presentations were
obtained in~\cite{Z1,Z2,Bi1,S}. As well as being interesting in their
own right, braid groups have played an important r\^ole in many
branches of mathematics, for example in topology, geometry, algebra
and dynamical systems, and notably in the study of knots and
links~\cite{BZ}, of the mapping class groups~\cite{Bi2,Bi3}, and
of configuration spaces~\cite{CG,FH}. The reader may
consult~\cite{Bi2,Han,MK,R} for some general references on the
theory of braid groups.

Let $M$ be a connected manifold of dimension~$2$ (or
\emph{surface}), perhaps with boundary. Further, we shall suppose that
$M$ is homeomorphic to a compact $2$-manifold with a finite (possibly
zero) number of points removed from its interior. We recall two
(equivalent) definitions of surface braid groups. The first is that
due to Fox. Let $F_n(M)$ denote the \emph{$n\up{th}$ configuration
space} of $M$, namely the set of all ordered $n$-tuples of distinct
points of $M$:
\begin{equation*}
F_n(M)=\set{(x_1,\ldots,x_n)}{\text{$x_i\in M$ and $x_i\neq x_j$ if $i\neq
j$}}. 
\end{equation*}
Since $F_n(M)$ is a subspace of the $n$-fold Cartesian product of $M$ with itself, the topology on $M$ induces a topology on $F_n(M)$. Then we define the \emph{$n$-string pure (or unpermuted) braid group $\pnm$} of $M$ to be:
\begin{equation*}
\pnm=\pi_1(F_n(M)).
\end{equation*}
There is a natural action of the symmetric group $S_n$ on $F_n(M)$ by
permutation of coordinates, and the resulting orbit space $F_n(M)/\sn$
shall be denoted by $D_n(M)$. The fundamental group $\pi_1(D_n(M))$ is
called the \emph{$n$-string (full) braid group} of $M$, and shall be
denoted by $B_n(M)$. Notice that the projection $F_n(M) \to D_n(M)$ is
a regular $n!$-fold covering map. It is well known that $B_n$ is isomorphic to $B_n(\dt)$ and $P_n\cong P_n(\dt)$, where $\dt$ is the closed $2$-disc.

The second definition of surface braid groups is geometric. Let
$\mathcal{P}=\brak{p_1, \dots, p_n}$ be a set of $n$ distinct points of $M$. A
\emph{geometric braid} of $M$ with basepoint $\mathcal{P}$ is a collection
$\beta=(\beta_1, \dots, \beta_n)$ of $n$ paths $\map{\beta}{[0,1]}[M]$ such
that:
\begin{enumerate}[(a)]
\item for all $i=1,\ldots,n$, $\beta_i(0)=p_i$ and $\beta_i(1)\in \mathcal{P}$.
\item for all $i,j=1,\ldots,n$ and $i\neq j$, and for all $t\in [0,1]$, $\beta_i(t)\neq \beta_j(t)$.
\end{enumerate}
Two geometric braids are said to be \emph{equivalent} if there exists a homotopy between them through geometric braids. The usual concatenation of paths induces a group operation on the set of equivalence classes of geometric braids. This group is isomorphic to $B_n(M)$, and does not depend on the choice of $\mathcal{P}$. The subgroup of \emph{pure} braids, satisfying additionally $\beta_i(1)=p_i$ for all $i=1,\ldots,n$, is isomorphic to $P_n(M)$. There is a natural surjective homomorphism $B_n(M)\to \sn$ which to a geometric braid $\beta$ associates the permutation $\pi$ defined by $\beta_i(1)= p_{\pi(i)}$. The kernel is precisely $P_n(M)$, and we thus obtain the following short exact sequence:
\begin{equation*}
1\to P_n(M)\to B_n(M) \to \sn\to 1.
\end{equation*}

\section{The Fadell-Neuwirth short exact sequence}

Let $m,n\in \N$ be positive integers such that $m>n$, and consider
the projection
\begin{align*}
p\colon\thinspace F_m(M) &\to F_n(M)\\
(x_1,\ldots,
x_n,\ldots ,x_m)&\mapsto  (x_1,\ldots,x_n).
\end{align*}
In~\cite{FaN}, Fadell and Neuwirth studied the map $p$, and showed
that it is a locally-trivial fibration. The fibre over a point
$(x_1,\ldots,x_n)$ of the base space is
$F_{m-n}(M\setminus\brak{x_1,\ldots,x_n})$ which may be considered to
be a subspace of the total space via the map 
\begin{equation*}
\map{i}{F_{m-n}(M\setminus\brak{x_1,\ldots,x_n})}[F_m(M)]
\end{equation*}
defined by
\begin{equation*}
i((y_1,\ldots,y_{m-n}))=(x_1,\ldots,x_n,y_1,\ldots,y_{m-n}). 
\end{equation*}
Then $p$ induces a group homomorphism $\map{p_{\ast}}{\pnm[m]}[\pnm]$,
which representing $\pnm[m]$ geometrically as a collection of $m$
strings, corresponds to forgetting the last $(m-n)$ strings.
\textbf{We adopt the convention throughout this paper, that unless
explicitly stated otherwise, all homomorphisms \( \pnm[m] \to \pnm \)
in the text will be this one}.

The fibration $\map{p}{F_m(M)}[F_n(M)]$ gives rise to a long exact sequence of homotopy groups of configuration spaces, from which we obtain the
\emph{Fadell-Neuwirth pure braid group short exact sequence}:
\begin{equation}\label{eq:split}
1 \to P_{m-n}(M\setminus\brak{x_1,\ldots,x_n}) \stackrel{i_{\ast}}{\to} \pnm[m] \stackrel{p_{\ast}}{\to}
\pnm \to 1,
\end{equation}
where $i_{\ast}$ is the group homomorphism induced by  $i$, and $n\geq
3$ if $M$ is the $2$-sphere $\St$~\cite{Fa,FvB}, $n\geq 2$ if $M$ is
the real projective plane $\mathbb{R}P^2$~\cite{vB}, and $n\geq 1$
otherwise~\cite{FaN} (in each case, the condition on $n$ implies that
$F_n(M)$ is an Eilenberg-MacLane space).  This short exact sequence
plays a central r\^ole in the study of surface braid groups. It was
used by~\cite{PR} to study mapping class groups, in the work
of~\cite{GMP} on Vassiliev invariants for braid groups, as well as to
obtain presentations for surface pure braid
groups~\cite{Bi1,S,GG1,GG4}.

An interesting question is that of whether the Fadell-Neuwirth short
exact sequence~\reqref{split} splits. If the above conditions on $n$
are satisfied then the existence of a section for $p_{\ast}$ is
equivalent to that of a geometric section for $p$
(cf.~\cite{GG3,GG4}). In~\cite{A2}, Artin showed that if $M$ is the
plane then~\reqref{split} splits for all $n\in\N$. This implies that
$P_n$ may be expressed as a repeated semi-direct product of free
groups, which enables one to solve the word problem in the pure and
full Artin braid groups. The splitting problem has been studied for
other surfaces besides the plane. Fadell and Neuwirth gave various
sufficient conditions for the existence of a geometric section for $p$
in the general case~\cite{FaN}. For the sphere, it was known that
there exists a section on the geometric level~\cite{FvB}.  If $M$ is
the $2$-torus then Birman exhibited an explicit algebraic section
for~\reqref{split} for $m=n+1$ and $n\geq 2$~\cite{Bi1}. However, for
compact orientable surfaces without boundary of genus $g\geq 2$, she
posed the question of whether the short exact sequence~\reqref{split}
splits. In~\cite{GG1}, we provided a complete answer to this question:
\begin{thm}[\cite{GG1}]\label{th:fnsplit}
If $M$ is a compact orientable surface without boundary of genus $g\geq 2$, the short exact sequence~\reqref{split} splits if and only if $n=1$.
\end{thm}

\section[A generalisation of the Fadell-Neuwirth sequence]{A generalisation of the Fadell-Neuwirth short exact sequence}

As we mentioned above, the Fadell-Neuwirth short exact sequence is a
very important tool in the study of pure surface braid groups, but
unfortunately it does not generalise directly to the corresponding
full braid groups. However, by considering intermediate coverings
between $F_n(M)$ and $D_n(M)$, it is possible to extend it to certain
subgroups of $B_n(M)$~\cite{GG2}. A special case of this construction may be formulated as follows. Let $m,n\in\N$, and let $D_{m,n}(M)$ denote the   quotient space of $F_{m+n}(M)$ by the action of the subgroup $\sn[m]\times \sn[n]$ of $\sn[m+n]$. Then we obtain a fibration $D_{m,n}(M)\to D_m(M)$, defined by forgetting the last $n$ coordinates. We set $B_{m,n}(M)=\pi_1(D_{m,n}(M))$, sometimes termed a \emph{`mixed'} braid group. As in the pure braid group case, we obtain a generalisation of the short exact sequence of Fadell and Neuwirth:
\begin{equation}\label{eq:fnsurface}
1\to B_n(M\setminus\brak{x_1,\ldots,x_m}) \to B_{m,n}(M)\stackrel{p_{\ast}}{\to} B_m(M)\to 1,
\end{equation}
where again we take $m\geq 3$ if $M=\St$, $m\geq 2$ if $M=\mathbb{R}P^2$ and $m\geq 1$ otherwise. Once more, unless explicitly stated, all homomorphisms $B_{m,n}(M)\to B_m(M)$ in the text will be this one.

\section{The braid groups of the sphere}

The braid groups of the sphere and the real projective plane are of particular interest, notably because they have non-trivial centre (which is also the case for the Artin braid groups), and torsion elements. The braid groups of the sphere were studied during the 1960's~\cite{Fa,FvB,vB,GVB}: let us recall briefly some of their properties.

If $\dt\subseteq \St$ is a topological disc, there is a group
homomorphism $\map {\iota}{B_n(\dt)}[B_n(\St)]$ induced by the
inclusion. If $\beta\in B_n(\dt)$ then we shall denote its image
$\iota(\beta)$ simply by $\beta$. It is well known that $B_n(\St)$ is
generated by $\sigma_1,\ldots,\sigma_{n-1}$ which are subject to the
following relations:
\begin{equation}\label{eq:presnbns}
\left.
\begin{aligned}
&\text{$\si{i}\si{j}=\si{j}\si{i}$ if $\lvert i-j\rvert\geq 2$ and $1\leq i,j\leq n-1$}\\
&\text{$\si{i}\si{i+1}\si{i}=\si{i+1}\si{i}\si{i+1}$ for all $1\leq i\leq n-2$, and}\quad\\
&\text{$\si 1\cdots \si {n-2}\si [2]{n-1}\si {n-2}\cdots\si 1=1$.}
\end{aligned}
\right\}
\end{equation}
In what follows, the third relation will be referred to as the
\emph{surface relation} of $B_n(\St)$. It follows from this
presentation and \req{presartin} that $B_n(\St)$ is a quotient of
$B_n$. The first three sphere braid groups are finite: $B_1(\St)$ is trivial, $B_2(\St)$ is cyclic of order~$2$, and $B_3(\St)$ is a $ZS$-metacyclic group (a group whose Sylow subgroups, commutator subgroup and commutator quotient group are all cyclic) of order~$12$.

If $n\geq 3$, the so-called `full twist' $\ft$ braid of $B_n(\St)$, defined by 
\begin{equation*}
\ft= (\sigma_1\cdots\sigma_{n-1})^n,
\end{equation*}
generates the centre $Z(B_n(\St))$ of $B_n(\St)$, and is a torsion element of order $2$. Using Seifert fibre space theory, Murasugi characterised the torsion elements of $B_n(\St)$: they are all conjugates of powers of the three elements $\sigma_1\cdots \sigma_{n-2} \sigma_{n-1}$, $\sigma_1\cdots \sigma_{n-2} \sigma_{n-1}^2$ and $\sigma_1\cdots \sigma_{n-3} \sigma_{n-2}^2$ which are respectively $n\up{th}$, $(n-1)\up{th}$ and $(n-2)\up{th}$ roots of $\ft$~\cite{M}.

In~\cite{GG4}, we studied the short exact sequence~\reqref{fnsurface} in the case $M=\St$ of the sphere:
\begin{equation}\label{eq:seqbnm}
1\to B_n(\St\setminus\brak{x_1,\ldots,x_m}) \to B_{m,n}(\St)\stackrel{p_{\ast}}{\to} B_m(\St)\to 1,
\end{equation}
and proved the following results:
\begin{thm}[\cite{GG4}]\mbox{}
\begin{enumerate}[(a)] 
\item\label{it:bmn4a} The short exact sequence
\begin{equation*}
1\to B_n(\St\setminus\brak{x_1,x_2,x_3}) \to B_{3,n}(\St)\stackrel{p_{\ast}}{\to} B_3(\St)\to 1
\end{equation*}
splits if and only if $n\equiv 0,2\bmod 3$.
\item\label{it:bmn4b} Let $m\geq 4$. If the homomorphism $\map {p_{\ast}}{B_{m,n}(\St)}[B_m(\St)]$ admits a section then there exist $\epsilon_1, \epsilon_2\in\brak{0,1}$ such that:
\begin{equation*}
n\equiv \epsilon_1 (m-1)(m-2)-\epsilon_2 m(m-2) \bmod{m(m-1)(m-2)}. 
\end{equation*}
\end{enumerate}
\end{thm}
An open question is whether the necessary condition in
part~(\ref{it:bmn4b}) is also sufficient. If $n\geq 4$ then $B_n(\St)$
is infinite, and it follows from the proof of part~(\ref{it:bmn4a})
that $B_n(\St)$ contains an isomorphic copy of the finite group
$B_3(\St)$ of order $12$ if and only if $n\not\equiv 1 \bmod 3$. We
have recently shown that $B_n(\St)$ contains an isomorphic copy of the quaternion group $\mathcal{Q}_8$ of order $8$ if and only if $n$ is even~\cite{GG5}. The realisation of finite subgroups in $B_n(\St)$ and $B_n(\R P^2)$ seems an interesting problem which we are pursuing.

\section{Braid group series and motivation for their study}

If $G$ is a group, then we recall that its \emph{lower central series}
$\brak{\Gamma_i(G)}_{i\in\N}$ is defined inductively by
$\Gamma_1(G)=G$, and $\Gamma_{i+1}(G)=[G, \Gamma_i(G)]$ for all
$i\in\N$, and its \emph{derived series}
$\brak{G^{(i)}}_{i\in\N\cup\brak{0}}$ is defined inductively by
$G^{(0)}=G$, and $G^{(i)}=[G^{(i-1)}, G^{(i-1)}]$ for all $i\in\N$.
One may check easily that $\Gamma_i(G)\supseteq \Gamma_{i+1}(G)$ and
$G^{(i-1)}\supseteq G^{(i)}$ for all $i\in\N$, and for all $j\in\N$,
$j>i$, $\Gamma_j(G)$ (resp.\ $G^{(j)}$) is a normal subgroup of
$\Gamma_i(G)$ (resp.\ $G^{(i)}$).  Notice that $\Gamma_2(G)=G^{(1)}$
is the \emph{commutator subgroup} of $G$. The \emph{Abelianisation} of
the group $G$, denoted by $\gpab$ is the quotient $G/\Gamma_2(G)$; the
\emph{Abelianisation} of an element $g\in G$ is its
$\Gamma_2(G)$-coset in $\gpab$. The group $G$ is said to be
\emph{perfect} if $G=G^{(1)}$, or equivalently if $\gpab=\brak{1}$.
Following P.~Hall, for any group-theoretic property $\mathcal{P}$, a
group $G$ is said to be \emph{residually $\mathcal{P}$} if for any
(non-trivial) element $x \in G$, there exists a group $H$ with the
property  $\mathcal{P}$ and a surjective homomorphism
$\map{\phi}{G}[H]$ such that  $\phi(x) \not=1$. It is well known that
a group $G$ is \emph{residually nilpotent} (respectively \emph{residually soluble}) if and only if $\bigcap_{i \geq 1}\, \Gamma_i(G)=\{ 1\}$ (respectively $\bigcap_{i \geq 0}\, G^{(i)}=\{ 1\}$). If $g,h\in G$ then $[g,h]=gh g^{-1}h^{-1}$ will denote their commutator, and we shall use the symbol $g\comm h$ to mean that $g$ and $h$ commute.

Our main aim in this monograph is to study the lower central and
derived series of the braid groups of the sphere and the punctured
sphere. This was motivated in part by the study of the problem of the
existence of a section for the short exact sequences~\reqref{split}
and~\reqref{fnsurface}. To obtain a positive answer, it suffices of
course to exhibit an explicit section (although this may be easier
said than done!). However, and in spite of the fact that we possess
presentations of surface braid groups, in general it is very difficult
to prove directly that such an extension does not split. One of the
main methods that we used to prove the non-splitting of~\reqref{split}
for $n\geq 2$ and of~\reqref{seqbnm} for $m\geq 4$ was based on the
following observation: let $1\to K\to G\to Q\to 1$ be a split
extension of groups, where $K$ is a normal subgroup of $G$, and let
$H$ be a normal subgroup of $G$ contained in $K$. Then the extension
$1\to K/H\to G/H\to Q\to 1$ is also split. The condition on $H$ is
satisfied for example if $H$ is an element of either the lower central
series or the derived series of $K$. In~\cite{GG1}, considering the
extension~\reqref{split} with $n\geq 3$, we showed that it was
sufficient to take  $H=\Gamma_2(K)$ to prove the non-splitting of the
quotiented extension, and hence that of the full extension. In this
case, the kernel $K/\Gamma_2(K)$ is Abelian, which simplifies somewhat
the calculations in $G/H$. This was also the case in~\cite{GG4} for
the extension~\reqref{seqbnm} with $m\geq 4$. However, for the
extension~\reqref{split} with $n=2$, it was necessary to go a stage
further in the lower central series, and take $H=\Gamma_3(K)$. From
the point of view of the splitting problem, it is thus helpful to know
the lower central and derived series of the braid groups occurring in
these group extensions. But these series are of course interesting in
their own right, and help us to understand better the structure of
surface braid groups. 

Let us remark that braid groups of the punctured disc were studied in~\cite{Lam} in relation with the study of knots in handlebodies, and were used by Bigelow to understand the Lawrence-Krammer representation in his proof of the linearity of the Artin braid groups~\cite{Big}. Furthermore, during our study of the braid groups of the $2$- and $3$-punctured sphere, we will also come across some of the \emph{Artin and affine Artin groups} (also known as \emph{generalised braid groups}), notably those of types $B$ and $\widetilde{C}$~\cite{Bri,T}.

The lower central series of groups and their successive quotients
$\Gamma_i/\Gamma_{i+1}$ are isomorphism invariants, and have been
widely studied using commutator calculus, in particular for free
groups of finite rank~\cite{H,MKS}. Falk and Randell, and
independently Kohno investigated the lower central series of the pure
braid group $P_n$, and were able to conclude that $P_n$ is residually
nilpotent~\cite{FR1,Ko}. Falk and Randell also studied the lower
central series of generalised pure braid groups~\cite{FR2,FR3}. 

Using the Reidemeister-Schreier rewriting process, Gorin and Lin
obtained a presentation of the commutator subgroup of $B_n$ for $n\geq
3$~\cite{GL} (see \reth{gorinlin}). For $n\geq 5$, they were able to
infer that $(B_n)^{(1)}=(B_n)^{(2)}$, and so $(B_n)^{(1)}$ is perfect.
From this it follows that $\Gamma_2(B_n)=\Gamma_3(B_n)$, hence $B_n$
is not residually nilpotent. If $n=3$ then they showed that
$(B_3)^{(1)}$ is a free group of rank~$2$, while if $n=4$, they proved
that $(B_4)^{(1)}$ is a semi-direct product of two free groups of
rank~$2$. By considering the action, one may see that
$(B_4)^{(1)}\supsetneqq (B_4)^{(2)}$. The work of Gorin and Lin on
these series was motivated by the study of almost periodic solutions
of algebraic equations with almost periodic coefficients.

\section{Statement of the main results}

\rechap{sphere} is devoted to determining the lower central series
of the braid groups of the sphere.  In \reth{lcsbn}, we show that for
all $n\geq 2$, the lower central series is constant from the
commutator subgroup onwards. As in the case of the disc, the case
$n=4$ is particularly interesting: $\Gamma_2(B_4(\St))$ is a
semi-direct product of the quaternion group of order~$8$ by the free group of rank~$2$. Here is the main theorem of \rechap{sphere}:
\begin{thm}\label{th:lcsbn} 
For all $n\geq 2$, the lower central series of $B_n(\St)$ is constant
from the commutator subgroup onwards: $\Gamma_m(B_n(\St))=
\Gamma_2(B_n(\St))$ for all $m\geq 2$. The subgroup $\Gamma_2(B_n(\St))$
is as follows:   
\begin{enumerate}
\item If $n=1,2$ then $\Gamma_2(B_n(\St))=\brak{1}$.
\item If $n=3$ then $\Gamma_2(B_n(\St)) \cong\Z_3$. Thus $B_3(\St)\cong \Z_3\rtimes \Z_4$, the action being the non-trivial one.
\item If $n=4$ then $\Gamma_2(B_4(\St))$ admits a presentation of the following form:
\begin{enumerate}
\item[\underline{\textbf{generators:}}] $g_1,g_2,g_3$, where in terms of the
usual generators of $B_4(\St)$, $g_1=\sigma_1^2\sigma_2 \sigma_1^{-3}$,
$g_2=\sigma_1^3\sigma_2 \sigma_1^{-4}$ and $g_3=\sigma_3 \sigma_1^{-1}$.
\item[\underline{\textbf{relations:}}]
\begin{gather*}
g_3^4=1 \\
\text{$[g_3^2,g_i]=1$ for $i=1,2$}\\
[g_3,g_2g_1]=1\\
g_1^{-1}g_3^{-1}g_1=g_2g_3g_2^{-1}\\ g_1^{-1}g_3^{-1}g_1=g_1g_3g_1^{-1}g_3. \
\end{gather*}
\end{enumerate}
Furthermore, 
\begin{equation*}
\Gamma_2(B_4(\St))\cong \mathcal{Q}_8\rtimes \F[2](a,b),
\end{equation*}
where $\mathcal{Q}_8= \setang{x,y}{x^2=y^2,\; xyx^{-1}=y^{-1}}$ is the quaternion group of order~$8$, and $\F[2](a,b)$ is the free group of rank~$2$ on two generators $a$ and $b$. The action is given by:
\begin{align*}
&\phi(a)(x)=y &  &\phi(a)(y)=xy\\
&\phi(b)(x)=yx && \phi(b)(y)=x.
\end{align*}
\item In the cases $n=5$ and $n\geq 6$, a presentation for $\Gamma_2(B_n(\St))$ is given in \rechap{prescom}, by Propositions~\ref{prop:g2b5} and~\ref{prop:g2b6} respectively.
\end{enumerate}
\end{thm}
The lower central series of $B_n(\St)$ is thus completely determined. In particular, if $n\geq 3$ then $B_n(\St)$ is not residually nilpotent.

\medskip

In \rechap{dssphere}, we study the derived series of $B_n(\St)$. As in
the case of the disc, $(B_n(\St))^{(1)}$ is perfect if $n\geq 5$, in
other words, the derived series of $B_n(\St)$ is constant from
$(B_n(\St))^{(1)}$ onwards. The cases $n=1,2,3$ are straightforward,
and the groups $B_n(\St)$ are finite and soluble. In the case $n=4$,
we make use of the semi-direct product decomposition of
$(B_4(\St))^{(1)}$ obtained in \reth{lcsbn}. \repr{gammasemi}
describes the structure of the commutator subgroup of a general
semi-direct product, and shall be applied frequently throughout this
monograph. This will enable us to show that from $(B_4(\St))^{(4)}$
onwards, the derived series of $B_4(\St)$ coincides with that of the
free group of rank~$2$. We also determine some of the derived series
quotients of $B_4(\St)$:
\begin{thm}\label{th:dsbn}
The derived series of $B_n(\St)$ is as follows.
\begin{enumerate}[(a)]
\item If $n=1,2$ then $(B_n(\St))^{(1)}=\brak{1}$.
\item If $n=3$ then $(B_3(\St))^{(1)}\cong\Z_3$ and $(B_3(\St))^{(2)}=\brak{1}$.
\item Suppose that $n=4$. Then:
\begin{enumerate}[(i)]
\item $(B_4(\St))^{(1)}=\Gamma_2(B_4(\St))$ is given by part~(\ref{it:lcs4}) of \reth{lcsbn}; it is isomorphic to the semi-direct product $\mathcal{Q}_8\rtimes \F[2]$. Further, $B_4(\St)/(B_4(\St))^{(1)}$ is isomorphic to $\Z_6$.
\item  $(B_4(\St))^{(2)}$ is isomorphic to the semi-direct product $\mathcal{Q}_8\rtimes \F[2]^{(1)}$, where $(\F[2])^{(1)}$ is the commutator subgroup of the free group $\F[2]=\F[2](a,b)$ of rank~$2$ on two generators $a,b$. The action of $(\F[2])^{(1)}$ on $\mathcal{Q}_8$ is the restriction of the action of $\F[2](a,b)$ given in part~(\ref{it:lcs4}) of \reth{lcsbn}. Further, 
\begin{equation*}
\text{$(B_4(\St))^{(1)}/(B_4(\St))^{(2)}\cong \Z^2$, and $B_4(\St)/(B_4(\St))^{(2)} \cong \Z^2\rtimes \Z_6$,}
\end{equation*}
where the action of the generator $\overline{\sigma}$ of $\Z_6$ on $\Z^2$ is given by left multiplication by the matrix $\left( \begin{smallmatrix}
0 & 1\\
-1 & 1
\end{smallmatrix}\right)$.
\item $(B_4(\St))^{(3)}$ is a subgroup of $P_4(\St)$ isomorphic to the direct product $\Z_2\times (\F[2])^{(2)}$. Further,
\begin{equation*}
\text{$(B_4(\St))^{(2)}/(B_4(\St))^{(3)}\cong (\Z_2 \times \Z_2)\times (\F[2])^{(1)}/(\F[2])^{(2)}$.}
\end{equation*}
\item $(B_4(\St))^{(m)}\cong (\F[2])^{(m-1)}$ for all $m\geq 4$. Further, 
\begin{equation*}
(B_4(\St))^{(3)}/(B_4(\St))^{(4)}\cong \Z_2\times (\F[2])^{(2)}/(\F[2])^{(3)},
\end{equation*}
and for $m\geq 4$, 
\begin{equation*}
(B_4(\St))^{(m)}/(B_4(\St))^{(m+1)}\cong (\F[2])^{(m-1)}/(\F[2])^{(m)}.
\end{equation*}
\end{enumerate}
\item If $n\geq 5$ then $(B_n(\St))^{(2)}=(B_n(\St))^{(1)}$, so $(B_n(\St))^{(1)}$ is perfect. A presentation of $(B_n(\St))^{(1)}$ is given in Propositions~\ref{prop:g2b5} and~\ref{prop:g2b6}.
\end{enumerate}
\end{thm}
In particular, the derived series of $B_n(\St)$ is thus completely determined (up to knowing the derived series of the free group $\F[2]$ of rank~$2$, see \rerem{freederive}).

\medskip

\rechap{lcsbmn} deals with the lower central and derived series of braid groups of the punctured sphere $B_m(\St\setminus\brak{x_1,\ldots, x_n})$, $n\geq 1$, and is divided into eight sections, according to the respective values of $m$ and $n$. In \repr{presbetanm} (\resec{presbmn}), we recall a presentation of these groups obtained in~\cite{GG4}. In \resec{gorinlin}, we consider the case $n=1$, and show that $B_m(\St\setminus\brak{x_1})$ is isomorphic to $B_n(\dt)$ (\repr{iso}). In \repr{gorinlin}, we study the series of $B_n(\dt)$ in further detail, thus extending the results of Gorin and Lin:
\begin{prop}\label{prop:gorinlin}
Let $m\geq 1$. Then:
\begin{enumerate}[(a)]
\item For all $s\geq 3$, $\Gamma_s(B_m(\dt))= \Gamma_2(B_m(\dt))$.
\item If $m=1,2$ then $(B_m(\dt))^{(s)}=\brak{1}$ for all $s\geq 1$.
\item\label{it:b3deriv2}  If $m=3$ then the derived series of $(B_3(\dt))^{(1)}$ is that of the free group $\F[2](u,v)$ on two generators $u$ and $v$, where $u=\sigma_2\sigma_1^{-1}$ and $v=\sigma_1 u\sigma_1^{-1}= \sigma_1 \sigma_2\sigma_1^{-2}$.
Further,
\begin{equation*}
B_3(\dt)\left/(B_3(\dt))^{(2)}\right. \cong \Z^2\rtimes \Z,
\end{equation*}
where $\Z^2$ is the free Abelian group generated by the respective Abelianisations $\overline{u}$ and $\overline{v}$ of $u$ and $v$, and the action is given by $\sigma \cdot \overline{u}=\overline{v}$ and $\sigma \cdot \overline{v}=-\overline{u}+\overline{v}$, where $\sigma$ is a generator of $\Z$.
\item If $m=4$ then 
\begin{gather*}
\text{$(B_4(\dt))^{(1)}\left/(B_4(\dt))^{(2)} \right. \cong \Z^2$, and}\\
(B_4(\dt))^{(2)}\cong \F[2](a,b)\rtimes \Gamma_2(\F[2](u,v)),
\end{gather*}
where $a=\sigma_3\sigma_1^{-1}$ and $b= uau^{-1}= \sigma_2\sigma_3
\sigma_1^{-1} \sigma_2^{-1}$. 
\end{enumerate}
\end{prop}
Hence the lower central series (respectively derived series) of $B_m(\dt)$ is completely determined for all $m\geq 1$ (respectively for all $m\neq 4$; for the case $m=3$, this is again up to knowing the derived series of $\F[2]$). 

\medskip

In the difficult case of the derived series of $B_4(\dt)$, we then go on to describe some of the higher order terms and the successive derived series quotients:
\begin{prop}\label{prop:b4disc2ab}
\begin{equation*}
\left(B_4(\dt)\right)^{(2)}\left/\left(B_4(\dt)\right)^{(3)}\right. \cong  \Z_2\times \Z_2\times
\gpab[{(\Gamma_2(\F[2](u,v)))}].
\end{equation*}
\end{prop}

\begin{prop}\label{prop:b4deriv}
$(B_4(\dt))^{(3)}\cong \F[5](z_1,\ldots,z_5) \rtimes (\F[2](u,v))^{(2)}$.
\end{prop}
The action for this semi-direct product will be described by equations~\reqref{acuxi} and~\reqref{acvxi}. From this, we may obtain the Abelianisation of $(B_4(\dt))^{(3)}$:
\begin{prop}\label{prop:b4derivab}
\begin{align*}
\gpab[{\left((B_4(\dt))^{(3)}\right)}]  &=(B_4(\dt))^{(3)}\left/(B_4(\dt))^{(4)}\right.\\
& \cong \Z^3 \times \Z_{18}\times \Z_{18} \times (\F[2](u,v))^{(2)}\left/(\F[2](u,v))^{(3)}.\right.
\end{align*}
\end{prop}
This result suggests that the derived series of $B_4(\dt)$ is highly non trivial. In principle, using the semi-direct product structure of $(B_4(\dt))^{(3)}$ and \repr{gammasemi}, it is possible to discover further terms of the derived series, but in practice, the calculations become very hard. The main results of \resec{gorinlin} are summed up in Table~\ref{table:bmdisc}.

\medskip

\begin{table}
\renewcommand{\arraystretch}{1.22}
\begin{tabular}{|c||c|c|c|}\hline
values  &  &  & \\ 
of $m$ & \raisebox{1.5ex}[0pt]{series/group} & \raisebox{1.5ex}[0pt]{result} &\raisebox{1.5ex}[0pt]{reference} \\
\hline\hline
$\forall m\geq 1$ & lower central &$ \Gamma_3(G)= \Gamma_2(G)$ & \\ \cline{1-3}
$\forall m\geq 5$ & derived & $G^{(2)}=G^{(1)}$ & \\ \cline{1-3}
 & derived & $G^{(i)}=(\F[2])^{(i-1)}$, $\forall i\geq 1$ & \raisebox{1.5ex}[0pt]{\cite{GL} (see}  \\ \cline{2-3}
 & $\Gamma_2(G)$ & $\F[2]$ & \raisebox{1.5ex}[0pt]{\reth{gorinlin})}  \\ \cline{2-3}
\raisebox{1.5ex}[0pt]{$m=3$}  &  $G^{(1)}/G^{(2)}$ & $\Z^2$ & \\ \cline{2-4}
   & $G/G^{(2)}$ & $\Z^2\rtimes \Z $ & \repr{gorinlin}\\ \hline      
   & $\Gamma_2(G)$ & $\F[2]\rtimes \F[2]$& \cite{GL} (see \\ \cline{2-3}   
  & $G^{(1)}/G^{(2)}$ & $\Z^2$ & \reth{gorinlin})\\ \cline{2-4}  
   & $G^{(2)}$ & $\F[2]\rtimes (\F[2])^{(1)}$ & \repr{gorinlin}\\ \cline{2-4}   
\raisebox{1.5ex}[0pt]{$m=4$}   & $G^{(2)}/G^{(3)}$ & $ \Z_2\times \Z_2\times \gpab[{((\F[2])^{(1)})}]$ & \repr{b4disc2ab}\\ \cline{2-4}
  & $G^{(3)}$ & $\F[5]\rtimes (\F[2])^{(2)}$ & \repr{b4deriv}\\ \cline{2-4}  
   & $G^{(3)}/G^{(4)}$ & $\Z^3\times \Z_{18}\times \Z_{18} \times \gpab[{((\F[2])^{(2)})}]$ & \repr{b4derivab}\\ \hline   
\end{tabular}\medskip
\caption{\label{table:bmdisc}Summary of results of \resec{gorinlin}, \rechap{lcsbmn} concerning the lower central and derived series of $G=B_m(\dt)$. For the semi-direct product actions, one should consult the corresponding reference.}
\end{table}

In \resec{bmn1n}, we comment briefly on the case $m=1$ which is that of a free group of rank $n-1$.
From \resec{m3n2} of \rechap{lcsbmn} onwards, we suppose that $n\geq 2$. If $m\geq 3$ (resp.\ $m\geq 5$) the lower central series (resp.\ the derived series) of $B_m(\St\setminus\brak{x_1,\ldots, x_n})$ is constant from the commutator subgroup onwards. Once more, for the derived series, $m=4$ represents a challenging case. Nevertheless, we are able to determine some of the derived series quotients. The main theorem of \resec{m3n2} is as follows:
\begin{thm}\label{th:lcdsbmsn}
Let $n\geq 2$. Then:
\begin{enumerate}[(a)]
\item If $m\geq 3$ then
\begin{equation*}
\Gamma_3(B_m(\St\setminus\brak{x_1,\ldots, x_n}))= \Gamma_2(B_m(\St\setminus\brak{x_1,\ldots, x_n})).
\end{equation*}
\item If $m\geq 5$ then
\begin{equation*}
(B_m(\St\setminus\brak{x_1,\ldots, x_n}))^{(2)}= (B_m(\St\setminus\brak{x_1,\ldots, x_n}))^{(1)}.
\end{equation*}
\item If $m=4$ then 
\begin{equation*}
B_4(\St\setminus\brak{x_1,\ldots, x_n})\left/(B_4(\St\setminus\brak{x_1,\ldots, x_n}))^{(2)}\right. \cong \left( \Z^2\rtimes \Z\right) \times \Z^{n-1}
\end{equation*}
where the semi-direct product structure is that of part~(\ref{it:b3deriv2}) of \repr{gorinlin},
and
\begin{equation*}
(B_4(\St\setminus\brak{x_1,\ldots, x_n}))^{(1)}\left/(B_4(\St\setminus\brak{x_1,\ldots, x_n}))^{(2)}\right. \cong \Z^2.
\end{equation*}
Alternatively, 
\begin{equation*}
B_4(\St\setminus\brak{x_1,\ldots, x_n})/(B_4(\St\setminus\brak{x_1,\ldots, x_n}))^{(2)}\cong \Z^2 \rtimes \Z^n,
\end{equation*}
where $\Z^2\cong (B_4(\St\setminus\brak{x_1,\ldots, x_n}))^{(1)}/(B_4(\St\setminus\brak{x_1,\ldots, x_n}))^{(2)}$ is the free Abelian group with basis $\brak{\overline{u}, \overline{v}}$, $\Z^n\cong \gpab[B_4(\St\setminus\brak{x_1,\ldots, x_n})]$ has basis $\brak{\sigma,\rho_1, \ldots,\rho_{n-1}}$, and the action is given by 
\begin{align*}
\sigma\cdot \overline{u} &= \overline{v} & \sigma\cdot \overline{v} &= -\overline{u}+ \overline{v}\\
\rho_i\cdot \overline{u} &= \overline{u} & \rho_i\cdot \overline{v} &= \overline{v}
\end{align*}
for all $1\leq i\leq n-1$.
\end{enumerate}
\end{thm}

\begin{table}
\renewcommand{\arraystretch}{1.22}
\begin{tabular}{|c||c|c|c|}\hline
values  &  &  & \\ 
of $m$ & \raisebox{1.5ex}[0pt]{series/group} & \raisebox{1.5ex}[0pt]{result} &\raisebox{1.5ex}[0pt]{reference} \\
\hline\hline
$\forall m\geq 3$ & lower central &$ \Gamma_3= \Gamma_2$ & \\ \cline{1-3}
$\forall m\geq 5$ & derived & $G^{(2)}=G^{(1)}$ & \\ \cline{1-3}
 &  & $ (\Z^2\rtimes \Z) \times \Z^{n-1}$ & \reth{lcdsbmsn} \\ \cline{3-3}
$m=4$ & \raisebox{1.5ex}[0pt]{$G/G^{(2)}$} & $\Z^2 \rtimes \Z^n$ &  \\ \cline{2-3}
 
 & $G^{(1)}/G^{(2)}$ & $\Z^2$ &\\ \hline
\end{tabular}\medskip
\caption{\label{table:bngeq2sphere}Summary of results of \resec{m3n2}, \rechap{lcsbmn} concerning the lower central and derived series of $G=B_m(\St\setminus\brak{x_1,\ldots, x_n})$, $m\geq 3$, $n\geq 2$. }
\end{table}

So if $n\geq 2$, the lower central and derived series of the braid
group $B_m(\St\setminus\brak{x_1,\ldots, x_n})$ are completely
determined, with the exception of a small number of values of $m$: for
the lower central series, they consist of just $m=2$, and for the
derived series, $m=2,3$ and $4$.

\medskip 

The case $m\geq 2$ and $n=2$ is considered in Sections~\ref{sec:mgeq2},~\ref{sec:lcgsb22} and~\ref{sec:affineatil}. Applying the results of \repr{iso}, one may see that $B_m(\St \setminus \brak{x_1,x_2})$ is isomorphic to the $m$-string braid group $B_m(\mathbb{A})$ of the annulus $\mathbb{A}=[0,1] \times \mathbb{S}^1$, and is thus an Artin group of type~$B_m$. In \repr{aug}, \resec{mgeq2}, we prove the following general result concerning the structure of $\Gamma_2(B_m(\St \setminus \brak{x_1,x_2}))$:
\begin{prop}\label{prop:aug}
Let $m\geq 2$. Then:
\begin{enumerate}[(a)]
\item $B_m(\St\setminus\brak{x_1,x_2})\cong \F[m] \rtimes B_m(\dt)$, where the action $\phi$ is given by the Artin representation of $B_m(\dt)$ as a subgroup of $\aut{\F[m]}$ (see \req{sigij}).
\item $\Gamma_2(B_m(\St\setminus\brak{x_1,x_2})) \cong \ker{\rho} \rtimes \Gamma_2(B_m(\dt))$, where
\begin{equation*}
\map{\rho}{\F[m](A_{2,3}, \ldots, A_{2,m+2})}[\Z]
\end{equation*}
is the augmentation homomorphism, and the action is that induced by $\phi$ (the generators $A_{i,j}$ are described in \repr{presbetanm}).
\end{enumerate}
\end{prop}
The semi-direct product structure allows us to determine some derived series quotients:
\begin{prop}\label{prop:b3z4}
\begin{equation*}
\left(B_3\left(\St\setminus\brak{x_1,x_2}\right)\right)^{(1)} \left/ \left(B_3\left(\St\setminus\brak{x_1,x_2}\right)\right)^{(2)}\right. \cong \Z^4.
\end{equation*}
\end{prop}

\begin{prop}\label{prop:z4z2}
\begin{equation*}
B_3\left(\St\setminus\brak{x_1,x_2}\right)\left/ \left(B_3\left(\St\setminus\brak{x_1,x_2}\right)\right)^{(2)}\right. \cong \Z^4 \rtimes \Z^2,
\end{equation*}
where $\Z^4$ has a basis $\brak{\widetilde{\alpha_0}, \widetilde{\beta_0}, \widetilde{u}, \widetilde{v}}$, $\Z^2$ has a basis $\brak{\sigma,\rho_1}$, and the action is given by:
\begin{align*}
\sigma&\cdot \widetilde{u}=\widetilde{v} & \sigma&\cdot \widetilde{v}=-\widetilde{u}+ \widetilde{v}\\
\sigma&\cdot \widetilde{\alpha_0}=\widetilde{\beta_0} &\sigma&\cdot \widetilde{\beta_0}= \widetilde{\beta_0} -\widetilde{\alpha_0}\\
\rho_1&\cdot \widetilde{\alpha_0}=\widetilde{\alpha_0} &  \rho_1&\cdot \widetilde{\beta_0}=\widetilde{\beta_0}\\
\rho_1&\cdot \widetilde{u}=  -\widetilde{\alpha_0} -\widetilde{u} +\widetilde{v}&
\rho_1&\cdot \widetilde{v}=-\widetilde{\beta_0}-\widetilde{u}.
\end{align*}
\end{prop}

We then give an alternative proof of \repr{b3z4}, showing along the
way that the commutator subgroup of 
$B_3(\St\setminus\brak{x_1,x_2})$ is the semi-direct product of a
given infinite rank subgroup of a free group of rank~$5$ by a free
group of rank~$2$ (see \repr{f5f2}). 

\medskip

In \resec{lcgsb22}, we study the lower central series of $B_2(\St\setminus\brak{x_1,x_2})$ (which is one of the outstanding cases not covered by \reth{lcdsbmsn}). Using an exact sequence due to Stallings (see \req{stallings}), we prove the following:
\begin{cor}\label{cor:b2b2}
$\Gamma_2(B_2(\St\setminus\brak{x_1,x_2})) \cong \Gamma_2(\F[2](a,b)) \rtimes \Z$, where the action of $\Z$ on $\Gamma_2(\F[2](a,b))$ is given by conjugation by $b^{-1}a$.
\end{cor}
The group $B_2(\St\setminus\brak{x_1,x_2})$ is particularly fascinating, not least because it may be interpreted in many different ways: as the $2$-string braid group $B_2(\mathbb{A})$ of the annulus  (and so as the Artin group of type~$B_2$), and as the Baumslag-Solitar group $\operatorname{BS}(2,2)$, for example (see \rerems{b2s2}). It is also a one-relator group with non-trivial (infinite cyclic) centre, which applying results of Kim and McCarron~\cite{KMc,McC} implies that:
\begin{prop}\label{prop:b2b2resid}
$B_2(\St\setminus\brak{x_1,x_2})$ is residually nilpotent and residually a finite $2$-group.
\end{prop}
Further, using the fact that the quotient of $B_2(\St\setminus\brak{x_1,x_2})$ by its centre is isomorphic to the free product $\Z_2\ast \Z$, we prove that apart from the first term, the lower central series of these two groups coincide, and applying results of Gaglione and Labute~\cite{Ga,La} which describe  the lower central series of certain free products of cyclic groups, we are able to determine completely the lower central series (in terms of that of $\Z_2\ast \Z$), as well as the successive lower central series quotients of $B_2(\St\setminus\brak{x_1,x_2})$ in an explicit manner: 
\begin{thm}\label{th:lcsb2}
For all $i\geq 2$, $\Gamma_i(B_2(\St\setminus\brak{x_1,x_2})) \cong \Gamma_i(\Z_2\ast \Z)$, and:
\begin{align*}
\Gamma_i(B_2(\St\setminus\!\brak{x_1,x_2}))/ \Gamma_{i+1}(B_2(\St\setminus\!\brak{x_1,x_2})) &\cong \Gamma_i(\Z_2\ast \Z)/ \Gamma_{i+1}(\Z_2\ast \Z)\\
& \cong \underbrace{\Z_2\oplus \cdots \oplus \Z_2}_{\text{$R_i$ times}},
\end{align*}
where
\begin{equation*}
R_i=\sum_{j=0}^{i-2}\; \left( \sum_{\substack{k\mid i-j\\ k>1}}\; \mu\left( \frac{i-j}{k}\right) \frac{k\alpha_k}{i-j} \right),
\end{equation*}
$\mu$ is the M\"obius function, and
\begin{equation*}
\alpha_k=\frac{1}{k} \left( \tr{\begin{pmatrix}
0 & -1\\
-1 & 1
\end{pmatrix}^k}-1\right).
\end{equation*}
\end{thm}
From this, we may see (\reco{dsb22}) that apart from the first term, the derived series of $B_2(\St\setminus\brak{x_1,x_2})$ is that of $\pi(\F[2])$, where $\map{\pi}{\F[2]=\Z\ast \Z}[\Z_2 \ast \Z]$ is the homomorphism obtained by taking the first factor modulo~$2$.

\begin{table}
\renewcommand{\arraystretch}{1.22}
\begin{tabular}{|c||c|c|c|}\hline
values  &  &  & \\ 
of $m$ & \raisebox{1.5ex}[0pt]{series/group} & \raisebox{1.5ex}[0pt]{result} &\raisebox{1.5ex}[0pt]{reference} \\
\hline\hline
$\forall m\geq 2$ & $\Gamma_2(G)$ & $(\F[m])^{(1)}\rtimes \Gamma_2(B_m(\dt))$ & \repr{aug}\\ \cline{1-4}
 & lower central & $\Gamma_i= \Gamma_i(\Z_2\ast\Z)$, $i\geq 2$ & \\ \cline{2-3}
 & lower central & $ \Gamma_i(\Z_2\ast\Z)/\Gamma_{i+1}(\Z_2\ast\Z)$ &  \\ \cline{3-3} 
 & quotients &  & \raisebox{1.5ex}[0pt]{\reth{lcsb2}} \\ 
  & $\Gamma_i(G)/\Gamma_{i+1}(G)$ & \raisebox{1.5ex}[0pt]{$\oplus_{j=1}^{R_i}\; \Z_2$} &\\
 \cline{2-4}
\raisebox{1.5ex}[0pt]{$m=2$} & $\Gamma_2(G)$ & $(\F[2])^{(1)}\rtimes \Z$ & \reco{b2b2}\\ \cline{2-4}
& $\Gamma_2(G)$ & $\F[\infty]$ & \reco{gam2b2s2}\\ \cline{2-4}
 & $\Gamma_2(G)/\Gamma_3(G)$ & $\Z_2$ & \\ \cline{2-3}
& $\Gamma_3(G)$ & $\F[\infty]$ & \raisebox{1.5ex}[0pt]{\repr{gam3b2}}\\ \cline{1-4}
& $\Gamma_2(G)$ & $\F[\infty] \rtimes \F[2]$ & \repr{f5f2}\\ \cline{2-4}
$m=3$ & $G^{(1)}/G^{(2)}$ & $\Z^4$ & \repr{b3z4}\\ \cline{2-4}
& $G/G^{(2)}$ & $\Z^4\rtimes \Z^2$ & \repr{z4z2}\\  \hline   
$m=4$ & $G^{(1)}/G^{(2)}$ & $\Z^2$ & \\ \cline{1-3}
$m\geq 5$ & $G^{(1)}/G^{(2)}$ & $\Z$ & \raisebox{1.5ex}[0pt]{\reco{dsannu}}\\ \hline
\end{tabular}\medskip
\caption{\label{table:b2n2sphere}Summary of results of Sections~\ref{sec:mgeq2},~\ref{sec:lcgsb22} and~\ref{sec:affineatil} of \rechap{lcsbmn} concerning the lower central and derived series of $G=B_m(\St\setminus\brak{x_1, x_2})$, $m\geq 2$. In each case, $\F[\infty]$ is a given free group of countable infinite rank.}
\end{table}

\medskip

In \resec{affineatil}, we consider the more general case of the $m$-string braid group $B_m(\St\setminus\brak{x_1, x_2}))$, $m\geq 3$, which we know to be isomorphic to the $m$-string braid group $B_m(\mathbb{A})$ of the annulus. With this interpretation, Kent and Peifer gave a nice presentation of this group (\repr{kent}) from which they were able to conclude that $B_m(\mathbb{A})$ is a semi-direct product of the affine Artin group $\widetilde{A}_{m-1}$ by $\Z$ (\reco{kent})~\cite{KP}. Applying \repr{gammasemi} once more,  we obtain in \repr{gamma2annu} a presentation of $\Gamma_2(B_m(\St\setminus\brak{x_1, x_2}))$ (which as we shall see, is isomorphic to $\Gamma_2(\widetilde{A}_{m-1})$), from which we may deduce:
\begin{cor}\label{cor:dsannu}
Let $m\geq 3$. Then 
\begin{equation*}
\left(B_m \left(\St\setminus\brak{x_1,x_2}\right)\right)^{(1)}\left/ \left(B_m\left(\St\setminus\brak{x_1,x_2}\right)\right)^{(2)}\right. \cong
\begin{cases}
\Z^4 & \text{if $m=3$}\\
\Z^2 & \text{if $m=4$}\\
\Z & \text{if $m\geq 5$.}
\end{cases}
\end{equation*}
\end{cor}
The main results of Sections~\ref{sec:mgeq2},~\ref{sec:lcgsb22} and~\ref{sec:affineatil} of \rechap{lcsbmn} are summed up in Table~\ref{table:b2n2sphere}.

\medskip

In \resec{bmn3} of \rechap{lcsbmn}, we consider $B_m(\St\setminus\brak{x_1,x_2,x_3})$, $m\geq 2$, which is also one of the outstanding cases for the derived series not covered by \reth{lcdsbmsn}. This group is isomorphic to the affine Artin group of type $\widetilde{C}_m$ for which little seems to be known~\cite{All}. Despite the existence of nice presentations for this group~\cite{BG}, we were not able to describe satisfactorily the commutator subgroup even for $m=2$. We obtain however some partial results, notably in \repr{b2n3} the fact that the successive lower central series quotients of $B_2(\St\setminus\brak{x_1,x_2,x_3})$ are finite direct sums of $\Z_2$, which generalises part of \reth{lcsb2}, as well as  for all $i\geq 1$ and $m\geq 2$, $(B_m(\St\setminus\brak{x_1,x_2,x_3}))^{(i)}$ is a semi-direct product of some group $K_i$ by $(B_m(\dt))^{(i)}$ (\repr{bmn3i}).

\medskip

Finally in \rechap{prescom}, we give presentations of the commutator subgroups $\Gamma_2(B_n(\St))$ of the sphere braid groups for $n\geq 4$, and in the case $n=4$, in \repr{presg2b4} we derive the presentation of $\Gamma_2(B_4(\St))$ given in \reth{lcsbn}(\ref{it:lcs4}).

\section{Extension to surfaces of higher genus}

Since work on this paper started, one of the authors, in collaboration with P.~Bellingeri and S.~Gervais has undertaken the study of the lower central series of braid groups of orientable surfaces, with and without boundary, of genus $g\geq 1$~\cite{BGG}. We remark that some of the techniques appearing in this monograph were used subsequently in that paper. It is worth stating the corresponding results of~\cite{BGG} which contrast somewhat with those obtained here for the sphere and punctured sphere. 

\begin{thm}[\cite{BGG}]
Let $M$ be a compact, connected orientable surface without boundary, of genus $g\geq 1$, and let $m\geq 3$. Then:
\begin{enumerate}[(a)]
\item  $\Gamma_1(B_m(M))/\Gamma_2(B_m(M)) \cong
  \Z^{2g} \oplus \Z_2$. 
\item $\Gamma_2(B_m(M))/\Gamma_3(B_m(M)) \cong \Z_{n-1+g}$.
\item $\Gamma_3(B_m(M))=\Gamma_4(B_m(M))$. Moreover, $\Gamma_3(B_m(M))$ is perfect for $m\geq 5$.
\end{enumerate}
\end{thm}
This implies that braid groups of compact, connected orientable surfaces without boundary may be distinguished by their lower central series (indeed by the first two lower central quotients).

\begin{thm}[\cite{BGG}]
Let $g\geq 1$, $q\geq 1$ and $m\geq 3$. Let $M$ be a
compact, connected orientable surface of genus $g$ with $q$ boundary
components. Then:
\begin{enumerate}[(a)]
\item $\Gamma_1(B_m(M))/ \Gamma_2(B_m(M))\cong \Z^{2g+q-1} \oplus \Z_2$.
\item $\Gamma_2(B_m(M))/ \Gamma_3(B_m(M))\cong \Z$.
\item $\Gamma_3(B_m(M))= \Gamma_4(B_m(M))$. Moreover, $\Gamma_3(B_m(M))$ is perfect for $m\geq 5$.
\end{enumerate}
\end{thm}
Thus if $m\geq 3$ and if $M$ a compact surface (with or without boundary) of genus $g\geq 1$, since $\Gamma_3(B_m(M))\neq \brak{1}$, $B_m(M)$ is not residually nilpotent. Moreover, we observe similar phenomena to those seen in \reth{lcdsbmsn} for the punctured sphere (stability of the lower central series for $m\geq 3$, perfectness of the $\Gamma_i(B_m(M))$ for $m\geq 5$). However, they occur one stage further, not from the commutator subgroup onwards, but from $\Gamma_3$ onwards.

Just as for $B_2(\St\setminus\brak{x_1,x_2})$, the $2$-string braid
groups represent a very difficult and interesting case. In the case of
the $2$-torus $\mathbb{T}^2$, we prove that its $2$-string braid group
is residually nilpotent. Further, arguing as in the proof of
\reth{lcsb2}, we show that apart from the first term, the lower
central series of $B_2(\mathbb{T}^2)$ and $\Z_2 \ast \Z_2 \ast \Z_2$
coincide, and by applying Gaglione's results, we may also determine
explicitly all of their successive lower central series quotients. More precisely:
\begin{thm}[\cite{BGG}]\mbox{}\label{th:bggtor}
\begin{enumerate}
\item $B_2(\mathbb{T}^2)$ is residually nilpotent.
\item For all $i\geq 2$:
\begin{enumerate}
\item $\Gamma_i(B_2(\mathbb{T}^2))\cong \Gamma_i(\Z_2 \ast \Z_2 \ast \Z_2)$.
\item $\Gamma_i(B_2(\mathbb{T}^2))/\Gamma_{i+1}(B_2(\mathbb{T}^2))$ is isomorphic to the direct sum of $R_i$ copies of $\Z_2$, where:
\begin{equation*}
R_i=\sum_{j=1}^{i-2} \left( \sum_{
\substack
{k \mid i-j\\
k>1}}
\; \mu \left( \frac{i-j}{k}\right) \frac{k\alpha_k}{i-j}
\right) \quad \text{and} \quad k\alpha_k=2^k+2(-1)^k.
\end{equation*}
\end{enumerate}
\end{enumerate}
\end{thm}
As in the case of the $2$-string braid group of the $n$-punctured sphere, $n\geq 3$, it seems to be very difficult even to describe the commutator subgroup of the $2$-string braid groups of orientable surfaces of higher genus.

\section*{Acknowledgements}

This work took place during the visit of the second author to the
Departmento de Matem\'atica do IME-Universidade de S\~ao Paulo during
the periods 5\up{th}--31\up{st}~August 2003,
9\up{th}~July--4\up{th}~August 2004, and
23\up{rd}~June--23\up{rd}~July 2005, and of the visit of the first
author to the Laboratoire de Math\'ematiques Emile Picard during the
period 30\up{th}~September--1\up{st}~November 2004. The first and
fourth of these visits were supported by the international Cooperation
Capes/Cofecub project number 364/01. The second visit was supported by
the `Accord franco-br\'esilien en math\'ematiques', and the third
visit by FAPESP.

\aufm{Daciberg Lima Gon\c{c}alves and John Guaschi}

\mainmatter

\chapter{The lower central series of $B_n(\St)$}\label{chap:sphere}

The main aim of this chapter is to prove \reth{lcsbn}, which describes the lower central series of $B_n(\St)$. This will be carried out in \resec{lcsbnst}. Before doing so, in \resec{general}, we state and prove some general results concerning the splitting of the short exact sequence~\reqref{ses} (\repr{bnabz}), as well as homological conditions for the stabilisation of the lower central series of a group (\relem{stallings}).

\section{Generalities}\label{sec:general}

Let $n\in\N$. Let $B_n(\St)$ denote the braid group of $\St$ on
$n$~strings, let $\bnab{n}= B_n(\St)/\Gamma_2(B_n(\St))$ denote the
Abelianisation of $B_n(\St)$, and let
$\map{\alpha}{B_n(\St)}[\bnab{n}]$ be the canonical projection.  Then
we have the following short exact sequence:
\begin{equation}\label{eq:ses}
\xymatrix{%
1\ar[r] & \Gamma_2(B_n(\St)) \ar[r] & B_n(\St) \ar[r]^-{\alpha} & \bnab{n} \ar[r] & 1.}
\end{equation} 

We first prove the following result which deals with the splitting of this short exact sequence. 
\begin{prop}\label{prop:bnabz}
Let $n\in\N$. 
\begin{enumerate}[(a)]
\item\label{it:bnabz1} $\bnab{n}=B_n(\St)/\Gamma_2(B_n(\St)) \cong \Z_{2(n-1)}$.
\item\label{it:bnabz2} The short exact sequence \reqref{ses} splits if and only if $n$ is odd, where the action on $\Gamma_2(B_n(\St))$ by a generator of $\Z_{2(n-1)}$ is given by conjugation by $\sigma_1\dotsc \sigma_{n-2}\sigma_{n-1}^2$.
\item If $n$ is even then $B_n(\St)$ is not isomorphic to the semi-direct product of a subgroup $K$ by $\Z_{2(n-1)}$.
\end{enumerate}
\end{prop}

\begin{proof}\mbox{}
\begin{enumerate}[(a)]
\item This follows easily from the presentation~\reqref{presnbns} of
the group $B_n(\St)$. The generators $\sigma_i$ of $B_n(\St)$ are all
identified by $\alpha$ to a single generator
$\widetilde{\sigma}=\alpha(\sigma_i)$ of $\Z_{2(n-1)}$.
\item\label{it:nodd} In order to construct a section, we consider the
elements of $B_n(\St)$ of order $2(n-1)$. According to Murasugi's
classification of the torsion elements of $B_n(\St)$~\cite{M}, these elements are precisely the conjugates of the
elements of the form $(\sigma_1\cdots\sigma_{n-2}\sigma_{n-1}^2)^r$,
where $r$ and $2(n-1)$ are coprime. Such an element projects to
$\widetilde{\sigma}^{rn}$ whose order is $2(n-1)/\pgcd{rn}{2(n-1)}$.
Since 
\begin{equation*}
\pgcd{rn}{2(n-1)}=\pgcd{n}{2(n-1)}=\pgcd{n}{2},
\end{equation*}
the result follows from \req{ses} and part~(\ref{it:bnabz1}).
\item Let $n\in\N$ be even. We first prove the following lemma: 

\begin{lem}\label{lem:hopf}
Let $G$ be a group whose Abelianisation $\gpab$ is Hopfian i.e.\ $\gpab$ is not isomorphic to any of its proper quotients. Suppose that there exists a group $H$ isomorphic to $\gpab$, a normal subgroup $K$ of $G$, and a split short exact sequence $1\to K\to G\to H\to 1$. Then $G\cong \Gamma_2(G)\rtimes \gpab$.
\end{lem}

\begin{proof}[Proof of \relem{hopf}]\mbox{}
Let $\map{\alpha}{G}[\gpab]$ denote Abelianisation, let $\map{\xi}{G}[H]$ denote the homomorphism in the given short exact sequence, and let $\map{s}{H}[G]$ be a section for $\xi$. Since $H\cong G/K$ is Abelian, it follows from standard properties of the commutator subgroup that $\Gamma_2(G)\subseteq K$. Hence we have the following commutative diagram:
\begin{equation*}
\xymatrix{%
1\ar[r] & \Gamma_2(G) \txt{~~}\ar@{^{(}->}[r] \ar@{^{(}->}[d] & G \ar[r]^-{\alpha}\ar@{=}[d] & \gpab \ar[r]  & 1\\
1\ar[r] & K \txt{~~}\ar@{^{(}->}[r] & G \ar@<4pt>[r]^-{\xi} &  H \ar[r] \ar@<4pt>@{-->}[l]^-{s} & 1,}
\end{equation*}
This extends to a commutative diagram of short exact sequences by taking $\map{\rho}{\gpab}[H]$ defined by $\rho(y)=\xi(x)$ for all $y\in \gpab$, where $x\in G$ is any element satisfying $\alpha(x)=y$. This homomorphism is well defined, and is surjective since $\xi$ and $\alpha$ are. But $\gpab\cong H$ is Hopfian by hypothesis, which implies that $\rho$ is an isomorphism. Hence $\alpha=\rho^{-1}\circ \xi$, and $s\circ \rho$ is a section for $\alpha$, which proves the lemma.
\end{proof}
By taking  $G=B_n(\St)$ and $K=\Z_{2(n-1)}$ in the statement of \relem{hopf}, if $B_n(\St)$ were a semi-direct product of $K$ with $H$ then this would contradict part~(\ref{it:nodd}). This completes the proof of \repr{bnabz}. \qedhere
\end{enumerate}
\end{proof}

\begin{rem} 
If $n$ is even, let us consider the natural projection $\map p{\Z_{2(n-1)}}[\Z_{n-1}]$. Then we have a short exact sequence:
\begin{equation*}
\xymatrix{%
1\ar[r] & \Gamma_2^{\ast}(B_n(\St)) \ar[r] & B_n(\St) \ar[r]^-{\alpha^{\ast}} &  \Z_{(n-1)} \ar[r] & 1.}
\end{equation*} 
where $\alpha^{\ast}=p\circ \alpha$, and $\Gamma_2^{\ast}(B_n(\St))$ is the kernel
of $\alpha^{\ast}$. It is not difficult to see that this short exact
sequence splits: a section is given by sending the generator of
$\Z_{(n-1)}$ to $(\sigma_1\dotsc \sigma_{n-2}\sigma_{n-1}^2)^{2^r}$,
where $2^r$ is the greatest power of $2$ dividing $n$.
\end{rem}

Let $G$ be a group which acts on a group $H$. Following~\cite[p.~67]{HMR}, we may define the commutator subgroup with respect to this action by
\begin{equation}\label{eq:hmr}
\Gamma_G(H)=\setang{(g\star h)\, kh^{-1}k^{-1}}{g\in G, h,k\in H},
\end{equation}
where $g\star h$ denotes the action of $g$ on $h$. We say that the action is \emph{perfect} if $\Gamma_G(H)=H$. Note that if $H$ is a normal subgroup of $G$ then $H\supseteq \Gamma_G(H)=[G,H] \supseteq [H,H]$ for the action of conjugation of $G$ on $H$. In particular, if $G=H$ then $\Gamma_G(H)=\Gamma_2(G)$ for the action of conjugation of $G$ on itself. If this action is perfect then the group $G$ is perfect.

\begin{lem}\label{lem:stallings}
Let $G$ be a group, and let $\gpab$ be its Abelianisation. Let $\map{\delta}{H_2(G,\Z)}[H_2(\gpab,\Z)]$ be the homomorphism induced by Abelianisation. Then 
\begin{equation*}
\Gamma_2(G)/\Gamma_3(G) \cong \coker{\delta}\cong H_0\left(\gpab, H_1\left(\Gamma_2(G), \Z\right)\right).
\end{equation*}
In particular:
\begin{enumerate}[(a)]
\item\label{it:stalla} $\Gamma_2(G)=\Gamma_3(G)$ if and only if $\delta$ is surjective.
\item\label{it:stallb} If $H_2(\gpab,\Z)$ is trivial then $\Gamma_n(G)= \Gamma_2(G)$ for all $n\geq 2$.
\item \label{it:stallc} If either the action (by conjugation) of $G$ on $\Gamma_2(G)$ or the action (by conjugation) of $\gpab$ on $H_1\left(\Gamma_2(G), \Z\right)$ is perfect then $\Gamma_n(G)= \Gamma_2(G)$ for all $n\geq 2$.
\end{enumerate}
\end{lem}

\begin{proof} 
Recall that if $1\to K\to G\to Q\to 1$ is an extension of groups then we have a $6$-term exact sequence 
\begin{equation}\label{eq:stallings}
H_2(G)\to H_2(Q)\to K/[G,K]\to H_1(G)\to H_1(Q)\to 1
\end{equation}
due to Stallings~\cite{Bro,McCl,St}. Applying this to the short exact sequence:
\begin{equation}\label{eq:gamma2}
1 \to \Gamma_2(G) \to G \to \gpab \to 1,
\end{equation}
we obtain:
\begin{equation*}
H_2(G,\Z) \stackrel{\delta}{\to} H_2(\gpab ,\Z) \to \Gamma_2(G)/\Gamma_3(G) \to H_1(G,\Z) \to \gpab \to 1.
\end{equation*} 
But $H_1(G,\Z) \to \gpab$ is an isomorphism, so this becomes 
\begin{equation*}
H_2(G,\Z) \stackrel{\delta}{\to} H_2(\gpab ,\Z) \to \Gamma_2(G)/\Gamma_3(G) \to 1.
\end{equation*}
Hence $\Gamma_2(G)/\Gamma_3(G) \cong \coker{\delta}$ which yields the
first isomorphism. To obtain the second,  we consider the
Lyndon-Hochschild-Serre spectral sequence~\cite{Bro,McCl} applied to
the short exact sequence~\reqref{gamma2}, for which the relevant terms
are $E_{(2,0)}^2= H_2(\gpab ,\Z)$ and $E_{(0,1)}^2= H_0\left(\gpab,
H_1\left(\Gamma_2(G), \Z\right)\right)$. Since $H_1(G)= H_1(\gpab)$,
the differential $\map{d_2}{E_{(2,0)}^2}[E_{(0,1)}^2]$ is surjective,
with kernel $E_{(2,0)}^{\infty}$. From the general definition of the
filtration of $H_2(G)$ given by the spectral sequence, we have a
surjection $H_2(G)\to E_{(2,0)}^{\infty}$, and hence the following
exact sequence:
\begin{equation*}
H_2(G)\to E_{(2,0)}^{\infty} \lhra E_{(2,0)}^2 \to E_{(0,1)}^2 \to 1.
\end{equation*}
Hence $\im{\delta}= E_{(2,0)}^{\infty}$, and 
\begin{equation*}
\coker{\delta}= E_{(2,0)}^2/\im{\delta}\cong E_{(0,1)}^2= H_0\left(\gpab, H_1\left(\Gamma_2(G), \Z\right)\right)
\end{equation*}
as required. From the first isomorphism, one may check that part~(\ref{it:stalla}) is satisfied. Part~(\ref{it:stallb}) then follows easily.

To prove part~(\ref{it:stallc}), if the action by conjugation of $G$ on $\Gamma_2(G)$ is perfect then $\Gamma_G(\Gamma_2(G))= [G,\Gamma_2(G)]=\Gamma_3(G)= \Gamma_2(G)$ and the result is clear. Now let us consider the action of $G$ on $H_1(\Gamma_2(G))= \gpab[(\Gamma_2(G))]$ given by conjugation, defined by $g\cdot \widetilde{h}=\widetilde{ghg^{-1}}$, where $g,h\in G$, and $\widetilde{~}$ denotes Abelianisation in $\Gamma_2(G)$. If $g\in \Gamma_2(G)$ then the induced action on $\gpab[(\Gamma_2(G))]$ is trivial, so the original action factors through $\gpab$, and we obtain an action of $\gpab$ on $\gpab[(\Gamma_2(G))]$ given by $\widetilde{g}\cdot \widetilde{h}=\widetilde{ghg^{-1}}$ ($\widetilde{g}$ denotes the Abelianisation of $g$ in $G$).  Suppose that this action is perfect, so that $\Gamma_{\gpab}(\gpab[(\Gamma_2(G))])= \gpab[(\Gamma_2(G))]$. Now 
\begin{equation*}
\Gamma_{\gpab}(\gpab[(\Gamma_2(G))])= [G, \Gamma_2(G)]/[\Gamma_2(G), \Gamma_2(G)]=\Gamma_3(G)/[\Gamma_2(G), \Gamma_2(G)],
\end{equation*}
and since $\Gamma_3(G) \subseteq \Gamma_2(G)$, it follows that $\Gamma_3(G)=\Gamma_2(G)$, which implies the result.
\end{proof}

\begin{rem}\label{rem:h2cyclic}
The hypothesis of part~(\ref{it:stallb}) of the lemma holds for example if $\gpab$ is cyclic. Recall that if $\gpab$ is finitely-generated then this condition is also necessary: if $H$ is a finitely-generated Abelian group satisfying $H_2(H,\Z)=\brak{0}$ then $H$ is cyclic.
\end{rem}

\section{The lower central series of $B_n(\St)$}\label{sec:lcsbnst}

Now we come to the main result of this chapter.

\begin{varthm}[\reth{lcsbn}] 
For all $n\geq 2$, the lower central series of $B_n(\St)$ is constant
from the commutator subgroup onwards: $\Gamma_m(B_n(\St))=
\Gamma_2(B_n(\St))$ for all $m\geq 2$. The subgroup $\Gamma_2(B_n(\St))$
is as follows: 
\begin{enumerate}
\item If $n=1,2$ then $\Gamma_2(B_n(\St))=\brak{1}$.
\item If $n=3$ then $\Gamma_2(B_n(\St)) \cong\Z_3$. Thus $B_3(\St)\cong \Z_3\rtimes \Z_4$, the action being the non-trivial one.
\item\label{it:lcs4} If $n=4$ then $\Gamma_2(B_4(\St))$ admits a presentation of the following form:
\begin{enumerate}
\item[\underline{\textbf{generators:}}] $g_1,g_2,g_3$, where in terms of the
usual generators of $B_4(\St)$, $g_1=\sigma_1^2\sigma_2 \sigma_1^{-3}$,
$g_2=\sigma_1^3\sigma_2 \sigma_1^{-4}$ and $g_3=\sigma_3 \sigma_1^{-1}$.
\item[\underline{\textbf{relations:}}]\label{it:g2bn4} 
\begin{gather}
g_3^4=1 \label{eq:g2b41}\\
[g_3^2,g_1]=1\label{eq:g2b42i}\\
[g_3^2,g_2]=1\label{eq:g2b42ii}\\
[g_3,g_2g_1]=1\label{eq:g2b43}\\
g_1^{-1}g_3^{-1}g_1=g_2g_3g_2^{-1}\label{eq:g2b44}\\ g_1^{-1}g_3^{-1}g_1=g_1g_3g_1^{-1}g_3. \label{eq:g2b45}
\end{gather}
\end{enumerate}
Furthermore, 
\begin{equation*}
\Gamma_2(B_4(\St))\cong \mathcal{Q}_8\rtimes \F[2](a,b),
\end{equation*}
where $\mathcal{Q}_8= \setang{x,y}{x^2=y^2,\; xyx^{-1}=y^{-1}}$ is the quaternion group of order~$8$, and $\F[2](a,b)$ is the free group of rank~$2$ on two generators $a$ and $b$. The action is given by:
\begin{align*}
&\phi(a)(x)=y &  &\phi(a)(y)=xy\\
&\phi(b)(x)=yx && \phi(b)(y)=x.
\end{align*}
\item In the cases $n=5$ and $n\geq 6$, a presentation for $\Gamma_2(B_n(\St))$ is given in \rechap{prescom}, by Propositions~\ref{prop:g2b5} and~\ref{prop:g2b6} respectively.
\end{enumerate}
\end{varthm}

\begin{proof} 
The first part of the theorem, $\Gamma_m(B_n(\St))= \Gamma_2(B_n(\St))$ for $m\geq 2$, follows from \relem{stallings}(\ref{it:stallb}) and \rerem{h2cyclic}. 

Now let us consider the rest of the theorem.
\begin{enumerate}
\item If $n=1,2$ then $B_n(\St)\cong \Z_n$, and the result follows easily.
\item Let $n=3$. Then $B_3(\St)$ is a $\text{ZS}$-metacyclic group (a
group whose Sylow subgroups, commutator subgroup and commutator
quotient group are all cyclic) of order~$12$~\cite{FvB}. It follows
from \repr{bnabz}(\ref{it:bnabz1}) that $\bnab{3}\cong \Z_4$, and
hence $\Gamma_2(B_3(\St))\cong \Z_3$. 

From \repr{bnabz}(\ref{it:bnabz2}), the short exact sequence \reqref{ses} splits, so $B_3(\St)\cong \Z_3\rtimes \Z_4$, and the action of the generator $\widetilde{\sigma}$ of $\bnab{3}$ on the generator $\rho$ of $\Z_3$ is given by $\widetilde{\sigma}\cdot\rho=\rho^{-1}$ i.e.\ the non-trivial action.

\item Let $n=4$. To obtain the given presentation of $\Gamma_2(B_4(\St))$, one applies the Reidemeister-Schreier rewriting process to the short exact sequence~\reqref{ses}. The calculations are deferred to \repr{presg2b4}, see \resec{presg2bn4} of \rechap{prescom}. 

Using this presentation, let us prove the second part of~(\ref{it:g2bn4}) of \reth{lcsbn}, that $\Gamma_2(B_4(\St))\cong \mathcal{Q}_8\rtimes \F[2](a,b)$. This will be achieved by the following two propositions.

\begin{prop}\label{prop:normq8}
The normal subgroup of $\Gamma_2(B_4(\St))$ generated by $g_3$ is isomorphic to a quotient of the quaternion group $\mathcal{Q}_8$.
\end{prop}

\begin{proof}\mbox{}
Let $N$ be the normal subgroup of $\Gamma_2(B_4(\St))$ generated by
$g_3$, and let $H$ be the subgroup of $\Gamma_2(B_4(\St))$ generated
by $g_3$ and $g_1g_3g_1^{-1}$. Clearly $H\subseteq N$. To prove the
converse, it suffices to show that if we conjugate $g_3$ and
$g_1g_3g_1^{-1}$ by $g_1^{\pm 1}$ and $g_2^{\pm 1}$, we obtain elements of $H$. This is a consequence of the following equalities:
\begin{align*}
g_2g_3g_2^{-1}&= g_1^{-1}g_3^{-1}g_1 \quad\text{by \req{g2b44}}\\
&=g_1g_3g_1^{-1}\cdot g_3 \quad\text{by \req{g2b45}}\\
g_1^2g_3g_1^{-2} &= g_3^{-1}\cdot g_1 g_3^{-1}g_1^{-1} \quad \text{by \req{g2b45}}\\
g_2g_1g_3g_1^{-1}g_2^{-1}&=g_3\quad\text{by \req{g2b43}}\\
g_2^{-1}g_3 g_2 &= g_1g_3g_1^{-1}\quad\text{by \req{g2b43}}\\
g_2^{-1}g_1g_3g_1^{-1} g_2&= g_2^{-1}g_1^{-1}g_3^{-1}g_1 g_3^{-1}g_2 \quad\text{by \req{g2b45}}\\
&= g_3 \cdot g_2^{-1}g_3^{-1} g_2\quad\text{by \req{g2b44}.}
\end{align*}
Hence $H=N$ is normal in $\Gamma_2(B_4(\St))$. Now $g_3^2=(g_1g_3g_1^{-1})^2$ by \req{g2b42i}, and $(g_1g_3g_1^{-1} g_3)^2=(g_1^{-1}g_3^{-1}g_1)^2=g_3^{-2}= g_3^2$ by equations~\reqref{g2b45} and~\reqref{g2b41}. By equations~\reqref{g2b41} and~\reqref{g2b42i} it thus follows that $[g_1g_3g_1^{-1}, g_3]=g_3^2$, and hence $g_1g_3g_1^{-1}\cdot g_3 g_1g_3^{-1}g_1^{-1}=g_3^3= g_3^{-1}$. So $g_1g_3g_1^{-1}$ and $g_3$ satisfy a set of defining relations of $\mathcal{Q}_8$, and thus $H$ is a quotient of $\mathcal{Q}_8$.
\end{proof}

\begin{prop}
With $H$ as defined as in the proof of \repr{normq8}, $H\cong \mathcal{Q}_8$, and $\Gamma_2(B_4(\St)) \cong \mathcal{Q}_8\rtimes \F[2](a,b)$, the action being given by $\phi(a)(x)=axa^{-1}=y$, $\phi(a)(y)=aya^{-1}=xy$, $\phi(b)(x)=bxb^{-1}=yx$ and $\phi(b)(y)=byb^{-1}=x$.
\end{prop}

\begin{proof}\mbox{}
Let $\mathcal{Q}_8$ be generated by $x$ and $y$, subject to the
relations $x^2=y^2$ and $xyx^{-1}= y^{-1}$. We remark that if $z\in
\mathcal{Q}_8$ and $w\in \F[2](a,b)$ then $wzw^{-1}= \phi(w)(z)$, and
$[z,w]= z\cdot \phi(w)(z^{-1})$. Consider the map 
\begin{equation*}
\map{\psi}{\brak{g_1,g_2,g_3}}[\mathcal{Q}_8\rtimes\FF_2(a,b)]
\end{equation*}
defined as follows: $\psi(g_1)=a$, $\psi(g_2)=b$ and $\psi(g_3)=x$. It
is straightforward to check that the images under $\psi$ of
relations~\reqref{g2b41}--\reqref{g2b43} hold in 
$\mathcal{Q}_8\rtimes \F[2](a,b)$. As for relation~\reqref{g2b44}, the
right-hand side yields $bxb^{-1}=\phi(b)(x)= yx$ from the definition
of the action, while the left-hand side yields $a^{-1}x^{-1}a=
\phi(a^{-1})(x^{-1})$. Now $\phi(a)(xy^{-1})=x^{-1}$, so
$\phi(a^{-1})(x^{-1})=xy^{-1}=yx$ in $\mathcal{Q}_8$. So
relation~\reqref{g2b44} is preserved under $\psi$. Finally, consider
relation~\reqref{g2b45}. From the previous relation, the left-hand
side yields $yx$. As for the right-hand side, we obtain
$\phi(a)(x)\cdot x=yx$ also. So $\psi$ extends to a homomorphism,
which we also call $\psi$, from $\Gamma_2(B_4(\St))$ into
$\mathcal{Q}_8\rtimes \F[2](a,b)$. Since $\psi(g_1g_3g_1^{-1})=y$,
this homomorphism is certainly surjective. Further, since the normal
subgroup $H$ of \repr{normq8} is generated by $g_3$ and
$g_1g_3g_1^{-1}$, it follows that $H$ is mapped surjectively onto
$\mathcal{Q}_8$. But $H$ is a quotient of $\mathcal{Q}_8$, and since
$\mathcal{Q}_8$ is finite, $H$ is isomorphic to $\mathcal{Q}_8$. This
proves the first part of the proposition. The induced map from the
quotient $\Gamma_2(B_4(\St))$ by $H$ (which is the normal subgroup
generated by $g_3$) into the quotient of $\mathcal{Q}_8\rtimes
\F[2](a,b)$ by $\mathcal{Q}_8$ is a surjective homomorphism from a
free group on two generators into a free group on two generators, so
is an isomorphism by the Hopfian property of free groups of finite rank. This completes the proof of the proposition, as
well as that of part~(\ref{it:g2bn4}) of \reth{lcsbn}.
\end{proof}

\item Now suppose that $n\geq 5$. The presentations are given in \rechap{prescom}, Propositions~\ref{prop:g2b5} and~\ref{prop:g2b6} respectively. This completes the proof of \reth{lcsbn}. \qedhere
\end{enumerate}
\end{proof}

\chapter{The derived series of $B_n(\St)$}\label{chap:dssphere}

In this chapter, we study the derived series of $B_n(\St)$. The aim is to prove the following result, which shows that for all $n\neq 3,4$, $(B_n(\St))^{(1)}$ is perfect. The difficult case is $n=4$, but using the semi-direct product structure of $(B_4(\St))^{(1)}$ obtained in \reth{lcsbn}, we shall be able to prove that the derived series of $B_4(\St)$ coincides from a certain point with that of the free group of rank~$2$. Before doing so, we state and prove \repr{gammasemi} which describes the commutator subgroup of a general semi-direct product.

\begin{varthm}[\reth{dsbn}]
The derived series of $B_n(\St)$ is as follows.
\begin{enumerate}[(a)]
\item \label{it:ds1} If $n=1,2$ then $(B_n(\St))^{(1)}=\brak{1}$.
\item \label{it:ds2} If $n=3$ then $(B_3(\St))^{(1)}\cong\Z_3$ and $(B_3(\St))^{(2)}=\brak{1}$.
\item  \label{it:ds3} Suppose that $n=4$. Then:
\begin{enumerate}[(i)]
\item \label{it:ds4gam} $(B_4(\St))^{(1)}=\Gamma_2(B_4(\St))$ is given by part~(\ref{it:lcs4}) of \reth{lcsbn}; it is isomorphic to the semi-direct product $\mathcal{Q}_8\rtimes \F[2]$. Further, $B_4(\St)/(B_4(\St))^{(1)}$ is isomorphic to $\Z_6$.
\item \label{it:ds4a} $(B_4(\St))^{(2)}$ is isomorphic to the semi-direct product $\mathcal{Q}_8\rtimes (\F[2])^{(1)}$, where $(\F[2])^{(1)}$ is the commutator subgroup of the free group $\F[2](a,b)$ of rank~$2$ on two generators $a,b$. The action of $(\F[2])^{(1)}$ on $\mathcal{Q}_8$ is the restriction of the action of $\F[2](a,b)$ given in part~(\ref{it:lcs4}) of \reth{lcsbn}. Further, 
\begin{equation*}
\text{$(B_4(\St))^{(1)}/(B_4(\St))^{(2)}\cong \Z^2$, and $B_4(\St)/(B_4(\St))^{(2)} \cong \Z^2\rtimes \Z_6$,}
\end{equation*}
where the action of the generator $\overline{\sigma}$ of $\Z_6$ on $\Z^2$ is given by left multiplication by the matrix $\left( \begin{smallmatrix}
0 & 1\\
-1 & 1
\end{smallmatrix}\right)$.
\item \label{it:ds4b} $(B_4(\St))^{(3)}$ is a subgroup of $P_4(\St)$ isomorphic to the direct product $\Z_2\times (\F[2])^{(2)}$. Further,
\begin{equation*}
\text{$(B_4(\St))^{(2)}/(B_4(\St))^{(3)}\cong (\Z_2 \times \Z_2)\times (\F[2])^{(1)}/(\F[2])^{(2)}$.}
\end{equation*}
\item \label{it:ds4c} $(B_4(\St))^{(m)}\cong (\F[2])^{(m-1)}$ for all $m\geq 4$. Further, 
\begin{equation*}
(B_4(\St))^{(3)}/(B_4(\St))^{(4)}\cong \Z_2\times (\F[2])^{(2)}/(\F[2])^{(3)},
\end{equation*}
and for $m\geq 4$, 
\begin{equation*}
(B_4(\St))^{(m)}/(B_4(\St))^{(m+1)}\cong (\F[2])^{(m-1)}/(\F[2])^{(m)}.
\end{equation*}
\end{enumerate}
\item \label{it:ds4} If $n\geq 5$ then $(B_n(\St))^{(2)}=(B_n(\St))^{(1)}$, so $(B_n(\St))^{(1)}$ is perfect. A presentation of $(B_n(\St))^{(1)}$ is given in Propositions~\ref{prop:g2b5} and~\ref{prop:g2b6}.
\end{enumerate}
\end{varthm}

\begin{rem}\label{rem:freederive}
In part~(\ref{it:ds3}) of \reth{dsbn} and also in what follows, we shall often refer to the derived series of $\F[2](a,b)$ as well as its quotients. We were not able to track down an explicit reference for them, but one may observe that for $i\geq 1$, $(\F[2](a,b))^{(i)}$ is a free group of infinite rank, and hence $(\F[2](a,b))^{(i)}/(\F[2](a,b))^{(i+1)}$ is a free Abelian group of infinite rank. A basis of $(\F[2](a,b))^{(1)}=\Gamma_2(\F[2](a,b))$ may be obtained as follows: considering the short exact sequence~\reqref{gamma2} with $G=\F[2](a,b)$, $(\F[2](a,b))^{(1)}$ may be identified with the fundamental group of the Cayley graph of $\F[2](a,b)$. Let $\mathcal{T}$ be a maximal tree  in this graph. For each $g\in \F[2](a,b)$, let $w_g$ be the word corresponding to the path in $\mathcal{T}$ between $e$ and $g$. Then a basis is given by the set of elements of the form $w_g [a,b] w_g^{-1}$, where $g$ runs over $\F[2](a,b)$. For example, the set $\brak{a^pb^q [a,b] b^{-q}a^{-p}}_{p,q\in\Z}$ is a basis of $(\F[2](a,b))^{(1)}$. Since $\F[2](a,b)$ is residually nilpotent and $(\F[2](a,b))^{(i-1)} \subseteq \Gamma_i(\F[2](a,b))$, it follows that $\bigcap_{i\geq 0} \; (\F[2](a,b))^{(i)}=\brak{1}$ and $\F[2](a,b)$ is residually soluble.
\end{rem}

We obtain easily the following corollary of \reth{dsbn}:
\begin{cor}\label{cor:residb4}
Let $n\in\N$. Then $B_n(\St)$ is residually soluble if and only if $n\leq 4$.
\end{cor}

\begin{proof}[Proof of \reco{residb4}]
Recall that a group $G$ is residually soluble if and only if $\bigcap_{i\geq 0} \; G^{(i)}=\brak{1}$. If $n=1,2,3$, this is obvious, and if $n=4$, the residual solubility of $B_4(\St)$ follows from that of $\F[2](a,b)$. For $n\geq 5$, the result also follows easily, since $(B_n(\St))^{(1)}$ is non trivial.
\end{proof}

Before proving \reth{dsbn}, let us state and prove the following proposition which describes the commutator subgroup of a semi-direct product. This result will be used frequently throughout the rest of this paper.

\begin{prop}\label{prop:gammasemi}
Let $G,H$ be groups, and let $\map{\phi}{G}[\aut{H}]$ be an action of $G$ on $H$. Let $\widehat{H}$ be the subgroup of $H$ generated by the elements of the form $\phi(g)(h)\cdot h^{-1}$, where $g\in G, h\in H$, and let $L$ be the subgroup of $H$ generated by $\Gamma_2(H)$ and $\widehat{H}$. Then $\phi$ induces an action (also denoted by $\phi$) of $\Gamma_2(G)$ on $L$, and $L\rtimes_{\phi} \Gamma_2(G) = \Gamma_2(H\rtimes_{\phi} G)$. In particular, $\Gamma_2(H\rtimes_{\phi} G)$ is the subgroup generated by $\Gamma_2(H),\Gamma_2(G)$ and $\widehat{H}$.
\end{prop}

\begin{rem}\label{rem:gammagh}
We claim that $L$ is none other than the commutator subgroup
$\Gamma_G(H)$ defined by \req{hmr} with respect to the given action.
To see this, recall that $\Gamma_G(H)$ is the subgroup of $H$
generated by the elements of the form $\phi(g)(h)\cdot kh^{-1}k^{-1}$,
where $g\in G$ and $h,k\in H$. Taking $g=e$ (respectively $k=e$), it
follows that $\Gamma_G(H) \supseteq \Gamma_2(H)$ (respectively
$\Gamma_G(H) \supseteq \widehat{H}$), and hence $L \subseteq
\Gamma_G(H)$. Conversely, $\phi(g)(h)\cdot kh^{-1}k^{-1}=
\phi(g)(h)h^{-1}\cdot hkh^{-1}k^{-1}\in L$, so $\Gamma_G(H) \subseteq
L$, and the claim is proved. Note further that if $h'\in H$ then there
exists $h''\in H$ such that $\phi(g)(h'')=h'$, so
\begin{equation*}
h'\left(\phi(g)(h)\cdot h^{-1}\right)h'^{-1}= \phi(g)(h''h) (h''h)^{-1}\cdot  (\phi(g)(h'')h''^{-1})^{-1}.
\end{equation*}
It follows that $\widehat{H}$ and $L$ are normal in $H$.  In particular, $\Gamma_G(H)$ is normal in $H$.
\end{rem}

\begin{proof}[Proof of \repr{gammasemi}]
From now on, we shall identify each subgroup $H_1$ of $H$ (respectively each subgroup $G_1$ of $G$) with the corresponding subgroup $\set{(h,1)}{h\in H_1}$ (respectively $\set{(1,g)}{g\in G_1}$) of $H\rtimes_{\phi} G$ without further comment. The group operation in $H\rtimes_{\phi} G$ shall be written as:
\begin{equation*}
\text{$(h,g)\star (h',g')= (h\ldotp \phi(g)(h),gg')$, where $(h,g),(h',g')\in H\rtimes_{\phi} G$.}
\end{equation*}
The subgroup $L$ is normal in $H$ by \rerem{gammagh}. Let us show that $\phi$ induces an action (also denoted by $\phi$) of $G$ on $L$. Let $g\in G$. Since $\phi(g)([h_1,h_2])= [\phi(g)(h_1), \phi(g)(h_2)]\in \Gamma_2(H)$ for all $h_1,h_2\in H$, and 
\begin{equation*}
\phi(g)(\phi(g')(h)h^{-1}) =  \phi(gg')(h)h^{-1}\ldotp
(\phi(g)(h)h^{-1})^{-1} \in \widehat{H}
\end{equation*}
for all $h\in H$ and $g'\in G$, it follows that $\phi(g)(L)\subseteq L$. Clearly $\phi(g)$ is injective. The surjectivity of $\phi(g)$ (restricted to $L$) may be deduced from the following observations:
\begin{enumerate}
\item if $i=1,2$ and $h'_i\in H$ then there exists $h_i\in H$ such that $\phi(g)(h_i)=h'_i$, and hence $\phi(g)([h_1,h_2])= [h'_1,h'_2]$.
\item If $g'\in G$ and $h,h'\in H$ then 
\begin{equation*}
\phi(g)\left(\phi(g^{-1}g')(h)h^{-1}\ldotp h\left(\phi(g^{-1})(h^{-1})h\right)h^{-1}\right)=
\phi(g')(h)h^{-1}.
\end{equation*}
\end{enumerate}
Thus $\phi$ induces an action (also denoted by $\phi$) of $\Gamma_2(G)$ on $L$, and $L\rtimes_{\phi} \Gamma_2(G)$ is a subgroup of $H\rtimes_{\phi} G$.

Clearly any element of $\Gamma_2(H)$ (respectively $\Gamma_2(G)$) may be written as an element of $\Gamma_2(H\rtimes_{\phi} G)$. Further, if $g\in G$ and $h\in H$ then
\begin{equation*}
\left[(1,g), (h,1)\right]=(\phi(g)(h),1) \star (h^{-1},1)= (\phi(g)(h) h^{-1},1),
\end{equation*}
and thus every element of $\widehat{H}$ may be written as an element of $\Gamma_2(H\rtimes_{\phi} G)$. This proves that $L\rtimes_{\phi} \Gamma_2(G) \subseteq \Gamma_2(H\rtimes_{\phi} G)$.

To see the converse, notice that the commutator of two elements $(h_1,g_1),(h_2,g_2)\in H\rtimes_{\phi} G$ may be written as:
\begin{multline*}
[(h_1,g_1),(h_2,g_2)]=\bigl(h_1\ldotp \phi(g_1)(h_2)\ldotp \phi(g_1g_2g_1^{-1})(h_1^{-1})\ldotp\\
\phi([g_1,g_2])(h_2^{-1}),[g_1,g_2]\bigr).
\end{multline*} 
The second factor belongs clearly to $\Gamma_2(G)$. The first factor is of the form:
\begin{multline*}
[h_1,h_2] \ldotp h_2h_1h_2^{-1} \left(\phi(g_1)(h_2) h_2^{-1}\right) h_2h_1^{-1}h_2^{-1} \ldotp\\
 h_2h_1 \left(\phi(g_1g_2g_1^{-1})(h_1^{-1}) h_1\right) h_1^{-1}h_2^{-1}\ldotp h_2 \left(\phi([g_1,g_2])(h_2^{-1}) h_2\right) h_2^{-1},
\end{multline*}
which is a product of elements of $L$. Hence $\Gamma_2(H\rtimes_{\phi} G) \subseteq L\rtimes_{\phi} \Gamma_2(G)$, and the proposition follows.
\end{proof}

We now prove the main result of this chapter.
\begin{proof}[Proof of \reth{dsbn}]
Cases~(\ref{it:ds1}) and~(\ref{it:ds2}) follow directly from \reth{lcsbn}.

Now consider case~(\ref{it:ds4}), i.e.\ $n\geq 5$. 
Let $H\subseteq (B_n(\St))^{(1)}$ be a normal subgroup of $B_n(\St)$ such that $A=(B_n(\St))^{(1)}/H$ is Abelian (notice that this condition is satisfied if 
$H=(B_n(\St))^{(2)}$). Let
\begin{equation*}
\left\{ \begin{aligned}
\pi \colon\thinspace B_n(\St) & \to B_n(\St)/H\\
\beta &\mapsto \overline{\beta}
\end{aligned}\right.
\end{equation*}
denote the canonical projection. So the Abelianisation homomorphism $\map{\alpha}{B_n(\St)}[{\bnab{n}}]$ of \rechap{sphere} factors through $B_n(\St)/H$ i.e.\ there exists a (surjective) homomorphism $\map{\widehat{\alpha}}{B_n(\St)/H}[\bnab{n}]$ satisfying $\alpha= \widehat{\alpha}\circ \pi$. So we have the following short exact sequence:
\begin{equation*}
\xymatrix{%
1\ar[r] & A \ar[r] & B_n(\St)/H \ar[r]^{\widehat{\alpha}} & \bnab{n} \ar[r] & 1.}
\end{equation*} 
Now $\overline{\sigma_1},\ldots, \overline{\sigma_{n-1}}$ generate $B_n(\St)/H$, but since $\alpha(\sigma_i)=\alpha(\sigma_1)$ for $1\leq i\leq n-1$, it follows that $\widehat{\alpha}(\overline{\sigma_i})= \widehat{\alpha}(\overline{\sigma_1})$, and so there exists $t_i\in A$ such that $\overline{\sigma_i}=t_i\overline{\sigma_1}$.

We now apply $\pi$ to each of the relations of \req{presnbns} of $B_n(\St)$. First suppose that $3\leq i\leq n-1$. Since $\sigma_i$ commutes with $\sigma_1$, we have that
\begin{equation*}
\overline{\sigma_1}\cdot t_i\overline{\sigma_1} =t_i\overline{\sigma_1}\cdot \overline{\sigma_1},
\end{equation*}
and hence $t_i$ commutes with $\overline{\sigma_1}$.

Now let $4\leq i\leq n-1$. Since $\sigma_i$ commutes with $\sigma_2$, we obtain
\begin{equation*}
t_i\overline{\sigma_1}\cdot t_2\overline{\sigma_1} =t_2\overline{\sigma_1}\cdot t_i\overline{\sigma_1}.
\end{equation*}
Since $A$ is Abelian, it follows from the previous paragraph that $t_2$ commutes with $\overline{\sigma_1}$. Applying this to the image of the relation $\sigma_1\sigma_2\sigma_1= \sigma_2\sigma_1\sigma_2$, under $\pi$, we see that $t_2=t_2^2$, and hence $t_2=1$.

Next, if $i\geq 2$ then the relation $\sigma_i\sigma_{i+1}\sigma_i=
\sigma_{i+1}\sigma_i\sigma_{i+1}$ implies that $t_i=t_{i+1}$, and so
$t_2=\dotsc=t_{n-1}=1$. Hence $\overline{\sigma_1}=\overline{\sigma_2}=\dotsc
=\overline{\sigma_{n-1}}$. Thus $B_n(\St)/H$ is cyclic, generated by
$\overline{\sigma_1}$, and finite of order not greater than $2(n-1)$, because
the surface relation $\sigma_1\dotsc\sigma_{n-2}\sigma_{n-1}^2
\sigma_{n-2}\dotsc \sigma_1=1$ projects to $\overline{\sigma_1}^{2(n-1)}=1$.
Since $\widehat{\alpha}$ is surjective and $\bnab{n}\cong \Z_{2(n-1)}$, we
conclude that $\widehat{\alpha}$ is an isomorphism, so $B_n(\St)/H\cong
\Z_{2(n-1)}$, and $A=(B_n(\St))^{(1)}/H$ is trivial. In particular 
\begin{equation*}
(B_n(\St))^{(2)}=\bigl[
(B_n(\St))^{(1)},(B_n(\St))^{(1)} \bigr]= (B_n(\St))^{(1)},
\end{equation*}
in other words, $(B_n(\St))^{(1)}$ is perfect.

Now consider case~(\ref{it:ds3}), so $n=4$. Recall that part~(\ref{it:ds4gam})
was proved in \reth{lcsbn} and \repr{bnabz}. To obtain $(B_4(\St))^{(2)}$, it
suffices to observe that for the action of $\F[2](a,b)$ on $\mathcal{Q}_8$, the
subgroup $\widehat{\mathcal{Q}_8}$ defined in \repr{gammasemi} is
$\mathcal{Q}_8$ (which is the case, since by \reth{lcsbn}(\ref{it:lcs4}),
$\phi(b)(x)x^{-1}=y$ and $\phi(a)(y)y^{-1}=x$). So $(B_4(\St))^{(2)}$ is
generated by $\mathcal{Q}_8$ and $(\F[2])^{(1)}$, $(B_4(\St))^{(2)}\cong
\mathcal{Q}_8\rtimes (\F[2])^{(1)}$, and the action is the restriction of that
of $\F[2](a,b)$ on $\mathcal{Q}_8$, which proves the first part
of~(\ref{it:ds3})(\ref{it:ds4a}).

To determine $(B_4(\St))^{(3)}$, we first have to describe the subgroup
$\widehat{\mathcal{Q}_8}$ for the action of $(\F[2])^{(1)}$ on $\mathcal{Q}_8$.
By \reth{lcsbn}(\ref{it:lcs4}), if $B=[a,b]\in (\F[2](a,b))^{(1)}$ then the
automorphism $\phi(B)$ satisifies $\phi(B)(z)=x^2\cdot z$ for $z\in\brak{x,y}$
(recall that $x^2=y^2$). Since $(\F[2](a,b))^{(1)}$ is the subgroup of
$\F[2](a,b)$ normally generated by $B$, and the centre $\ang{x^2}$ of
$\mathcal{Q}_8$ is invariant under $\aut{\mathcal{Q}_8}$, it follows that
$\widehat{\mathcal{Q}_8}=\ang{x^2}$. So $(B_4(\St))^{(3)}$ is isomorphic to the
semi-direct product of $\Z_2$ by $(\F[2])^{(2)}$. But the action is trivial, and
so the product is direct. This proves the first part
of~(\ref{it:ds3})(\ref{it:ds4b}).

For $m\geq 4$, the subgroup $(B_4(\St))^{(m)}$ is clear from the description of $(B_4(\St))^{(3)}$, and hence we obtain the first part of~(\ref{it:ds3})(\ref{it:ds4c}).

We now analyse various quotients of the form
$B_4(\St)/(B_4(\St))^{(m)}$ and $(B_4(\St))^{(m-1)}/(B_4(\St))^{(m)}$
for several values of $m$. For the quotient
$B_4(\St)/(B_4(\St))^{(m)}$, we shall consider the case $m=2$ (the
case $m=1$ is given by \repr{bnabz}(\ref{it:bnabz1})). For 
$(B_4(\St))^{(m-1)}/(B_4(\St))^{(m)}$, we consider the cases $m\geq 2$
(the case $m=1$ was considered in \repr{bnabz}(\ref{it:bnabz1})). If
$m>4$, the problem reduces to the corresponding problem for the free
group on two generators. 

We adopt the notation used above in the case $n\geq 5$, and again we suppose
that $H\subseteq (B_n(\St))^{(1)}$ is a normal subgroup of $B_n(\St)$
such that  $A= (B_n(\St))^{(1)}/H$ is Abelian. So we have a short
exact sequence:
\begin{equation*}
\entrymodifiers={!! <0pt, .8ex>+}
\xymatrix{%
1\ar[r] & A \ar[r] & B_4(\St)/H \ar[r]^{\widehat{\alpha}} & {\underbrace{\bnab{4}}_{\Z_6}} \ar[r] & 1.}
\end{equation*}
Now $\overline{\sigma_1}, \overline{\sigma_2}, \overline{\sigma_3}$
generate $B_4(\St)/H$. As above, for $i=2,3$ we set
$\overline{\sigma_i}=t_i\overline{\sigma_1}$, where $t_i\in A$, and
we apply $\pi$ to the relations of $B_4(\St)$. The fact that
$\sigma_1$ commutes with $\sigma_3$ implies that $t_3$ commutes with
$\overline{\sigma_1}$. The relation $\sigma_1\sigma_2\sigma_1=
\sigma_2\sigma_1\sigma_2$ implies that:
\begin{equation}\label{eq:eq1a}
\overline{\sigma_1}t_2\overline{\sigma_1}^{-1} = t_2\cdot \overline{\sigma_1}^2 t_2 \overline{\sigma_1}^{-2}.
\end{equation}

Now consider the relation $\sigma_2\sigma_3\sigma_2= \sigma_3\sigma_2\sigma_3$. We have that:
\begin{equation*}
t_2\overline{\sigma_1}\cdot t_3\overline{\sigma_1}\cdot t_2 \overline{\sigma_1} = t_3\overline{\sigma_1}\cdot t_2\overline{\sigma_1}\cdot t_3 \overline{\sigma_1},
\end{equation*}
and so 
\begin{equation*}
t_2 \overline{\sigma_1}^2 t_2=t_3 \overline{\sigma_1} t_2 \overline{\sigma_1},
\end{equation*}
since $A$ is Abelian and $t_3$ commutes with $\overline{\sigma_1}$. Thus:
\begin{equation*}
t_2 \overline{\sigma_1}^2 t_2= t_3 t_2 \overline{\sigma_1}^2 t_2
\end{equation*}
from \req{eq1a}, and so $t_3=1$. We conclude that $B_4(\St)/H$ is generated by $\overline{\sigma_1}$ and $t_2 \overline{\sigma_1}$.

Finally, we consider the image of the surface relation under $\pi$. Using \req{eq1a}, note first that:
\begin{align}
\overline{\sigma_1}^3t_2 \overline{\sigma_1}^{-3}& = \overline{\sigma_1}(t_2^{-1}\cdot \overline{\sigma_1} t_2 \overline{\sigma_1}^{-1})\overline{\sigma_1}^{-1}= \overline{\sigma_1}t_2^{-1}\overline{\sigma_1}^{-1}\cdot \overline{\sigma_1}^2 t_2 \overline{\sigma_1}^{-2}\notag\\
&= \overline{\sigma_1}t_2^{-1}\overline{\sigma_1}^{-1}\cdot t_2^{-1} \overline{\sigma_1}t_2\overline{\sigma_1}^{-1}=t_2^{-1}, \label{eq:eq2a}
\end{align}
since $A$ is normal and Abelian. Thus $\sigma_1\sigma_2\sigma_3^2\sigma_2\sigma_1=1$ implies that:
\begin{align*}
1 &=\overline{\sigma_1}\cdot t_2 \overline{\sigma_1}\cdot \overline{\sigma_1}^2 \cdot t_2 \overline{\sigma_1}\cdot \overline{\sigma_1}= \overline{\sigma_1} t_2 \overline{\sigma_1}^{-1}\cdot \overline{\sigma_1} (\overline{\sigma_1}^3 t_2 \overline{\sigma_1}^{-3}) \overline{\sigma_1}^{-1}\cdot \overline{\sigma_1}^6\\
&= \overline{\sigma_1} t_2 \overline{\sigma_1}^{-1}\cdot \overline{\sigma_1} t_2^{-1} \overline{\sigma_1}^{-1}\cdot \overline{\sigma_1}^6= \overline{\sigma_1}^6
\end{align*}
from \req{eq2a}. 

Recall that $\Gamma_2(B_4(\St))$ is the normal subgroup of $B_4(\St)$ generated
by the commutators of the generators of $B_4(\St)$. Hence $A$ is the normal
subgroup of $B_4(\St)/H$ generated by $[\overline{\sigma_1},
t_2\overline{\sigma_1}]= \overline{\sigma_1}t_2 \overline{\sigma_1}^{-1} \cdot
t_2^{-1}$. Since $A$ is Abelian and $t_2\in A$, the action of conjugation on $A$
by $t_2$ is trivial. From \req{eq2a}, the action of $\overline{\sigma_1}^3$ on
$t_2$ yields $t_2^{-1}$. Further,
\begin{equation*}
\overline{\sigma_1} (\overline{\sigma_1} t_2\overline{\sigma_1}^{-1} t_2^{-1})\overline{\sigma_1}^{-1}= \overline{\sigma_1}^2 t_2 \overline{\sigma_1}^{-2}\cdot \overline{\sigma_1} t_2^{-1} \overline{\sigma_1}^{-1}=t_2^{-1}
\end{equation*}
from \req{eq1a}, and since
\begin{equation*}
\overline{\sigma_1}^2 (\overline{\sigma_1} t_2\overline{\sigma_1}^{-1} t_2^{-1})\overline{\sigma_1}^{-2}= \overline{\sigma_1} t_2^{-1} \overline{\sigma_1}^{-1},
\end{equation*}
it follows that $A$ is the Abelian group generated by $\overline{\sigma_1}t_2\overline{\sigma_1}^{-1} t_2^{-1}$, $t_2$ and $\overline{\sigma_1} t_2 \overline{\sigma_1}^{-1}$, and thus by $t_2$ and $\overline{\sigma_1}t_2\overline{\sigma_1}^{-1}$.

Let $\widetilde{\sigma}=\alpha(\sigma_1)$ denote the generator of $\bnab{4}$. Let $M=\begin{pmatrix}
0 & 1\\
-1 & 1
\end{pmatrix}$; notice that $M$ is of order~$6$. We now let $\bnab{4}\cong \Z_6$ act on $\Z^2$ as follows:
\begin{equation*}
\widetilde{\sigma}\cdot 
\begin{pmatrix}
X_1\\
X_2
\end{pmatrix}
=M\begin{pmatrix}
X_1\\
X_2
\end{pmatrix}=
\begin{pmatrix}
X_2\\
X_2-X_1
\end{pmatrix},
\end{equation*}
and so we may form the associated semi-direct product $\Z^2 \rtimes \Z_6$. We now consider the following homomorphism:
\begin{align*}
\psi\thinspace\colon B_4(\St) & \to \Z^2 \rtimes \Z_6\\
\sigma_1,\sigma_3 & \mapsto \left( \begin{pmatrix}
0\\
0
\end{pmatrix}, \widetilde{\sigma}\right)\\
\sigma_2 & \mapsto \left( \begin{pmatrix}
1\\
0
\end{pmatrix}, \widetilde{\sigma}\right).
\end{align*}
We then check that $\psi$ is well defined: clearly
$\psi(\sigma_1\sigma_3)= \psi(\sigma_3\sigma_1)$. To see that
$\psi(\sigma_1\sigma_2\sigma_1)= \psi(\sigma_2\sigma_1\sigma_2)$ (and
that $\psi(\sigma_3\sigma_2\sigma_3)=
\psi(\sigma_2\sigma_3\sigma_2)$), 
\begin{align*}
\psi(\sigma_1)\cdot \psi(\sigma_2)\cdot\psi(\sigma_1) &= \left( \begin{pmatrix}
0\\
0
\end{pmatrix}, \widetilde{\sigma}\right) \cdot
\left( \begin{pmatrix}
1\\
0
\end{pmatrix}, \widetilde{\sigma}\right) \cdot
\left( \begin{pmatrix}
0\\
0
\end{pmatrix}, \widetilde{\sigma}\right)\\
&=\left( \begin{pmatrix}
0\\
0
\end{pmatrix}, \widetilde{\sigma}\right) \cdot
\left( \begin{pmatrix}
1\\
0
\end{pmatrix}, \widetilde{\sigma}^2\right)
=\left( \begin{pmatrix}
0\\
-1
\end{pmatrix}, \widetilde{\sigma}^3\right).
\end{align*}
Similarly,
\begin{align*}
\psi(\sigma_2)\cdot \psi(\sigma_1)\cdot\psi(\sigma_2) &= 
\left( \begin{pmatrix}
1\\
0
\end{pmatrix}, \widetilde{\sigma}\right) \cdot
\left( \begin{pmatrix}
0\\
0
\end{pmatrix}, \widetilde{\sigma}\right)\cdot 
\left( \begin{pmatrix}
1\\
0
\end{pmatrix}, \widetilde{\sigma}\right)\\
&=\left( \begin{pmatrix}
1\\
0
\end{pmatrix}, \widetilde{\sigma}\right) \cdot
\left( \begin{pmatrix}
0\\
-1
\end{pmatrix}, \widetilde{\sigma}^2\right)
=\left( \begin{pmatrix}
0\\
-1
\end{pmatrix}, \widetilde{\sigma}^3\right).
\end{align*}
As for the surface relation,
\begin{align*}
\psi(\sigma_1\sigma_2\sigma_3^2\sigma_2\sigma_1)&=
\left( \begin{pmatrix}
0\\
-1
\end{pmatrix}, \widetilde{\sigma}^2\right)\cdot
\left( \begin{pmatrix}
0\\
0
\end{pmatrix}, \widetilde{\sigma}^2\right)\cdot
\left( \begin{pmatrix}
1\\
0
\end{pmatrix}, \widetilde{\sigma}^2\right)\\
&= \left( \begin{pmatrix}
0\\
-1
\end{pmatrix}, \widetilde{\sigma}^2\right)\cdot \left( \begin{pmatrix}
-1\\
-1
\end{pmatrix}, \widetilde{\sigma}^4\right)\\
&=\left( \begin{pmatrix}
0\\
0
\end{pmatrix}, \widetilde{\sigma}^6\right)= \left( \begin{pmatrix}
0\\
0
\end{pmatrix}, 1\right)
\end{align*}
as required. Since $\psi(\sigma_1)=\left( \begin{pmatrix}
0\\
0
\end{pmatrix}, \widetilde{\sigma}\right)$, $\psi(\sigma_2\sigma_1^{-1})=\left( \begin{pmatrix}
1\\
0
\end{pmatrix}, 1\right)$ and $\psi([\sigma_1^{-1},\sigma_2])=\left( \begin{pmatrix}
0\\
1
\end{pmatrix}, 1\right)$, we see that $\psi$ is surjective.

Now let $H=(B_4(\St))^{(2)}$, and let
$\map{\delta}{\Z^2 \rtimes \Z_6}[\Z_6]$ denote the projection onto the
second factor. Since $\Z_6$ is Abelian, it follows that
$\delta(\psi(x))$ is trivial for all  $x\in (B_4(\St))^{(1)}$, so
$\psi(x)$ belongs to the $\Z^2$-factor. Hence
$H= \bigl[(B_4(\St))^{(1)}, (B_4(\St))^{(1)} \bigr]\subseteq
\ker{\psi}$, and thus $\psi$ factors through $A=B_4(\St)/H$, inducing
a (surjective) homomorphism
$\map{\widehat{\psi}}{B_4(\St)/H}[\Z^2\rtimes \Z_6]$. From the
following commutative diagram of short exact sequences,
\begin{equation*}
\xymatrix{%
1 \ar[r] & A=\Gamma_2(B_4(\St))/H \ar[r] \ar[d]_{\widehat{\psi} \bigr\rvert_A} & B_4(\St)/H \ar[r]^{\widehat{\alpha}} \ar[d]_{\widehat{\psi}} & \bnab{4} \ar[r] \ar@{=}[d] & 1\\
1 \ar[r] & \Z^2 \ar[r] & \Z^2\rtimes \Z_6 \ar[r]^-{\delta} & \Z_6 \ar[r] & 1,}
\end{equation*}
the surjectivity of $\widehat{\psi}$ implies that of 
$\map{\widehat{\psi}\bigr\rvert_A}{A}[\Z^2]$. But $A$ is an Abelian
group generated by  $\brak{t_2,\overline{\sigma}t_2
\overline{\sigma}^{-1}}$, so $\widehat{\psi} \bigr\rvert_A$ is an
isomorphism, and by the $5$-Lemma, $\widehat{\psi}$ is too. Hence:
\begin{equation*}
(B_4(\St))^{(1)}\bigl/ (B_4(\St))^{(2)} \cong \Z^2\bigr. \quad\text{and}\quad B_4(\St)\bigl/(B_4(\St))^{(2)}\bigr.  \cong \Z^2\rtimes \Z_6.
\end{equation*}
In fact the first of these two equations may be obtained directly
since we know that $(B_4(\St))^{(1)}\cong \mathcal{Q}_8 \rtimes
\F[2]$, and $(B_4(\St))^{(2)}$ is isomorphic to the subgroup
$\mathcal{Q}_8 \rtimes (\F[2])^{(1)}$ of $\mathcal{Q}_8 \rtimes \F[2]$,
so $(B_4(\St))^{(1)}/(B_4(\St))^{(2)}\cong \F[2]/(\F[2])^{(1)}\cong
\Z^2$. Similarly, $(B_4(\St))^{(2)}/(B_4(\St))^{(3)}\cong (\Z_2 \times
\Z_2)\times (\F[2])^{(1)}/(\F[2])^{(2)}$,
$(B_4(\St))^{(3)}/(B_4(\St))^{(4)} \cong \Z_2\times
(\F[2])^{(2)}/(\F[2])^{(3)}$, and for $m\geq 4$, 
\begin{equation*}
(B_4(\St))^{(m)}\bigl/(B_4(\St))^{(m+1)}\bigr. \cong (\F[2])^{(m-1)}\bigl/(\F[2])^{(m)}.\bigr.
\end{equation*}
This proves the remaining parts of~(\ref{it:ds3}), and thus completes the proof of \reth{dsbn}.
\end{proof}

\chapter[Lower central and derived series of $\bmmn{m}{n}$]{The lower central and derived series of $\bmmn{m}{n}$}\label{chap:lcsbmn}

In this chapter, the aim is to determine the lower central and derived series of the $m$-string braid group of the $n$-punctured sphere $\bmmn{m}{n}$, $n\geq 1$ according to the values of $m$ and $n$. In \resec{presbmn}, we begin by giving a presentation of this group. In \resec{gorinlin}, we deal with the case $n=1$ which corresponds to the Artin braid groups, and extend the results of Gorin and Lin. The case $m=1$  which is that of the fundamental group of the $n$-punctured sphere is dealt with in \resec{bmn1n}. From \resec{m3n2} onwards, we suppose that $n\geq 2$. In \resec{m3n2}, we prove \reth{lcdsbmsn}, which if $m\geq 3$ (respectively $m\geq 5$) shows that the lower central series (respectively the derived series) of $\bmmn{m}{n}$ is constant from the commutator subgroup onwards. In Sections~\ref{sec:mgeq2},~\ref{sec:lcgsb22} and~\ref{sec:affineatil}, we study the case $n=2$ which corresponds to that of the braid groups of the annulus (which are isomorphic to the Artin groups of type~$B$). The main results of these three sections are \repr{aug}, \reco{b2b2}, \repr{b2b2resid}, \reth{lcsb2}
and \reco{dsannu}. In \resec{bmn3}, we study $B_m(\St\setminus\brak{x_1,x_2,x_3})$, $m\geq 2$, which is isomorphic to the affine Artin group of type $\widetilde{C}_m$, and we prove Propositions~\ref{prop:b2n3} and~\ref{prop:bmn3i}.

\section{A presentation of $\bmmn{m}{n}$, $n\geq 1$}\label{sec:presbmn}

Let $q\in\N$. If $1\leq i<j\leq q$, let
$A_{i,j}=\sigma_{j-1}\cdots\sigma_{i+1} \sigma_i^2\sigma_{i+1}^{-1}
\cdots \sigma_{j-1}^{-1}\in P_q(\St)$ which geometrically corresponds
to a twist of the $j\up{th}$ string about the $i\up{th}$ string, with
all other strings remaining vertical. It is well known that the
$A_{i,j}$ generate $P_q(\St)$.

The following presentation of $\bmmn{m}{n}$ was derived in~\cite{GG4} using standard results concerning presentations of group extensions~\cite{J} (see also~\cite{Lam,Ma} for other presentations).
\begin{prop}[\cite{GG4}]\label{prop:presbetanm}
Let $m\geq 1$ and $n\geq 1$. The following constitutes a presentation of the group $B_m(\St\setminus\brak{x_1,\ldots,x_n})$:
\begin{enumerate}
\item[\underline{\textbf{generators:}}] $A_{i,j}$, where $1\leq i\leq n$ and $n+1\leq j\leq n+m$, and $\sigma_k$, $1\leq k\leq m-1$.
\item[\underline{\textbf{relations:}}] for $1\leq i,k\leq n$, $n+1\leq j<l\leq n+m$ but $j\leq n+m$ if $l$ is absent, and $1\leq r,s\leq m-1$,
\begin{gather*}
\text{$A_{i,j}A_{k,l}A_{i,j}^{-1}=A_{k,l}$ if $k<i$}\\
A_{i,j}A_{i,l}A_{i,j}^{-1}=A_{j,l}^{-1}A_{i,l}A_{j,l}\\
A_{i,j}^{-1}A_{i,l}A_{i,j}= A_{i,l}A_{j,l}A_{i,l}A_{j,l}^{-1}A_{i,l}^{-1}\\
\text{$A_{i,j}A_{k,l}A_{i,j}^{-1}=A_{j,l}^{-1}A_{i,l}^{-1}A_{j,l} A_{i,l} A_{k,l}A_{i,l}^{-1}A_{j,l}^{-1}A_{i,l} A_{j,l}$ if $i<k$}\\
\text{$A_{i,j}^{-1}A_{k,l}A_{i,j}= A_{i,l}A_{j,l}A_{i,l}^{-1}A_{j,l}^{-1}A_{k,l} A_{j,l} A_{i,l}A_{j,l}^{-1}A_{i,l}^{-1}$ if $i<k$}\\
A_{1,n+m}\cdots A_{n,n+m}\sigma_{m-1}\cdots\sigma_2\sigma_1^2 \sigma_2\cdots \sigma_{m-1}=1\\
\text{$\sigma_r\sigma_{s}=\sigma_{s}\sigma_{r}$ if $\lvert r-s\rvert\geq 2$}\\
\text{$\sigma_{r}\sigma_{r+1}\sigma_{r}= \sigma_{r+1}\sigma_{r}\sigma_{r+1}$}\\
\text{$\sigma_r A_{i,j}\sigma_r^{-1}=A_{i,j}$ if $r\neq j-n-1,j-n$}\\ 
\text{$\sigma_{j-n} A_{i,j}\sigma_{j-n}^{-1}= A_{i,j+1}$ if $n+1\leq j\leq n+m-1$}.
\end{gather*}
In the above relations, if $n+1\leq j<l\leq n+m$ then $A_{j,l}$ (which does not appear in the list of generators) should be rewritten as:
\begin{equation*}
A_{j,l}=\sigma_{l-n-1}\ldots \sigma_{j-n+1}\sigma_{j-n}^2 \sigma_{j-n+1}^{-1}\ldots \sigma_{l-n-1}^{-1}. \tag*{\mbox{\qed}}
\end{equation*}
\end{enumerate}
\end{prop}

\begin{rems}\mbox{}
\begin{enumerate}[(a)]
\item Geometrically, we think of the $n$ punctures labelled as points
from $1$ to $n$, and the basepoints of the $m$ strings as points
labelled from $n+1$ to $n+m$. The generator $A_{i,j}$ corresponds
geometrically to a twist of the $(j-n)\up{th}$ string about the
$i\up{th}$ puncture, with all other strings remaining vertical.
\item This presentation was derived in~\cite{GG4} for $m\geq 1$ and $n\geq 3$ (see Proposition~9 of that paper). But it is also correct for $n=1,2$. Indeed, to obtain the result, a presentation of $P_m(\St\setminus\brak{x_1,\ldots,x_n})$ was derived (Proposition~7 of~\cite{GG4}) using the fact that there is a split short exact sequence
\begin{multline*}
1\to P_1(\St\setminus\brak{x_1,\ldots,x_n,x_{n+1}}) \to P_2(\St\setminus\brak{x_1,\ldots,x_n}) \to\\ P_1(\St\setminus\brak{x_1,\ldots,x_n}) \to 1,
\end{multline*}
which is the case for all $n\geq 1$ (as $\pi_2(\St\setminus\brak{x_1,\ldots,x_n})=\brak{1}$). To prove Proposition~9 of~\cite{GG4}, we then apply standard techniques to the short exact sequence 
\begin{equation*}
1\to P_m(\St\setminus\brak{x_1,\ldots,x_n})\to B_m(\St\setminus\brak{x_1,\ldots,x_n}) \to \sn[m]\to 1.
\end{equation*}
\end{enumerate}
\end{rems}

From this presentation, we may obtain easily the Abelianisation of  $\bmmn{m}{n}$:
\begin{prop}[\cite{GG4}, Proposition~11]\label{prop:abbm} The Abelianisation of $\bmmn{m}{n}$ is a free Abelian group of rank~$n$.\qed
\end{prop}

\section[The case $n=1$: series of Artin's braid groups]{The case $n=1$: lower central and derived series of Artin's braid groups $B_m(\dt)$}\label{sec:gorinlin}

As we shall see below, the case $n=1$ corresponds to that of Artin's braid groups. In \reth{gorinlin}, we recall Gorin and Lin's results,  which we extend in \repr{gorinlin}, notably obtaining descriptions of some of the derived series elements and quotients of $B_m(\dt)$ for $m=3,4$. We begin by proving the following proposition which will allow us to identify certain types of braid groups.

\begin{prop}\label{prop:iso}\mbox{}
\begin{enumerate}[(a)]
\item\label{it:isoa} Let $z_0\in \Int{\dt}$ and $m\geq 2$. Then $P_m(\dt)\cong P_{m-1}(\dt\setminus \brak{z_0})$.
\item\label{it:isob} Let $m\in\N$, let $x_0\in\St$, and let $Y\subseteq \St\setminus \brak{x_0}$ be a finite set. Then the inclusion $\St\setminus\brak{x_0}\subseteq \dt$ induces an isomorphism
\begin{equation*}
B_m\left(\St\setminus\left( Y\cup\brak{x_0}\right)\right)\cong B_m(\dt\setminus Y).
\end{equation*} 
\item\label{it:isoc} Let $z_0\in \Int{\dt}$ and $m\geq 1$. Then $B_{m,1}(\dt)\cong B_m(\dt\setminus \brak{z_0})$.
\item\label{it:isod} Let $m\in \N$ and $(x_1,\ldots,x_m)\in F_m(\Int{\dt})$. Then
\begin{equation*}
B_{m,1}(\dt)\cong \pi_1(\dt\setminus \brak{x_1,\ldots,x_m})\rtimes B_m(\dt).
\end{equation*}
\end{enumerate}
\end{prop}

\begin{rems}\mbox{}\label{rem:annulus}
\begin{enumerate}
\item Part~(\ref{it:isoa}) of \repr{iso} is a manifestation of the Artin combing operation~\cite{A2,Bi2,Han}: any geometric pure braid of the disc is equivalent to a pure braid whose first string is vertical.
\item\label{it:compact} Taking $Y=\vide$ in part~(\ref{it:isob}), and noting that homeomorphic spaces have isomorphic braid groups leads to the well-known isomorphism $B_m\cong B_m\left(\St\setminus\brak{x_0} \right)\cong B_m(\dt)$.
\item Part~(\ref{it:isoa}), and parts~(\ref{it:isoc}) and~(\ref{it:isod}) describe respectively the pure braid groups and full braid groups of the annulus. Since the latter are isomorphic to the Artin groups of type~$B$~\cite{Cr}, we recover part~(2) of Proposition~2.1 of~\cite{CP}.
\item In part~(\ref{it:isod}), the action is given by the well-known Artin representation of the Artin braid group as a subgroup of $\aut{\F[n]}$~\cite{A1,Bi2,Han}, and may be described as follows:
let $\sigma_1,\ldots, \sigma_{m-1}$ denote the standard generators of $B_m(\dt)$, and let $A_1, \ldots, A_m$ denote those of $\pi_1(\dt\setminus \brak{x_1,\ldots,x_m}, x_{m+1})$. Then:
\begin{equation*}
\sigma_i A_j\sigma_i^{-1}=
\begin{cases}
A_{i+1} &\text{if $j=i$}\\
A_{i+1}^{-1} A_i A_{i+1} &\text{if $j=i+1$}\\
A_j &\text{otherwise.}\\
\end{cases}
\end{equation*}
This was used by Chow~\cite{Ch,Han} to obtain a presentation of Artin's pure braid group, and may be applied to the study of the Nielsen equivalence problem for fixed points of surface homeomorphisms~\cite{Gu}.
\end{enumerate}
\end{rems}

\begin{proof}[Proof of \repr{iso}]\mbox{}
\begin{enumerate}[(a)]
\item Consider the following Fadell-Neuwirth short exact sequence for the disc:
\begin{equation*}
\xymatrix{%
1\ar[r] & P_{m-1}(\dt\setminus \brak{z_0}) \ar[r] & P_m(\dt)\ar[r] & P_1(\dt)\ar[r] & 1.}
\end{equation*}
Since $P_1(\dt)$ is trivial, it follows that the kernel is equal to $P_m(\dt)$, and the result follows.
\item Let $m\in\N$ and $(x_1,\ldots, x_m)\in F_m(\St\setminus(Y\cup \brak{x_0}))$. Set $X=\brak{x_1, \ldots,x_m}$ and $Y'=Y \cup\brak{x_0}$. The inclusion $\St\setminus\brak{x_0}\subseteq \dt$ induces an isomorphism of the free groups $\pi_1(\St\setminus(X\cup Y\cup\brak{x_0}))$ and $\pi_1(\dt\setminus(X\cup Y))$. Consider the following commutative diagram of short exact sequences:
\begin{equation*}
\vcenter{\xymatrix@C=0.55cm{%
1\ar[r] & \pi_1(\St\setminus (X\cup Y')) \ar[r] \ar[d]_{\cong} & P_{m+1}(\St\setminus Y') \ar[r] \ar[d] & P_m(\St\setminus Y') \ar[r] \ar[d] & 1\\
1\ar[r] & \pi_1(\dt\setminus(X\cup Y)) \ar[r] & P_{m+1}(\dt\setminus Y) \ar[r] &  P_m(\dt\setminus Y) \ar[r] & 1.}}
\end{equation*}
Applying induction on~$m$ and the $5$-Lemma, it follows that $P_m(\St\setminus Y') \cong P_{m+1}(\dt\setminus Y)$. By commutativity of the following diagram of short exact sequences
\begin{equation*}
\xymatrix{%
1\ar[r] & P_m(\St\setminus Y') \txt{~~}\ar@{^{(}->}[r] \ar[d]_{\cong} & B_m(\St\setminus Y') \ar[r] \ar[d] & \sn[m] \ar[r] \ar@{=}[d] & 1\\
1\ar[r] & P_m(\dt\setminus Y) \txt{~~}\ar@{^{(}->}[r] & B_m(\dt\setminus Y) \ar[r] &  \sn[m] \ar[r] & 1,}
\end{equation*} 
and the $5$-Lemma, we see that $B_m(\St\setminus Y')\cong B_m(\dt\setminus Y)$, which proves part~(\ref{it:isob}).
\item From the generalised Fadell-Neuwirth short exact sequence, we have that:
\begin{equation*}
\xymatrix{%
1\ar[r] & B_m(\dt\setminus\brak{z_0}) \ar[r] & B_{m,1}(\dt) \ar[r] & B_1(\dt) \ar[r] & 1.}
\end{equation*}
The result then follows easily.
\item Consider the following generalised Fadell-Neuwirth short exact sequence:
\begin{equation*}
\xymatrix@C=0.75cm{%
1\ar[r] & \pi_1(\dt\setminus\brak{x_1,\ldots,x_m})\ar[r] & B_{m,1}(\dt) \ar^{p_{\ast}}[r] & B_m(\dt) \ar[r] & 1.}
\end{equation*}
Since $p_{\ast}$ admits a section given by the obvious inclusion $B_m(\dt)\to B_{m,1}(\dt)$, the result again follows easily.\qedhere
\end{enumerate}
\end{proof}

So by \rerem{annulus}(\ref{it:compact}), $B_m(\St\setminus\brak{x_1})$
and $B_m(\dt)$ may be identified with Artin's braid group $B_m$. The
series of such groups were previously studied by Gorin and
Lin~\cite{GL}. For all $m\geq 1$, they determined presentations for
$\Gamma_2(B_m(\dt))$, from which they were able to deduce that:
\begin{thm}[\cite{GL}]\mbox{}\label{th:gorinlin}
\begin{enumerate}[(a)]
\item The commutator subgroups $\Gamma_2(B_m(\dt))$ are finitely presented.
\item\label{it:gorinlin3} $\Gamma_2(B_3(\dt))$ is a free group $\F[2](u,v)$ on two generators $u$ and $v$.
\item\label{it:gorinlin4} $\Gamma_2(B_4(\dt))$ is a semi-direct product of the free group $\F[2](a,b)$ by $\F[2](u,v)$, the action (denoted by $\phi$) being given by: 
\begin{equation}\label{eq:gorinlin4}
\left.\begin{aligned}
&\phi(u)(a)=uau^{-1}=b & &\phi(u)(b)=ubu^{-1}=b^2a^{-1}b\\
&\phi(v)(a)=vav^{-1}=a^{-1}b & &\phi(v)(b)=vbv^{-1}=(a^{-1}b)^3 a^{-2}b.
\end{aligned}\right\}
\end{equation}
\item For all $m\geq 5$, the derived subgroup $(B_m(\dt))^{(1)}$ is perfect, i.e.\ $(B_m(\dt))^{(s)}=(B_m(\dt))^{(1)}$ for all $s\geq 2$. 
\end{enumerate}
\end{thm}
We now go on to extend their results.

\begin{varthm}[\repr{gorinlin}]
Let $m\geq 1$. Then:
\begin{enumerate}[(a)]
\item For all $s\geq 3$, $\Gamma_s(B_m(\dt))= \Gamma_2(B_m(\dt))$.
\item If $m=1,2$ then $(B_m(\dt))^{(s)}=\brak{1}$ for all $s\geq 1$.
\item\label{it:gorin3}  If $m=3$ then the derived series of $(B_3(\dt))^{(1)}$ is that of the free group $\F[2](u,v)$ on two generators $u$ and $v$, where $u=\sigma_2\sigma_1^{-1}$ and $v=\sigma_1 u\sigma_1^{-1}= \sigma_1 \sigma_2\sigma_1^{-2}$.
Further,
\begin{equation*}
B_3(\dt)/(B_3(\dt))^{(2)}\cong \Z^2\rtimes \Z,
\end{equation*}
where $\Z^2$ is the free Abelian group generated by the respective Abelianisations $\overline{u}$ and $\overline{v}$ of $u$ and $v$, and the action is given by $\sigma \cdot \overline{u}=\overline{v}$ and $\sigma \cdot \overline{v}=-\overline{u}+\overline{v}$, where $\sigma$ is a generator of $\Z$.
\item If $m=4$ then 
\begin{gather*}
\text{$(B_4(\dt))^{(1)}/(B_4(\dt))^{(2)} \cong \Z^2$, and}\\
(B_4(\dt))^{(2)}\cong \F[2](a,b)\rtimes \Gamma_2(\F[2](u,v)),
\end{gather*}
where $a=\sigma_3\sigma_1^{-1}$ and $b= uau^{-1}= \sigma_2\sigma_3
\sigma_1^{-1} \sigma_2^{-1}$. 
\end{enumerate}
\end{varthm}

\begin{proof}[Proof of \repr{gorinlin}]\mbox{}
\begin{enumerate}[(a)]
\item The result follows from \relem{stallings}, since $\gpab[(B_m(\dt))]\cong \Z$, and $H_2(\Z)=\brak{1}$.
\item Clear.
\item The first part is a direct consequence of \reth{gorinlin}(\ref{it:gorinlin3}).
For the second part, the short exact sequence
\begin{equation*}
1\to (B_3(\dt))^{(1)}\to B_3(\dt) \to \Z \to 1
\end{equation*}
splits, where $(B_3(\dt))^{(1)}\cong \F[2](u,v)$ and $\Z\cong \ang{\sigma}=\gpab[(B_3(\dt))]$, and a section is given by sending $\sigma$ onto $\sigma_1$. So 
\begin{equation*}
B_3(\dt) \cong (B_3(\dt))^{(1)}\rtimes \Z,
\end{equation*}
where the action is given by $\sigma \cdot u= \sigma_1
\cdot \sigma_2\sigma_1^{-1}\cdot  \sigma_1^{-1}=
v$ and $\sigma \cdot v= \sigma_1 \cdot \sigma_1\sigma_2\sigma_1^{-2}
\cdot \sigma_1^{-1}= \sigma_1 \sigma_2^{-1}\cdot \sigma_1 \sigma_2
\sigma_1^{-2}=u^{-1} v$. Then $\gpab[\left((B_3(\dt))^{(1)}\right)]=
(B_3(\dt))^{(1)}/(B_3(\dt))^{(2)}\cong \Z^2$ is a free Abelian group with basis $\brak{\widetilde{u},\widetilde{v}}$, and so it follows that
\begin{equation*}
B_3(\dt)/(B_3(\dt))^{(2)} \cong (B_3(\dt))^{(1)}/(B_3(\dt))^{(2)} \rtimes \Z\cong \Z^2\rtimes \Z,
\end{equation*}
with action given by $\sigma \cdot \overline{u}= \overline{v}$ and
$\sigma \cdot \overline{v}= -\overline{u}+ \overline{v}$ as
required.
\item From \reth{gorinlin}(\ref{it:gorinlin4}), we know that $(B_4(\dt))^{(1)}$ is a semi-direct product of the free group $\F[2](a,b)$ by $\F[2](u,v)$, where $a,b,u$ and $v$ are as defined in the statement of the proposition, and the action is given by \req{gorinlin4}. Under Abelianisation of $(B_4(\dt))^{(1)}$, we see that $a$ and $b$ are sent to the trivial element, and there are no other relations between $u$ and $v$ other than the fact that they commute. So 
\begin{equation*}
(B_4(\dt))^{(1)}/(B_4(\dt))^{(2)} \cong \Z^2.
\end{equation*}

To see that $(B_4(\dt))^{(2)}\cong \F[2](a,b)\rtimes \Gamma_2(\F[2](u,v))$, we apply 
\repr{gammasemi}. Since $(uau^{-1})a^{-1}=ba^{-1}$ and $(ubu^{-1})b^{-1}=b^2a^{-1}$, it follows that $a,b\in L$, where $L$ is the subgroup generated by $\Gamma_2(\F[2](a,b))$ and the normal subgroup $[\F[2](u,v), \F[2](a,b)]$ of $\F[2](a,b)$ generated by all elements of the form $(ghg^{-1})h^{-1}$, where $g\in \F[2](u,v)$ and $h\in \F[2](a,b)$. So $L=\F[2](a,b)$, and the result follows. \qedhere
\end{enumerate}
\end{proof}

Hence the lower central series of $B_m(\dt)$ is determined for all $m\in\N$, in particular, if $m\geq 3$ then $B_m(\dt)$ is not residually nilpotent; and the derived series of $B_m(\dt)$ is determined for all $m\neq 4$. In this case, it remains to determine the higher derived subgroups and their quotients. For the next step, by \repr{gammasemi}, 
\begin{equation*}
(B_4(\dt))^{(3)}= [(B_4(\dt))^{(2)}, (B_4(\dt))^{(2)}]\cong K\rtimes (\F[2](u,v))^{(2)},
\end{equation*}
where $K$ is the subgroup of $\F[2](a,b)$ generated by $\Gamma_2(\F[2](a,b))$ and the normal subgroup of $\F[2](a,b)$ generated by the elements of the form $\phi(g)(h)h^{-1}$, where $g\in (\F[2](u,v))^{(1)}$ and $h\in \F[2](a,b)$. 

Let $N$ be the normal subgroup of $\F[2](a,b)$ generated by $[a,b]$, $a^2$ and $b^2$.\label{page:free5} It may be interpreted as the kernel of the homomorphism $\map{\psi}{\F[2](a,b)}[\Z_2\times \Z_2]$ which to a word $w=w(a,b)$ associates the exponent sums modulo~$2$ of $w$ relative to $a$ and $b$ respectively. 

In order to determine $K$ we need to investigate the action of $[u,v]$ and its conjugates on $\F[2](a,b)$. One can check that $\phi(u^{-1})(a)= ab^{-1}a^2$, $\phi(u^{-1})(b)= a$, $\phi(v^{-1})(a)=ab^{-1}a^3$ and $\phi(v^{-1})(b)= ab^{-1}a^4$, then that:
\begin{equation}\label{eq:uvab}
\left.\begin{aligned}
\phi([u,v])(a)a^{-1} &= ab^{-2}(ab^{-1})^4 a^{-1}\\  \phi(u[u,v]u^{-1})(a)a^{-1} &= 
(ab^{-3}(ab^{-2})^4)^2 ab^{-2} a^{-1}\\
\phi([u,v])(b)b^{-1} &= ab^{-2}(ab^{-1})^5 b^{-1}\\ \phi(u[u,v]u^{-1})(b)b^{-1} &= ab^{-3}(ab^{-2})^5 b^{-1}.
\end{aligned}\right\}
\end{equation}
Clearly $K$ contains $[a,b]$. Further,  a calculation shows that the element
\begin{equation*}
(\phi([u,v])(a)a^{-1})^{-3}(\phi([u,v])(b)b^{-1})^2b^{-2}
\end{equation*} belongs to $\Gamma_2(\F[2](a,b))$, and so to $K$. Hence $b^2$ belongs to $K$ too. Considering the element
\begin{equation*}
(\phi([u,v])(a)a^{-1})^{-4}(\phi([u,v])(b)b^{-1})^3,
\end{equation*}
we infer similarly that $a^2\in K$. Now $K$ is normal in $\F[2](a,b)$ and contains $[a,b]$, $a^2$ and $b^2$, so it contains $N$. We claim that $N=K$. Since $\psi$ factors through Abelianisation, we see that $\Gamma_2(\F[2](a,b))\subseteq N$. Let $w\in \F[2](u,v)$. If $\eta\in \brak{a,b}$ then 
\begin{equation*}
\psi(\phi(w)(\eta^2))=\psi(w\eta^2w^{-1})= 2\psi(\phi(w)(\eta))= (0,0).
\end{equation*}
Also,
\begin{equation*}
\psi(\phi(w)([a,b]))= \psi([waw^{-1}, wbw^{-1}])= (0,0).
\end{equation*}
This implies that $\psi(w)(N)\subseteq N$, so $\psi(w)$ induces an endomorphism $\widetilde{\phi(w)}$ of $\Z_2\times \Z_2$ satisfying $\psi\circ \phi(w)= \widetilde{\phi(w)}\circ \psi$. The surjectivity of $\psi$ and $\phi(w)$ imply that $\widetilde{\phi(w)}$ is an automorphism. Furthermore, $\widetilde{\phi(w_1)}\circ \widetilde{\phi(w_2)}= \widetilde{\phi(w_1w_2)}$ for all $w_1,w_2\in \F[2](u,v)$. Using the above relations for $\phi([u,v])$, we see that $\widetilde{\phi([u,v])}= \id$. Hence for all $w\in \Gamma_2(\F[2](u,v))$, $\widetilde{\phi(w)}= \id$, so $\psi(\phi(w)(h)h^{-1})= (0,0)$ for all $h\in \F[2](a,b)$. This implies that $K\subseteq N$, which proves the claim.

Finally, we Abelianise $(B_4(\dt))^{(2)}\cong \F[2](a,b)\rtimes (\F[2](u,v))^{(1)}$. To the commutativity relations between $a,b,u$ and $v$, one needs to add the relators the Abelianisation of the relators $\phi(w)(h)h^{-1}$ where $w\in (\F[2](u,v))^{(1)}$ and $h\in \F[2](a,b)$, and in particular of relations~\reqref{uvab}, from which we obtain $a^2=b^2=1$. But these are the only extra relations: since $\phi(w)(h)h^{-1}\in\ker{\psi}$, it follows from that form of $N$ that the Abelianised relations are products of powers of $a^2$ and $b^2$. We thus obtain:
\begin{varthm}[\repr{b4disc2ab}]
\begin{equation*}
\left(B_4(\dt)\right)^{(2)}\left/\left(B_4(\dt)\right)^{(3)}\right. \cong \left( \Z_2\times \Z_2\right)\times
\gpab[{\left(\Gamma_2(\F[2](u,v))\right)}].\tag*{\mbox{\qed}}
\end{equation*}
\end{varthm}

Using the Reidemeister-Schreier rewriting process~\cite{MKS}, we may obtain a presentation of $N$. Let $X=\brak{a,b}$ be a generating set of $\F[2](a,b)$ and $U=\brak{1,a,ab,aba^{-1}}$ be a Schreier transversal. If $g\in \F[2](a,b)$, let $\overline{g}\in U$ denote its coset representative. A basis of $N$ is given by the set of elements of the form $ux(\overline{ux})^{-1}$ where $u\in U$ and $x\in X$ (we remove all occurrences of the trivial element). A simple calculation shows that $N$ is a free group of rank~$5$ with basis whose elements are given by $a^2$, $aba^2b^{-1}a^{-1}$, $bab^{-1}a^{-1}$, $ab^2a^{-1}$ and $b^2$. This may be transformed into the following basis: $z_1=a^2$, $z_2=b^2$, $z_3=(ab)^2$, $z_4=ba^2b^{-1}$ and $z_5=ab^2a^{-1}$. The action of $\F[2](u,v)$ on $N$ is given by equations~\reqref{acuxi} and~\reqref{acvxi}, see Tables~\ref{tab:uaction} and~\ref{tab:vaction} (we have used the relations $[a,b]=ababb^{-1}a^{-2}b^{-1}=z_3z_2^{-1}z_4^{-1}$, and $(ba^{-1})^2=z_2z_3^{-1}z_5$).
\addtocounter{table}{3}
\begin{table}
\begin{equation}\label{eq:acuxi}
\left.
\begin{aligned}
uz_1u^{-1} & = (uau^{-1})^2=z_2\\
uz_2u^{-1} & = (ubu^{-1})^2=b^2a^{-2}ab^2a^{-1}[a,b]b^2\\
&= z_2z_1^{-1}z_5 z_3z_2^{-1}z_4^{-1} z_2\\
uz_3u^{-1} & =(uabu^{-1})^2=b^4b^{-1}a^{-1}b^{-1}a^{-1} ab^4a^{-1}[a,b]b^2\\
&= z_2^2z_3^{-1} z_5^2 z_3z_2^{-1}z_4^{-1} z_2\\
uz_4u^{-1} & =b^2a^{-2}ab^2a^{-1}a^2b^{-2}= z_2z_1^{-1}z_5z_1z_2^{-1}\\
uz_5u^{-1} & =b^4b^{-1}a^{-1}b^{-1}a^{-1} ab^4a^{-1}=z_2^2z_3^{-1}z_5^2
\end{aligned}\right\}
\end{equation}
\medskip
\caption{\label{tab:uaction}The action of $u$ on the basis $z_1,\ldots, z_5$ of $N$}
\end{table}
\begin{table}
\begin{equation}\label{eq:acvxi}
\left.
\begin{aligned}
vz_1v^{-1} =& (vav^{-1})^2= (a^{-1}b)^2= a^{-2}(ab)^2b^{-2}ba^{-2}b^{-1}b^2\\
=& z_1^{-1}z_3z_2^{-1} z_4^{-1}z_2\\
vz_2v^{-1} =& a^{-2}(ab)^2 b^{-2}ba^{-2}b^{-1}b^2 a^{-2} (ab)^2b^{-2} ba^{-2}b^{-1} b^2\cdot\\ & (b^{-1}a^{-1})^2 ab^2a^{-1} b^2 (b^{-1}a^{-1})^2 ab^2a^{-1} ba^{-2}b^{-1} b^2\\
=& z_1^{-1}z_3 z_2^{-1} z_4^{-1}z_2 z_1^{-1} z_3z_2^{-1}z_4^{-1} z_2z_3^{-1}z_5 z_2z_3^{-1}z_5 z_4^{-1} z_2\\
=& (vz_1v^{-1})^2z_3^{-1}z_5 z_2z_3^{-1}z_5 z_4^{-1} z_2\\
vz_3v^{-1} =& (vabv^{-1})^2= (a^{-1}b)^4a^{-2}ba^{-1}ba^{-1}ba^{-1}b a^{-1}b a^{-2}b\\
=& (vz_1v^{-1})^2 z_1^{-1}(z_2z_3^{-1}z_5)^2z_4^{-1}z_2\\
=& (vz_1v^{-1})^2 z_1^{-1}z_2 vz_1^{-2}z_2v^{-1}\\
vz_4v^{-1} =& (a^{-1}b)^4 b^{-1}a^{-1}b^{-1}a^{-1}ab^2a^{-1}b^2 (a^{-1}b)^{-4}\\
=& (vz_1v^{-1})^2 z_3^{-1}z_5z_2 (vz_1v^{-1})^{-2}\\
vz_5v^{-1} =& (a^{-1}b)^4  a^{-2} b a^{-1}ba^{-1}ba^{-1}b  a^{-1}\\
=& (vz_1v^{-1})^2 z_1^{-1} (z_2z_3^{-1}z_5)^2
\end{aligned}\right\}
\end{equation}
\medskip
\caption{\label{tab:vaction}The action of $v$ on the basis $z_1,\ldots, z_5$ of $N$}
\end{table}
Hence:
\begin{varthm}[\repr{b4deriv}]
$(B_4(\dt))^{(3)}\cong \F[5](z_1,\ldots,z_5) \rtimes (\F[2](u,v))^{(2)}$, where the action is that induced by the action of $\F[2](u,v)$ on $\F[5](x_1,\ldots,x_5)$ given by equations~\reqref{acuxi} and~\reqref{acvxi}.\qed
\end{varthm}

From this, we may determine the Abelianisation $\gpab[((B_4(\dt))^{(3)})]$ of $(B_4(\dt))^{(3)}$:
\begin{varthm}[\repr{b4derivab}]
\begin{align*}
\gpab[{\left((B_4(\dt))^{(3)}\right)}]  &=(B_4(\dt))^{(3)}\left/(B_4(\dt))^{(4)}\right.\\
& \cong \Z^3 \times \Z_{18}\times \Z_{18} \times (\F[2](u,v))^{(2)}\left/(\F[2](u,v))^{(3)}.\right.
\end{align*}
\end{varthm}

\begin{proof}
The action of $\F[2](u,v)$ on $\F[5](z_1,\ldots,z_5)$ is by conjugation which leaves  $\Gamma_2(\F[5](z_1,\ldots,z_5))$ invariant. It thus induces an action of $\F[2](u,v)$ on 
$\gpab[{\F[5](z_1,\ldots,z_5)}]= \Z^5= \Z^5[Z_1,\ldots,Z_5]$,
where for $i=1,\ldots,5$, $Z_i$ is the image of $z_i$ under Abelianisation. For $w\in \F[2](u,v)$, let $M_w$ denote the matrix of this action with respect to the basis $(Z_1,\ldots, Z_5)$ of $\Z^5$, and let 
\begin{equation*}
\Lambda_w=\im{M_w-I_5}.
\end{equation*}
By equations~\reqref{acuxi} and~\reqref{acvxi},
\begin{equation*}
U=M_u=
\left( \begin{smallmatrix}
0 & -1 & 0 & 0 & 0\\
1 & 1 & 2 & 0 & 2\\
0 & 1 & 0 & 0 & -1\\
0 & -1 & -1 & 0 & 0\\
0 & 1 & 2 & 1 & 2
\end{smallmatrix}\right),\quad
U^{-1}=M_{u^{-1}}=
\left( \begin{smallmatrix}
1  &   1  &   2  &   2   &  0\\
-1 &    0  &   0 &    0  &   0\\
1  &   0  &   0  &  -1  &   0\\
1  &   0  &   2  &   2 &    1\\
-1  &   0  &  -1  &   0  &   0
\end{smallmatrix}\right),
\end{equation*}
and
\begin{equation*}
V=M_v=\left( \begin{smallmatrix}
-1 & -2 & -3 & 0 & -3\\
0 & 2 & 3 & 1 & 2\\
1 & 0 & 0 & -1 & 0\\
-1 & -3 & -3 & 0 & -2\\
0 & 2 & 2 & 1 & 2
\end{smallmatrix}\right),\quad
V^{-1}=M_{v^{-1}}=
\left( \begin{smallmatrix}
2  &   2   &  5  &   2  &   3\\
0  &  -1  &  -1  &  -1  &   0\\
0   &  1  &   0  &   0  &  -1\\
2  &   2  &   4   &  2   &  3\\
-1  &  -1  &  -1  &   0  &   0
\end{smallmatrix}\right).
\end{equation*}
Then
\begin{equation*}
C= M_{[u,v]}=\left( \begin{smallmatrix}
3 & 3 & 5 & 2 & 3\\
-3 & -3 & -7 & - 3 & -4\\
0 & 0 & 1 & 0 & 0\\
2 & 3 & 5 & 3 & 3\\
-3 & -4 & -7 & - 3 & -3
\end{smallmatrix}\right), \; 
C^{-1}= M_{[u,v]^{-1}}=\left( \begin{smallmatrix}
-3 & -3 & -7 & -4 & -3\\
3 & 3 & 5 & 3 & 2\\
0 & 0 & 1 & 0 & 0\\
-4 & -3 & -7 & - 3 & -3\\
3 & 2 & 5 & 3 & 3
\end{smallmatrix}\right).
\end{equation*}

Let $L$ be the subgroup of $\Z^5$ generated by the $\Lambda_w$, where $w\in (\F[2](u,v))^{(2)}$. The action of $\F[2](u,v)$ on $\Z^5$ restricts to an action of $(\F[2](u,v))^{(2)}$ on $\Z^5$. Since $L$ is generated by the relators $M_w(Z)-Z$, where $Z\in \Z^5$ and $w\in (\F[2](u,v))^{(2)}$, it follows that
\begin{equation}\label{eq:gogu}
\gpab[{\left((B_4(\dt))^{(3)}\right)}] \cong \Z^5\big/ L \times \left(\F[2](u,v)\right)^{(2)}\big/\left(\F[2](u,v)\right)^{(3)}.
\end{equation}
Let $c=[u,v]$ and $a=u[u,v]u^{-1}$. Consider first
the special case $w=[a,c]$, and set $\Sigma=\Lambda_{[a,c]}$. We claim
that:
\begin{enumerate}[(i)]
\item\label{it:claim2} $\Z^5/\Sigma\cong \Z^3\times \Z_{18}\times \Z_{18}$, and
\item\label{it:claim1} $L=\Sigma$.
\end{enumerate}
From this, it is obvious that $\Z^5/L\cong \Z^3\times \Z_{18}\times \Z_{18}$, and so the result follows from \req{gogu}.

To prove claim~(\ref{it:claim2}), one may check that 
\begin{equation*}
M_{[a,c]}= \left( \begin{smallmatrix}
-701 & -612 & -1314 & -702 & -612\\
1548 & 1351 & 2898 & 1548 & 1350\\
0 & 0 & 1 & 0 & 0\\
-702 & -612 & -1314 & -701 & -612\\
1548 & 1350 & 2898 & 1548 & 1351
\end{smallmatrix}\right).
\end{equation*}
So $\Sigma$ is the free Abelian group of rank~$2$ freely generated by $A_1=\left(
\begin{smallmatrix}
-702\\
1548\\
0\\
-702\\
1548
\end{smallmatrix}
\right)$ and $A_2=\left(
\begin{smallmatrix}
-612\\
1350\\
0\\
-612\\
1350
\end{smallmatrix}
\right)$, and $\Z^5/\Sigma$ has a finite presentation
\begin{equation*}
0\to \Sigma \stackrel{T}{\longrightarrow} \Z^5\to \Z^5/\Sigma \to 0,
\end{equation*}
where $T$ is the $\Z$-module homomorphism represented by the matrix $A=\left(
\begin{smallmatrix}
-702 & -612\\
1548 & 1350\\
0 & 0\\
-702 & -612\\
1548 & 1350
\end{smallmatrix}\right)$ relative to the bases $(A_1,A_2)$ and $(Z_1,\ldots, Z_5)$. Applying elementary row and column operations to $A$, and taking $P=\left(
\begin{smallmatrix}
2 & 1 & 0 & 0 & 0\\
-11 & -5 & 0 & 0 & 0\\
0 & 0 & 1 & 0 & 0\\
-1 & 0 & 0 & 1 & 0\\
0 & -1 & 0 & 0 & 1
\end{smallmatrix}\right)$ and $Q=\left(
\begin{smallmatrix}
1 & 7\\
-1 & -8
\end{smallmatrix}\right)$, we see that $PAQ=
\left(
\begin{smallmatrix}
18 & 0\\
0 & 18\\
0 & 0\\
0 & 0\\
0 & 0
\end{smallmatrix}\right)$, which gives the invariant factors of the Smith normal form of $A$~\cite{AW}. A new basis $W_1,\ldots, W_5$ of $\Z^5$ is obtained by taking $W_j=\sum_{i=1}^5\; (P^{-1})_{i,j} Z_i$, so $W_1=-5Z_1+11Z_2-5Z_4+11Z_5$, $W_2=-Z_1+2Z_2- Z_4+2Z_5$, $W_3=Z_3$, $W_4=Z_4$ and $W_5=Z_5$, and from the form of $PAQ$, it follows that in $\Z^5/\Sigma$, $18W_1=18W_2=0$, and that $W_3,W_4$ and $W_5$ are free generators. Thus $\Z^5/\Sigma\cong\Z^3 \times \Z_{18} \times \Z_{18}$, which proves claim~(\ref{it:claim2}).

We now set about proving claim~(\ref{it:claim1}). Since $[a,c]\in (\F[2](u,v))^{(2)}$, it is clear that $\Sigma\subseteq L$. For the converse, it suffices to check that for all $w\in (\F[2](u,v))^{(2)}$, $\Lambda_w= \im{M_w-I_5}\subseteq \ang{A_1,A_2}$. 

First note that $u$ and $v$ induce automorphisms of $\Sigma$; indeed, one may check that relative to the basis $(A_1,A_2)$, the matrix of $u$ is $\left(
\begin{smallmatrix}
-996 & -869\\
1145 & 999
\end{smallmatrix}\right)$, and that of $v$ is $\left(
\begin{smallmatrix}
18955 & 16531\\
-21731 & -18952
\end{smallmatrix}\right)$. So $M_w(\Sigma)=\Sigma$ for all $w\in \F[2](u,v)$.

Further, since for all $w\in \F[2](u,v)$, $M_w$ is an automorphism of $\Z^5$ which leaves $\Sigma$ invariant, if $y\in \F[2](u,v)$ satisfies $\Lambda_y=\im{M_y-I_5}\subseteq \Sigma$ then it follows that
\begin{align*}
\Lambda_{wyw^{-1}}=\im{M_{wyw^{-1}}-I_5}&=\im{M_w(M_y-I_5)M_{w^{-1}}}\\
&= \im{M_w(M_y-I_5)} \subseteq M_w(\Sigma)= \Sigma.
\end{align*}
So for our purposes, it will suffice to consider elements of $(\F[2](u,v))^{(2)}$ modulo conjugation by elements of $\F[2](u,v)$. Moreover, if $w_1,w_2\in \F[2](u,v)$  satisfy $\im{M_{w_i}-I_5}\subseteq \Sigma$ for $i=1,2$, then $\im{M_{w_1w_2}-I_5}\subseteq \Sigma$. This follows from the fact that for all $x\in\Z^5$, 
\begin{equation*}
(M_{w_1w_2}-I_5)(x)=M_{w_1}(M_{w_2}-I_5)(x)+(M_{w_1}-I_5)(x),
\end{equation*}
and the invariance of $\Sigma$ under $M_w$.

We now give a generating set for $(\F[2](u,v))^{(2)}$. A generating
set of $(\F[2](u,v))^{(1)}$ is given by the set of conjugates $wc^{\pm
1}w^{-1}$, where $w\in \F[2](u,v)$ (cf.\ \rerem{freederive}), and so
$(\F[2](u,v))^{(2)}$ is generated, up to conjugacy, by the set of 
commutators of the $wc^{\pm 1}w^{-1}$. So up to conjugacy, 
$(\F[2](u,v))^{(2)}$ is generated by the set of elements of the form
$\left[ c^{\epsilon_1}, tc^{\epsilon_2}t^{-1}\right]$, where
$\epsilon_1, \epsilon_2\in \brak{1,-1}$. By conjugating by
$c^{-1}tc^{\epsilon_2}t^{-1}$ if necessary, we may suppose that
$\epsilon_1=1$. By the remarks of the previous paragraph, it thus
suffices to show that $\im{M_y-I_5}\subseteq \ang{A_1,A_2}$ for
elements $y$ of the form $\left[ c, tc^{\epsilon_2}t^{-1}\right]$,
where $\epsilon_2\in \brak{1,-1}$. In order to do this, we shall now
calculate $M_y-I_5$ explicitly.

\begin{lem}\label{lem:gogustruc}
For all $t \in \F[2](u,v)$, $M_{tc^{\pm 1}t^{-1}}$ is of the form:
\begin{equation*}
A=\left( \begin{smallmatrix}
3m & 3n & 3m + 3n - 1 & 3m - 1 & 3n \\
 - 3p &  - 3m &  - 3m - 3p - 1 &  - 3p &  - 3m - 1 \\
0 & 0 & 1 & 0 & 0 \\
3m - 1 & 3n & 3m + 3n - 1 & 3m & 3n \\
 - 3p &  - 3m - 1 &  - 3m - 3p - 1 &  - 3p &  - 3m
\end{smallmatrix}\right),
\end{equation*}
where $m,n,p\in\Z$ and $np=m^2$.
\end{lem}

\begin{rem}\label{rem:goguinverse}
One may check easily that the inverse of this matrix is:
\begin{equation*}
A^{-1}=\left( \begin{smallmatrix}
- 3 m &  - 3 n &  - 3 m - 3 n - 1 &  - 3 m - 1 &  - 3 n \\
3 p & 3 m & 3 m + 3 p - 1 & 3 p & 3 m - 1 \\
0 & 0 & 1 & 0 & 0 \\
 - 3 m - 1 &  - 3 n &  - 3 m - 3 n - 1 &  - 3 m &  - 3 n \\
3 p & 3 m - 1 & 3 m + 3 p - 1 & 3 p & 3 m
\end{smallmatrix}\right).
\end{equation*}
So if $A$ satisfies the conditions of \relem{gogustruc} then it is inversible, and $A^{-1}$ also satisfies the conditions. Notice that $A$ may be obtained simply from $A^{-1}$ via the symmetry $(m,n,p)\mapsto (-m,-n,-p)$.
\end{rem}

\begin{proof}[Proof of \relem{gogustruc}.] We proceed by induction on the length $\ell(t)$ of the word $t$. If $\ell(t)=0$ then $t$ is the trivial element, and clearly $C=M_c$ and $C^{-1}=M_{c^{-1}}$ have the given structure. So suppose that $t$ has word length $\ell(t)\geq 0$, and that $M_{tc^{\pm 1}t^{-1}}$ has the given structure. By \rerem{goguinverse}, it suffices to prove the result for $M_{tct^{-1}}$. Setting $A=M_{tct^{-1}}$, a long but straightforward calculation shows that the respective conjugates of $M_{tct^{-1}}$ by $M_u$, $M_{u^{-1}}$, $M_v$ and $M_{v^{-1}}$ are:
\begin{gather*}
\left( \begin{smallmatrix}
9p - 3m & 3p & 12p - 1 - 3m & 9p - 1 - 3m
 & 3p \\
18m - 27p - 3n & 3m - 9p & 21m - 36p - 1 - 3
n & 18m - 27p - 3n & 3m - 1 - 9p \\
0 & 0 & 1 & 0 & 0 \\
9p - 1 - 3m & 3p & 12p - 1 - 3m & 9p - 3m
 & 3p \\
18m - 27p - 3n & 3m - 1 - 9p & 21m - 36p - 1
 - 3n & 18m - 27p - 3n & 3m - 9p
\end{smallmatrix}\right),\\
\left( \begin{smallmatrix}
9n - 3m&  - 18m + 3p + 27n &  - 1 + 36n - 21
m + 3p& 9n - 3m - 1&  
 - 18m + 3p + 27n \\
 - 3n& 3m - 9n&  - 1 - 12n + 3m&  - 3n
&  - 1 - 9n + 3m \\
0& 0& 1& 0& 0 \\
9n - 3m - 1&  - 18m + 3p + 27n &  - 1 + 36n
 - 21m + 3p& 9n - 3m&  
 - 18m + 3p + 27n \\
 - 3n&  - 1 - 9n + 3m&  - 1 - 12n + 3m& 
 - 3n& 3m - 9n
\end{smallmatrix}\right),\\
\left( \begin{smallmatrix}
\gamma_1 &  - 30m + 75p + 3n&  - 57
m + 135p + 6n - 1&  
 \gamma_1 - 1&  - 30m + 75p + 3n \\
 - 48p + 24m - 3n&  - \gamma_1 &  - 1
 - 108p + 51m - 6n&  
 - 48p + 24m - 3n&  - 1 - \gamma_1 \\
0& 0& 1& 0& 0 \\
 \gamma_1 - 1&  - 30m + 75p + 3n& 
 - 57m + 135p + 6n - 1&  
 \gamma_1 &  - 30m + 75p + 3n \\
 - 48p + 24m - 3n&  - 1 - \gamma_1& 
 - 1 - 108p + 51m - 6n&  
 - 48p + 24m - 3n&  - \gamma_1
\end{smallmatrix}\right),\\
\intertext{where $\gamma_1= - 27m + 60p + 3n$, and}
\left( \begin{smallmatrix}
\gamma_2 & -120m + 75p + 48n &  - 147m + 90
p - 1 + 60n&  - 1 + \gamma_2&  
-120m + 75p + 48n \\
 - 3p + 6m - 3n&  - \gamma_2 &  - 1 - 
18p + 33m - 15n&  
 - 3p + 6m - 3n&  - 1 - \gamma_2 \\
0& 0& 1& 0& 0 \\
 - 1 + \gamma_2 & -120m + 75p + 48n &  - 147m + 
90p - 1 + 60n&  \gamma_2 &  
-120m + 75p + 48n \\
 - 3p + 6m - 3n&  - 1 - \gamma_2 &  - 
1 - 18p + 33m - 15n&  
 - 3p + 6m - 3n&  - \gamma_2
\end{smallmatrix}\right),
\end{gather*}
where $\gamma_2=- 27m + 15p + 12n$.
One may then check that each of these matrices has the form of the statement of the lemma.
\end{proof}

We first consider the case $\epsilon_2=1$, so $y=\left[ c, tct^{-1}\right]$. With the matrix $M_{tct^{-1}}=A$ given by \relem{gogustruc}, a long but straightforward calculation shows once more that $M_{[c,tct^{-1}]}-I_5$ is of the form 
\begin{equation*}
\left( \begin{smallmatrix}
\alpha_1 & \alpha_2 & \alpha_1+\alpha_2 & \alpha_1 & \alpha_2\\
\beta_1 & \beta_2 & \beta_1+\beta_2 & \beta_1 & \beta_2\\
0 & 0 & 0 & 0 & 0\\
\alpha_1 & \alpha_2 & \alpha_1+\alpha_2 & \alpha_1 & \alpha_2\\
\beta_1 & \beta_2 & \beta_1+\beta_2 & \beta_1 & \beta_2
\end{smallmatrix}\right),
\end{equation*}
where
\begin{align*}
\begin{split}
\alpha_1 =&1278m^2 + 216m - 1836pm - 126p - 540
nm\\
 &- 90n + 648p^2 + 450pn
\end{split}\\ 
\begin{split}
\alpha_2 =& 1728nm - 756pn - 540n^{2} - 1080m^{2}\\
 & + 648pm - 72n + 180m - 108p
 \end{split} \displaybreak[0] \\ 
\begin{split}
\beta_1 =& - 1512m^{2} - 252m + 2160pm + 144p\\
 & + 648nm + 108n - 756p^{2} - 540pn
\end{split}\\
\begin{split}
\beta_2 =& - 2052nm + 882pn + 648n^{2} + 1278m^{2}\\
 & - 756pm - 216m + 126p + 90n.
\end{split}
\end{align*}
So $\im{M_{[c,tct^{-1}]}-I_5}$ is generated by the first two columns $C_1,C_2$ of $M_{[c,tct^{-1}]}-I_5$. It is necessary to show that each belongs to $\ang{A_1,A_2}$, in other words, that for $i=1,2$, there exist $\tau_i,\mu_i \in\Z$ such that $\tau_i A_1+\mu_i A_2=C_i$; these equations are equivalent to $\tau_i \left( \begin{smallmatrix}
-702\\
1548
\end{smallmatrix}\right)+
\mu_i \left( \begin{smallmatrix}
-612\\
1350
\end{smallmatrix}\right) = 
\left( \begin{smallmatrix}
\alpha_i\\
\beta_i
\end{smallmatrix}\right)$, and admit solutions (\emph{a priori} rational) of the form
\begin{equation*}
\begin{pmatrix}
\tau_i\\
\mu_i
\end{pmatrix} =-\frac{1}{324}
\begin{pmatrix}
1350 & 612\\
-1548 & -702
\end{pmatrix} 
\begin{pmatrix}
\alpha_i\\
\beta_i
\end{pmatrix}.
\end{equation*}
Substituting for $\alpha_i,\beta_i$, we obtain
 \begin{align*}
\begin{split}
\tau_1 =& - 2469m^2 - 424m + 3570pm + 253p + 1026nm + 171
n\\
& - 1272p^2 - 855pn
\end{split}\\
\begin{split}
\mu_1 =& 2830m^2 + 486m - 4092pm - 290p - 1176nm - 196n\\
&
 + 1458p^2 + 980pn
 \end{split}\\
\begin{split}
\tau_2 =& - 3324nm + 1484pn + 1026n^2 + 2086m^2 - 1272p
m + 130n\\
& - 342m + 212p
\end{split}\\
\begin{split}
\mu_2 =& 3810nm - 1701pn - 1176n^2 - 2391m^2 + 1458pm
 - 149n\\
& + 392m - 243p,
\end{split}
\end{align*}
and these solutions are clearly integers. Hence $C_1,C_2 \in \ang{A_1,A_2}$ as required.

To deal with the case $\epsilon_2=-1$, it suffices to invoke the observation of \rerem{goguinverse} concerning the symmetry between $A$ and $A^{-1}$. The above analysis holds, and we obtain the same solutions as above, but replacing everywhere $m,n$ and $p$ by $-m,-n$ and $-p$ respectively. This proves claim~(\ref{it:claim1}), and completes the proof of \repr{b4derivab}.
\end{proof}

We would now like to go a stage further, and determine $(B_4(\dt))^{(4)}$ and/or its Abelianisation. By applying \repr{gammasemi} to \repr{b4deriv}, $(B_4(\dt))^{(4)}$ is isomorphic to $M\rtimes \left(\F[2](u,v)\right)^{(3)}$, where $M$ is the subgroup of $\F[5](z_i)=\F[5](z_1,\ldots,z_5)$ generated by $\Gamma_2\left(\F[5](z_i)\right)$ and the normal subgroup generated by the elements of the form $\phi(g)(h)h^{-1}$, where $g\in \left(\F[2](u,v)\right)^{(2)}$ and $h\in \F[5](z_i)$. However, the complexity of finding a basis of $(\F[2](u,v))^{(2)}$ and calculating the action on $\F[5](z_i)$ makes it extremely difficult to obtain a description of $M$. In order to get some idea of the situation, we shall turn our attention to studying the semi-direct product $\F[5](z_i)\rtimes \F[2](u,v)$. In any case, the calculations that follow shall be used later in \resec{mgeq2} in order to study $\Gamma_2(B_3(\St\setminus\brak{x_1, x_2}))$.

From relations~\reqref{acuxi} and~\reqref{acvxi}, we have an action of $\F[2](u,v)$ on $\F[5](z_1,\ldots,z_5)$, and thus a semi-direct product $\F[5](z_1,\ldots,z_5) \rtimes \F[2](u,v)$. Let
\begin{equation}\label{eq:defeps} 
\map{\epsilon}{\F[5](z_1,\ldots,z_5) \rtimes \F[2](u,v)}[{(\gpab[{\F[5]}(z_1,\ldots,z_5)\rtimes {\F[2]}(u,v))]}]
\end{equation}
be Abelianisation. From relations~\reqref{acuxi} and~\reqref{acvxi}, it follows that 
\begin{equation*}
\gpab[{\left(\F[5](z_1,\ldots,z_5) \rtimes \F[2](u,v)\right)}]\cong \Z\oplus\Z^2,
\end{equation*}
where $\epsilon(u)=(0,1,0)$, $\epsilon(v)=(0,0,1)$, $\epsilon(z_i)=(1,0,0)$ if $i=1,2,3$ and $\epsilon(z_i)=(-1,0,0)$ if $i=4,5$.

Let $(Z_1,\ldots,Z_5)$, $(W_1,\ldots,W_5)$ be the bases of $\Z^5$ given in the proof of \repr{b4derivab}. Let 
\begin{equation*}
\map{\widetilde{\epsilon}}{\F[5](z_1,\ldots,z_5) \rtimes \left(\F[2](u,v)\right)^{(2)}}[{\gpab[{\left(\F[5](z_1,\ldots,z_5) \rtimes \F[2](u,v)^{(2)}\right)}]}]
\end{equation*}
be the restriction of $\epsilon$ to $\F[5](z_1,\ldots,z_5) \rtimes (\F[2](u,v))^{(2)}$. We identify this latter group with $\Z_{18}\times \Z_{18} \times \Z^3\times (\F[2](u,v))^{(2)}/(\F[2](u,v))^{(3)}$ via \repr{b4derivab}. Since $Z_1=2W_1-11W_2-W_4$, $Z_2=W_1-5W_2-W_5$ and $W_i=Z_i$ for $i=3,4,5$, as an element of $\Z_{18}\times \Z_{18} \times \Z^3$, we have $\widetilde{\epsilon}(z_1)= (\overline{2}, -\overline{11}, 0,-1,0)$, $\widetilde{\epsilon}(z_2)= (\overline{1}, -\overline{5},0, 0,-1)$, $\widetilde{\epsilon}(z_3)= (\overline{0}, \overline{0},1,0,0)$,
$\widetilde{\epsilon}(z_4)= (\overline{0}, \overline{0},0,1,0)$ and $\widetilde{\epsilon}(z_5)= (\overline{0}, \overline{0},0,0,1)$.
We thus obtain the following commutative diagram of short exact sequences:
\begin{equation}\label{eq:b4exst}
\begin{xy}*!C\xybox{\xymatrix{%
1 \ar[r] & (B_4(\dt))^{(4)} \ar[r] \ar[d] & G \ar[r]^(0.25){\widetilde{\epsilon}} \ar[d] & \gpab[\left((B_4(\dt))^{(3)}\right)] \ar[r] \ar[d]_{\xi} & 1\\
1 \ar[r] & \ker{\epsilon} \ar[r] & H \ar[r]^(0.4){{\epsilon}} & \Z\oplus\Z^2 \ar[r] & 1,}}
\end{xy}
\end{equation}
where we set
\begin{equation*}
\text{$G=\F[5](z_1,\ldots,z_5) \rtimes (\F[2](u,v))^{(2)}$ and $H=\F[5](z_1,\ldots,z_5) \rtimes \F[2](u,v)$.}
\end{equation*}
The first two vertical arrows are inclusions. Identifying
$\gpab[\left((B_4(\dt))^{(3)}\right)]$ with $\Z_{18}\times \Z_{18}
\times\Z^3\times (\F[2](u,v))^{(2)}/(\F[2](u,v))^{(3)}$ as above, the induced
homomorphism $\xi$ of the Abelianisations sends
$(\F[2](u,v))^{(2)}/(\F[2](u,v))^{(3)}$ and the $\Sigma$-cosets of $W_1$ and
$W_2$ onto the trivial element, and 
\begin{equation*}
\xi(W_3)= -\xi(W_4)= -\xi(W_5)=(1,0,0).
\end{equation*}

Let us determine $\ker{\epsilon}$. Since $\epsilon$ is Abelianisation, it follows from \repr{gammasemi} that 
\begin{equation*}
\ker{\epsilon}=\Gamma_2\bigl(\F[5](x_i) \rtimes \F[2](u,v)\bigr)= L\rtimes \Gamma_2\left(\F[2](u,v)\right),
\end{equation*}
where $L$ is the subgroup of $\F[5](z_i)$ generated by
$\Gamma_2\left(\F[5](z_i)\right)$ and the normal subgroup generated by the
elements of the form $\phi(g)(h)h^{-1}$, where $g\in \F[2](u,v)$ and $h\in
\F[5](z_i)$. Let $\map{\rho}{{\F[5](z_i)}}[\Z]$ be the restriction of $\epsilon$
to the first factor, in other words, 
\begin{equation}\label{eq:kerrho} 
\rho(z_i)=
\begin{cases}
1 & \text{if $i=1,2,3$}\\
-1 & \text{if $i=4,5$.}
\end{cases}
\end{equation}

\begin{prop}
$L=\ker{\rho}$.
\end{prop}

\begin{proof}
We first apply the Reidemeister-Schreier rewriting process in order to obtain a basis of $\ker{\rho}$. Taking $X=\brak{z_1,\ldots,z_5}$ as a basis of $\F[5](z_i)$, and $U=\brak{z_1^i}_{i\in\Z}$ as a Schreier transversal, one may check that a basis of $\ker{\rho}$ is given by 
\begin{equation*}
\brak{z_1^{-i}z_jz_1^{i-1}}_{i\in\Z,\, j\in\brak{2,3}} \bigcup \brak{z_1^{-i}z_jz_1^{i+1}}_{i\in\Z,\, j\in\brak{4,5}}.
\end{equation*}
These calculations are presented in Table~\ref{tab:basis}.
\begin{table}
\begin{tabular}{|*{4}{c}||c||*{4}{c}|}
\hline
$\cdots$  & $z_1^{-3}$ & \mbox{\rule[-6pt]{0cm}{7mm}$z_1^{-2}$} & $z_1^{-1}$ & & $1$ & $z_1$ & $z_1^2$ &  $\cdots$\\
\hline
& 1 & 1 & 1 & $z_1$ & 1 & 1 & \mbox{\rule[-6pt]{0cm}{7mm}1} &\\
&  $z_1^{-3}z_2z_1^2$ & $z_1^{-2}z_2z_1$ & $z_1^{-1}z_2$ & $z_2$ &
$z_2z_1^{-1}$ & $z_1z_2z_1^{-2}$ & \mbox{\rule[-6pt]{0cm}{7mm}$z_1^2z_2z_1^{-3}$}  &\\
&  $z_1^{-3}z_3z_1^2$ & $z_1^{-2}z_3z_1$ & $z_1^{-1}z_3$ & $z_3$ &
$z_3z_1^{-1}$ & $z_1z_3z_1^{-2}$ & \mbox{\rule[-6pt]{0cm}{7mm}$z_1^2z_3z_1^{-3}$}  &\\ 
&  $z_1^{-3}z_4z_1^4$ & $z_1^{-2}z_4z_1^3$ & $z_1^{-1}z_4z_1^2$ & $z_4$ &
$z_4z_1$ & $z_1z_4$ & \mbox{\rule[-6pt]{0cm}{7mm}$z_1^2z_4z_1^{-1}$}  &\\
&  $z_1^{-3}z_5z_1^4$ & $z_1^{-2}z_5z_1^3$ & $z_1^{-1}z_5z_1^2$ & $z_5$ &
$z_5z_1$ & $z_1z_5$ & \mbox{\rule[-6pt]{0cm}{7mm}$z_1^2z_5z_1^{-1}$}  &\\
\hline
\end{tabular}
\medskip
\caption{\label{tab:basis}Determination of a basis of $\ker{\rho}$}
\end{table}

We may now show that $\ker{\rho}\subseteq L$. Indeed, since $L$ is normal in $\F[5](z_i)$, and all basis elements of $\ker{\rho}$ are conjugates of $z_2z_1^{-1}$, $z_3z_1^{-1}$, $z_4z_1$ and $z_5z_1$ by powers of $z_1$, it suffices to show that these four elements belong to $L$. This can be done by studying \req{acuxi}. First, $z_2z_1^{-1}=\phi(u)(z_1)z_1^{-1}$, so $z_2z_1^{-1}\in L$. Next, 
\begin{equation*}
\phi(u)(z_2)z_2^{-1}= z_2z_1^{-1}\cdot z_5z_3z_2^{-1}z_4^{-1}\in L,
\end{equation*}
so $z_5z_3z_2^{-1}z_4^{-1}\in L$, and
$\phi(u)(z_5)z_5^{-1}= z_2^2z_3^{-1}z_5\in L$. Thus
\begin{equation*}
\phi(u)(z_3)z_3^{-1}= z_2^2z_3^{-1}z_5\cdot  z_5z_3z_2^{-1}z_4^{-1}\cdot z_2z_3^{-1}\in L,
\end{equation*}
so $z_2z_3^{-1}=z_2z_1^{-1}(z_3z_1^{-1})^{-1}\in L$, and hence $z_3z_1^{-1}\in L$. Since 
\begin{equation*}
\phi(u)(z_5)z_5^{-1}= z_2\cdot z_2z_3^{-1}\cdot z_5z_1 \cdot z_1^{-1}(z_2z_1^{-1})z_1 \cdot z_2^{-1}\in L,
\end{equation*}
and $L$ is normal in $\F[5](z_i)$, it follows that $z_5z_1\in L$. Finally, since $z_5z_3z_2^{-1}z_4^{-1}\in L$, and
\begin{equation*}
z_5z_3z_2^{-1}z_4^{-1}= z_5z_1\cdot z_1^{-1}(z_3z_1^{-1})z_1 \cdot z_1^{-1}(z_2z_1^{-1})^{-1}z_1 \cdot (z_4z_1)^{-1},
\end{equation*}
we have $z_4z_1\in L$. This proves that $\ker{\rho}\subseteq L$.

We now prove that $L\subseteq \ker{\rho}$. Clearly
$\Gamma_2(\F[5](z_i))\subset \ker{\rho}$, and since $\ker{\rho}$ is
normal in $\F[5](z_i)$, it suffices to prove that all elements of the
form $\phi(g)(h)h^{-1}$, where $g\in \F[2](u,v)$ and $h\in
\F[5](z_i)$, belong to $\ker{\rho}$, which we do by double induction.
Let $\ell$ denote the length function defined on elements of free
groups. If $\ell(g)=0$ or $\ell(h)=0$ then the result is clearly true.
If $\ell(g)=\ell(h)=1$ then one may check directly using
equations~\reqref{acuxi} and~\reqref{acvxi} that $\phi(g)(h)h^{-1} \in
\ker{\rho}$ for $g\in \brak{u,v}$. Suppose that $g\in \brak{u^{-1},
v^{-1}}$. In order to show that $\phi(g)(h) h^{-1}$ belongs to
$\ker{\rho}$, it suffices to show that its inverse
$h(\phi(g)(h))^{-1}$ belongs to it. Since $\phi(g^{-1})$ is an
automorphism, there exists $h_1\in \F[5](z_i)$ such that
$\phi(g^{-1})(h_1)=h$. Thus
\begin{equation*}
h(\phi(g)(h))^{-1}= \phi(g^{-1})(h_1)
(\phi(g)\circ\phi(g^{-1})(h))^{-1} =\phi(g^{-1})(h_1)h_1^{-1},
\end{equation*}
and the result follows from the case $g\in\brak{u,v}$.

First suppose that $\ell(g)=1$, and that the result is true for all $h$ of length less than or equal to $n\geq 1$. Let $h'\in \F[5](z_i)$ be such that $\ell(h')=n+1$. Set $h'=hz$, where $h,z\in \F[5](z_i)$, $\ell(h)=n$ and $\ell(z)=1$. Then
\begin{equation*}
\phi(g)(h'){h'}^{-1}= \phi(g)(h)h^{-1} \cdot h(\phi(g)(z)z^{-1}) h^{-1}.
\end{equation*}
By induction, both terms on the right-hand side belong to $\ker{\rho}$, and using the fact that $\ker{\rho}$ is normal in $\F[5](z_i)$, we see that the result holds for all $g$ of length one, and all $h$. 

Now suppose that the result is true for all $g$ of length less than or equal to $n\geq 1$, and all $h$. Let $g'\in \F[2](u,v)$ be such that $\ell(g')=n+1$. Set $g'=gy$, where $g,y\in \F[2](u,v)$, $\ell(g)=n$ and $\ell(y)=1$. Then 
\begin{equation*}
\phi(g')(h)h^{-1}= \phi(g)(\phi(y)(h)) (\phi(y)(h))^{-1}\cdot \phi(y)(h)h^{-1}.
\end{equation*}
Since $\phi(y)(h)\in \F[5](z_i)$, the result follows by induction. This completes the proof of the inclusion $L\subseteq \ker{\rho}$, and thus that of the proposition.
\end{proof}

\section{The lower central and derived series of $\bmmn{1}{n}$}\label{sec:bmn1n}

Let $m=1$ and $n\geq 1$. The group $\bmmn{1}{n}$ is the fundamental group of
$\St\setminus\brak{x_1,\ldots,x_n}$, and so is a free group on $n-1$ generators. So its lower central and derived series are those of free groups of finite rank. Further details about the lower central series of such groups may be found in~\cite{H,MKS}.

\section[series of $\bmmn{m}{n}$ for $m \geq 3$ and $n\geq 2$]{The lower central and derived series of $\bmmn{m}{n}$ for $m \geq 3$ and $n\geq 2$}\label{sec:m3n2}

In this section, we prove \reth{lcdsbmsn}, which tells us that if $n\geq 2$ then for most values of $m$, the lower central and derived series of $B_m(\St\setminus\brak{x_1,\ldots, x_n})$ are constant from the commutator subgroup onwards.
\begin{varthm}[\reth{lcdsbmsn}]
Let $n\geq 2$. Then:
\begin{enumerate}[(a)]
\item\label{it:itm3} If $m\geq 3$ then
\begin{equation*}
\Gamma_3(B_m(\St\setminus\brak{x_1,\ldots, x_n}))= \Gamma_2(B_m(\St\setminus\brak{x_1,\ldots, x_n})).
\end{equation*}
\item\label{it:itm5} If $m\geq 5$ then
\begin{equation*}
(B_m(\St\setminus\brak{x_1,\ldots, x_n}))^{(2)}= (B_m(\St\setminus\brak{x_1,\ldots, x_n}))^{(1)}.
\end{equation*}
\item\label{it:itm4} If $m=4$ then 
\begin{equation*}
B_4(\St\setminus\brak{x_1,\ldots, x_n})\big/(B_4(\St\setminus\brak{x_1,\ldots, x_n}))^{(2)}\cong \left( \Z^2\rtimes \Z\right) \times \Z^{n-1}
\end{equation*}
where the semi-direct product structure is that of part~(\ref{it:b3deriv2}) of \repr{gorinlin}, and
\begin{equation*}
(B_4(\St\setminus\brak{x_1,\ldots, x_n}))^{(1)}\big/(B_4(\St\setminus\brak{x_1,\ldots, x_n}))^{(2)}\cong \Z^2.
\end{equation*}
Alternatively, 
\begin{equation*}
B_4(\St\setminus\brak{x_1,\ldots, x_n})\big/(B_4(\St\setminus\brak{x_1,\ldots, x_n}))^{(2)}\cong \Z^2 \rtimes \Z^n,
\end{equation*}
where $\Z^2\cong (B_4(\St\setminus\brak{x_1,\ldots, x_n}))^{(1)}\big/(B_4(\St\setminus\brak{x_1,\ldots, x_n}))^{(2)}$ is the free Abelian group with basis $\brak{\overline{u}, \overline{v}}$, $\Z^n\cong \gpab[B_4(\St\setminus\brak{x_1,\ldots, x_n})]$ has basis $\sigma,\rho_1, \ldots,\rho_{n-1}$, and the action is given by 
\begin{align*}
\sigma\cdot \overline{u} &= \overline{v} & \sigma\cdot \overline{v} &= -\overline{u}+ \overline{v}\\
\rho_i\cdot \overline{u} &= \overline{u} & \rho_i\cdot \overline{v} &= \overline{v}
\end{align*}
for all $1\leq i\leq n-1$.
\end{enumerate}
\end{varthm}

\begin{rems}\mbox{}
\begin{enumerate}[(a)]
\item For the lower central series of $\bmmn{m}{n}$, the only case not covered by \reth{lcdsbmsn} is $m=2$ and $n\geq 2$; it will be discussed in Sections~\ref{sec:mgeq2} and~\ref{sec:lcgsb22}.
\item For the derived series of $\bmmn{m}{n}$, the outstanding cases are $n\geq 2$ and $m=2$, $m=3$ and $m=4$ (see Sections~\ref{sec:mgeq2},~\ref{sec:lcgsb22},~\ref{sec:affineatil} and~\ref{sec:bmn3}).
\end{enumerate}
\end{rems}

\begin{proof}[Proof of \reth{lcdsbmsn}]
The idea of much of the proof is similar to that of \reth{dsbn}. 
Let $m,n\geq 2$. Set $B_{m,n}=B_m(\St\setminus\brak{x_1,\ldots, x_n})$. Then we have a short exact sequence
\begin{equation*}
1\to \Gamma_2(B_{m,n})\to B_{m,n}\stackrel{\alpha}{\to} \gpab[(B_{m,n})] \to 1.  \end{equation*}
From \repr{abbm}, $\gpab[(B_{m,n})]$ is a free Abelian group of rank~$n$, generated by $\rho_1,\ldots,\rho_n,\sigma$, and subject to a single relation $\rho_1\cdots \rho_n\sigma^{2(m-1)}=1$. Taking the generators of $B_{m,n}$ given by \repr{presbetanm}, all of the $\sigma_i$ are identified to $\sigma$ by $\alpha$, and for each $1\leq i\leq n$, all of the $A_{i,j}$, $n+1\leq j\leq n+m$, are identified to $\rho_i$.

Let $H\subseteq \Gamma_2(B_{m,n})$ be a normal subgroup of $B_{m,n}$, and let 
\begin{equation*}
\left\{ \begin{aligned}
\pi \colon\thinspace B_{m,n} & \to B_{m,n}/H\\
\beta &\mapsto \overline{\beta}
\end{aligned}\right.
\end{equation*}
denote the canonical projection. Then $\alpha$ factors through $B_{m,n}/H$, and we have a short exact sequence of the form:
\begin{equation*}
\xymatrix{%
1\ar[r] & K \ar[r] & B_{m,n}/H \ar[r]^{\widehat{\alpha}} & \gpab[(B_{m,n})] \ar[r] & 1,} 
\end{equation*}
where $\alpha=\widehat{\alpha}\circ \pi$ and $K=\Gamma_2(B_{m,n})/H$.

In what follows, we shall impose one of the following two hypotheses:
\begin{enumerate}[(i)]
\item\label{it:kab} $K$ is Abelian.
\item\label{it:kcent} $K$ is central in $B_{m,n}/H$.
\end{enumerate}

\begin{rems}\mbox{}
\begin{enumerate}[(a)]
\item If $H=\Gamma_3(B_{m,n})$ then condition~(\ref{it:kcent}) is satisfied.
\item If $H=(B_{m,n})^{(2)}$ then condition~(\ref{it:kab}) is satisfied.
\item Clearly condition~(\ref{it:kcent}) implies condition~(\ref{it:kab}).
\end{enumerate}
\end{rems}

From \repr{presbetanm}, we conclude that $\overline{\sigma_1},\ldots, \overline{\sigma_{m-1}}$ and the $\overline{A_{i,j}}$, $1\leq i\leq n$, $n+1\leq j\leq n+m$, generate $B_{m,n}/H$. Since $\alpha$ identifies the $\sigma_k$ to $\sigma$, we see that $\widehat{\alpha}$ identifies the $\overline{\sigma_k}$ to $\sigma$. So for $2\leq k\leq m-1$, there exists $t_k\in K$ such that $\overline{\sigma_k}= t_k \overline{\sigma_1}$.

Suppose first that condition~(\ref{it:kcent}) is satisfied.
Then $t_k$ commutes with $\overline{\sigma_1}$. For $1\leq l\leq m-2$,
we deduce from applying $\pi$ to Artin's relations
$\overline{\sigma_l}\, \overline{\sigma_{l+1}} \,\overline{\sigma_l}=
\overline{\sigma_{l+1}}\,\overline{\sigma_l} \,
\overline{\sigma_{l+1}}$ that $t_2=1$ and $t_l=t_{l+1}$ if $l\geq 2$.
Thus $\overline{\sigma_1}=\cdots =\overline{\sigma_{m-1}}$. This
argument holds for all $m\geq 2$.

Now suppose instead that condition~(\ref{it:kab}) is
satisfied. Let $m\geq 5$. If $3\leq k\leq m-1$, then since
$\sigma_k$ commutes with $\sigma_1$, we have that
\begin{equation}\label{eq:oversig} 
\overline{\sigma_1}\cdot t_k\overline{\sigma_1}
=t_k\overline{\sigma_1}\cdot \overline{\sigma_1},
\end{equation}
and hence $t_k$ commutes with $\overline{\sigma_1}$.
Now let $4\leq l\leq m-1$ (such an $l$ exists). Since $\sigma_l$ commutes with $\sigma_2$, we obtain
\begin{equation*}
t_l\overline{\sigma_1}\cdot t_2\overline{\sigma_1} =t_2\overline{\sigma_1}\cdot t_l\overline{\sigma_1}.
\end{equation*}
But $K$ is Abelian, and thus it follows that $t_2$ commutes with
$\overline{\sigma_1}$. Applying this to the image under $\pi$ of the relation
$\sigma_1\sigma_2\sigma_1= \sigma_2\sigma_1\sigma_2$, we see that $t_2=t_2^2$,
and hence $t_2=1$. Finally, if $l\geq 2$ then the relation
$\sigma_l\sigma_{l+1}\sigma_l= \sigma_{l+1}\sigma_l\sigma_{l+1}$ implies that
$t_l=t_{l+1}$, and so $t_2=\cdots =t_{m-1}$. Hence
$\overline{\sigma_1}=\overline{\sigma_2}=\cdots =\overline{\sigma_{m-1}}$.

Let us now consider the $A_{i,j}$. In what follows, we suppose that
$m\geq 3$ and $\overline{\sigma_1}=\overline{\sigma_2}=\cdots
=\overline{\sigma_{m-1}}$ (which as we have just observed, is the case if either condition~(\ref{it:kcent}) holds, or if condition~(\ref{it:kab}) holds and additionally $m\geq 5$). By \repr{presbetanm}, $\sigma_{j-n}
A_{i,j}\sigma_{j-n}^{-1}= A_{i,j+1}$ where $n+1\leq j\leq n+m-1$, and
if $r\neq j-n-1,j-n$ then $\sigma_r A_{i,j}\sigma_r^{-1}=A_{i,j}$.
Since for all $n+1\leq j\leq n+m$, 
\begin{equation*}
A_{i,j}=\sigma_{j-n-1}\cdots
\sigma_1\cdot A_{i,n+1} \sigma_1^{-1}\cdots
\sigma_{j-n-1}^{-1},
\end{equation*}
we see by projecting into $B_{m,n}/H$, and using the condition $m\geq 3$ that
\begin{equation*}
\overline{A_{i,j}}= \overline{\sigma_{m-1}}^{j-n-1}
\overline{A_{i,n+1}}\, \overline{\sigma_{m-1}}^{-(j-n-1)}=
\overline{A_{i,n+1}}=\overline{\alpha_i},
\end{equation*}
where $\alpha_i=A_{i,n+1}$. Taking $j=n+2$, it follows that
$\overline{\alpha_i}$ commutes with $\overline{\sigma_1}$. Applying this to the
first relation of the presentation of $B_{m,n}$ given in~\repr{presbetanm}, it
follows that the $\overline{\alpha_i}$ commute pairwise. Projecting the
remaining relations for $B_{m,n}$ into $B_{m,n}/H$ give nothing new, except for
the surface relation which yields $\overline{\alpha_1}\cdots
\overline{\alpha_n}\, \overline{\sigma_1}^{2(m-1)}=1$. Hence $B_{m,n}/H$ is an
Abelian group generated by $\overline{\alpha_1},\ldots, \overline{\alpha_n}$ and $\overline{\sigma_1}$, in which the relation $\overline{\alpha_1}\cdots
\overline{\alpha_n}\, \overline{\sigma_1}^{2(m-1)}=1$ is satisfied. Since
$\widehat{\alpha}$ is surjective, we conclude by the Hopfian property of free
Abelian groups of finite rank that $\widehat{\alpha}$ is an isomorphism, so
$B_{m,n}/H\cong \Z^n$, and $H=\Gamma_2(B_{m,n})$.

Taking $m\geq 3$ and $H=\Gamma_3(B_{m,n})$, condition~(\ref{it:kcent}) is satisfied, and we conclude from the above arguments that $\Gamma_3(B_{m,n})= \Gamma_2(B_{m,n})$. This proves part~(\ref{it:itm3}) of the proposition.

Taking $m\geq 5$ and $H=(B_{m,n})^{(2)}$, condition~(\ref{it:kab}) is satisfied, and we conclude similarly that $(B_{m,n})^{(2)}=\Gamma_2(B_{m,n})= (B_{m,n})^{(1)}$, which proves part~(\ref{it:itm5}) of the proposition.

Now let us prove part~(\ref{it:itm4}) of the proposition. Let $m=4$
and $n\geq 2$, and making use of the previous notation, set
$H=(B_{4,n})^{(2)}$. Then condition~(\ref{it:kab}) holds, and indeed
$(B_{4,n})^{(1)}/H$ is Abelian. As in the case $m\geq 5$ (cf.\
\req{oversig}), we see that $t_3$ commutes with $\overline{\sigma_1}$.
From the remaining two Artin relations, we see that
$\overline{\sigma_1} t_2 \overline{\sigma_1}= t_2
\overline{\sigma_1}^2 t_2$ and $t_3 \overline{\sigma_1} t_2
\overline{\sigma_1} t_3= t_2 \overline{\sigma_1} t_3
\overline{\sigma_1} t_2$. But $t_3$ commutes with both $t_2$ and
$\overline{\sigma_1}$, hence the second equation reduces to $t_3
\overline{\sigma_1} t_2 \overline{\sigma_1}= t_2 \overline{\sigma_1}^2
t_2$. From the first equation, we see that $t_3=1$, in other words,
$\sigma_1$ and $\sigma_3$ are identified under $\pi$ to
$\overline{\sigma}$ say, and there is just one Artin relation of the
form $\overline{\sigma} \, \overline{\sigma_2}\,
\overline{\sigma}=\overline{\sigma_2} \, \overline{\sigma} \,
\overline{\sigma_2}$. Further, for $i=1,\ldots,n$, 
\begin{align*}
\overline{A_{i,n+2}}&= \overline{\sigma_1}\, \overline{A_{i,n+1}} \, \overline{\sigma_1}^{-1} = \overline{\sigma_3}\, \overline{A_{i,n+1}} \, \overline{\sigma_3}^{-1} = \overline{A_{i,n+1}},\\
\overline{A_{i,n+3}}&= \overline{\sigma_2}\, \overline{A_{i,n+2}} \, \overline{\sigma_2}^{-1} = \overline{\sigma_2}\, \overline{A_{i,n+1}} \, \overline{\sigma_2}^{-1} = \overline{A_{i,n+1}},\;\text{and}\\
\overline{A_{i,n+4}}&= \overline{\sigma_3}\, \overline{A_{i,n+3}} \, \overline{\sigma_3}^{-1} = \overline{\sigma_3}\, \overline{A_{i,n+1}} \, \overline{\sigma_3}^{-1} = \overline{A_{i,n+1}}.
\end{align*}
So for each $i=1,\ldots,n$, the $A_{i,j}$ are identified by $\pi$ to a single element $\overline{\alpha_i}$ which commutes with both $\overline{\sigma}$ and $\overline{\sigma_2}$. So $B_{4,n}/(B_{4,n})^{(2)}$ is generated by $\overline{\sigma}, \overline{\sigma_2}, \overline{\alpha_1},\ldots, \overline{\alpha_{n-1}}$, subject to the relations
\begin{align*}
\overline{\sigma} \, \overline{\sigma_2}\, \overline{\sigma}=\overline{\sigma_2} \, \overline{\sigma} \, \overline{\sigma_2}, && \overline{\sigma}, \overline{\sigma_2} \comm \overline{\alpha_i}, && \overline{\alpha_i},\comm \overline{\alpha_j}, && \text{for all $1\leq i,j\leq n-1$.}
\end{align*}
So there exist homomorphisms 
\begin{equation*}
\map{f}{B_{4,n}}[B_3(\dt)\times \Z^{n-1}]
\end{equation*}
given by 
\begin{align*}
f(\sigma_1)&=f(\sigma_3)=\sigma_1\\
f(\sigma_2)&=\sigma_2\\
f(A_{i,j})&=A_i\;\text{if $1\leq i\leq n-1$}\\
f(A_{n,j})&=A_{n-1}^{-1}\cdots A_1^{-1}\sigma_1^{-1} \sigma_2^{-1} \sigma_1^{-2} \sigma_2^{-1} \sigma_1^{-1},
\end{align*}
where $(A_1,\ldots, A_{n-1})$ is a basis of $\Z^{n-1}$, and 
\begin{equation*}
\map{\widetilde{\pi}}{B_3(\dt)\times \Z^{n-1}}[B_{4,n}/(B_{4,n})^{(2)}]
\end{equation*}
defined by 
\begin{align*}
\widetilde{\pi}(\sigma_1)&=\overline{\sigma}\\
\widetilde{\pi}(\sigma_2)&=\overline{\sigma_2}\\
\widetilde{\pi}(A_i)&=\overline{\alpha_i},
\end{align*}
and satisfying $\pi=\widetilde{\pi}\circ f$. Since $f$ is surjective, to prove the first part of~(\ref{it:itm4}), by \repr{gorinlin}(\ref{it:b3deriv2}), it suffices to show that $\ker{\widetilde{\pi}}= (B_3(\dt))^{(2)}$.

First let us show that $(B_3(\dt))^{(2)}\subseteq \ker{\widetilde{\pi}}$. Let $y\in (B_3(\dt))^{(2)}$. In particular, $y$ may be written as a word $w(\sigma_1,\sigma_2)$. Considering this word to be an element $x$ of $B_{4,n}$, since $y\in (B_3(\dt))^{(2)}$, we have that $x\in (B_{4,n})^{(2)}$ and $f(x)=y$. Since $\pi(x)=e$, it follows that $y\in \ker{\widetilde{\pi}}$.

Conversely, let $y\in \ker{\widetilde{\pi}}$. Since $f$ is surjective, there exists $x\in B_{4,n}$ such that $f(x)=y$, and so $x\in (B_{4,n})^{(2)}$. But since $(B_{4,n})^{(1)}$ is the normal subgroup of $B_{4,n}$ generated by the commutators $[\sigma_k,\sigma_l]$, $[\sigma_k, A_{i,j}]$ and $[A_{i,j}, A_{i',j'}]$, where $1\leq k,l\leq 3$, $1\leq i,i'\leq n$ and $1\leq j,j'\leq 4$, $f$ is surjective and $\Z^{n-1}$ is a direct factor of $B_3(\dt)\times \Z^{n-1}$, it follows that $f\left( (B_{4,n})^{(1)}\right) = (B_3(\dt))^{(1)}$, and thus $f\left( (B_{4,n})^{(2)}\right) = (B_3(\dt))^{(2)}$. In particular, $y\in (B_3(\dt))^{(2)}$. We thus conclude that $B_3(\dt)/(B_3(\dt))^{(2)} \times \Z^{n-1} \cong B_{4,n}/(B_{4,n})^{(2)}$, which proves the first part of~(\ref{it:itm4}).

We now move on to the second part of~(\ref{it:itm4}). Consider the homomorphism $B_3(\dt)\to B_{4,n}$ given by $\sigma_i\mapsto \sigma_i$. Since $g\left( (B_3(\dt))^{(i)}\right) \subseteq (B_{4,n})^{(i)}$ for all $i\in\N$, there is an induced homomorphism
\begin{equation*}
\map{g}{B_3(\dt)/(B_3(\dt))^{(2)}}[B_{4,n}/(B_{4,n})^{(2)}],
\end{equation*}
which sends the coset of $\sigma_i$ onto $\overline{\sigma_i}$, as well as its restriction 
\begin{equation*}
\map{g_1}{(B_3(\dt))^{(1)}/(B_3(\dt))^{(2)}}[(B_{4,n})^{(1)}/(B_{4,n})^{(2)}].
\end{equation*}
Similarly, since $B_4(\St\setminus \brak{x_1}) \cong B_4(\dt)$ by \repr{iso}, the surjective homomorphism $B_4(\St\setminus \brak{x_1,\ldots,x_n})\to B_4(\St\setminus \brak{x_1})$ given by closing up the $n-1$~punctures $x_2,\ldots, x_{n}$ induces a surjective homomorphism 
\begin{equation*}
\map{h}{B_{4,n}/(B_{4,n})^{(2)}}[B_4(\dt)/(B_4(\dt))^{(2)}],
\end{equation*}
which sends $\overline{\sigma_i}$ onto the coset of $\sigma_i$, as well as its restriction 
\begin{equation*}
\map{h_1}{(B_{4,n})^{(1)}/(B_{4,n})^{(2)}}[B_4(\dt)^{(1)}/(B_4(\dt))^{(2)}]. 
\end{equation*}
Hence we obtain the following commutative diagram:
\begin{equation*}
\xymatrix{%
(B_3(\dt))^{(1)}/(B_3(\dt))^{(2)} \ar[r]\ar[d]_{g_1}  &   B_3(\dt)/(B_3(\dt))^{(2)} \ar[r]\ar[d]_{g}    & \gpab[(B_3(\dt))]\ar[d]\\
(B_{4,n})^{(1)}/(B_{4,n})^{(2)}  \ar[r] \ar[d]_{h_1} &  B_{4,n}/(B_{4,n})^{(2)}  \ar[r] \ar[d]_{h} & \gpab[(B_{4,n})]  \ar[d] \\
(B_4(\dt))^{(1)}/(B_4(\dt))^{(2)} \ar[r]  &   B_4(\dt)/(B_4(\dt))^{(2)} \ar[r]   & \gpab[(B_4(\dt))].}
\end{equation*}
Note that the rows are all short exact sequences.

Now consider the first column. From~\cite{GL} and \repr{gorinlin}, we know that
$(B_3(\dt))^{(1)}/(B_3(\dt))^{(2)}$ and $(B_4(\dt))^{(1)}/(B_4(\dt))^{(2)}$ are
both free Abelian groups of rank~$2$, generated by their respective cosets of
$u=\sigma_2\sigma_1^{-1}$ and  $v= \sigma_1\sigma_2\sigma_1^{-2}$. By definition
of $g$ and $h$, it follows that $h_1\circ g_1$ is an isomorphism, sending the
$(B_3(\dt))^{(2)}$-coset of $u$ (respectively $v$) onto the
$(B_4(\dt))^{(2)}$-coset of $u$ (respectively $v$). Thus to prove that
$(B_{4,n})^{(1)}/(B_{4,n})^{(2)}\cong \Z^2$, it suffices to show that $g_1$ is
surjective. To see this, let $x\in (B_{4,n})^{(1)}/(B_{4,n})^{(2)}$. Since $x\in
B_{4,n}/(B_{4,n})^{(2)}$, it follows from above that 
\begin{equation*}
x=w(\overline{\sigma}, \overline{\sigma_2})\;\overline{\alpha_1}^{m_1}\cdots \overline{\alpha_{n-1}}^{m_{n-1}},
\end{equation*}
where $m_i\in\Z$ for $1\leq i\leq n-1$ and $w(\overline{\sigma}, \overline{\sigma_2})$ is a word in $\overline{\sigma}$ and $\overline{\sigma_2}$. Projecting into $\gpab[(B_{4,n})]$, since $\overline{\sigma}$ and $\overline{\sigma_2}$ map onto $\sigma$, and $\overline{\alpha_i}$ maps onto $\rho_i$, and furthermore, $\sigma,\rho_1, \ldots \rho_{n-1}$ generate freely $\gpab[(B_{4,n})]$, we see by exactness that the $m_i$ are all zero, in other words, $x=w(\overline{\sigma}, \overline{\sigma_2})$. Now take $z=w(\sigma_1,\sigma_2)\in B_3(\dt)/(B_3(\dt))^{(2)}$, so that $g(z)=x$. Projecting $z$ into $\gpab[(B_3(\dt))]$ yields zero by commutativity of the diagram (the homomorphism $\gpab[(B_3(\dt))]\to \gpab[(B_{4,n})]$ is injective), hence $z\in (B_3(\dt))^{(1)}/(B_3(\dt))^{(2)}$, and thus $g_1$ is surjective. So $(B_4(\St\setminus\brak{x_1,\ldots, x_n}))^{(1)}/(B_4(\St\setminus\brak{x_1,\ldots, x_n}))^{(2)}\cong \Z^2$, which proves the second part of part~(\ref{it:itm4}).

Finally, we prove the last part of part~(\ref{it:itm4}). Consider the short exact sequence
\begin{equation*}
1\to (B_{4,n})^{(1)}/(B_{4,n})^{(2)}  \to  B_{4,n}/(B_{4,n})^{(2)}  \stackrel{\widehat{\alpha}}{\to} \gpab[(B_{4,n})]\to 1.
\end{equation*}
Recall that $\gpab[(B_{4,n})]$ is a free Abelian group with basis $\brak{\sigma,\rho_1,\ldots, \rho_{n-1}}$, and that up to isomorphism, we may identify $B_3(\dt)/(B_3(\dt))^{(2)} \times \Z^{n-1}$ with $B_{4,n}/(B_{4,n})^{(2)}$, where the $\Z^{n-1}$-factor has a basis $\brak{\overline{\alpha_1},\ldots, \overline{\alpha_{n-1}}}$ for which $\widehat{\alpha}(\overline{\alpha_i})=\rho_i$. It follows that $\widehat{\alpha}$ admits a section given by $\sigma\mapsto \overline{\sigma}$ and $\rho_i\mapsto \overline{\alpha_i}$, and hence 
\begin{equation*}
B_{4,n}/(B_{4,n})^{(2)} \cong (B_{4,n})^{(1)}/(B_{4,n})^{(2)}\rtimes \gpab[(B_{4,n})].
\end{equation*}
Taking the basis $\brak{\overline{u}, \overline{v}}$ of $(B_{4,n})^{(1)}/(B_{4,n})^{(2)} \cong \Z^2$, the action is given by $\rho_i\cdot \overline{u}= \overline{\alpha_i}\, \overline{u}\, \overline{\alpha_i}^{-1}=\overline{u}$ and $\rho_i\cdot \overline{v}= \overline{\alpha_i}\, \overline{v}\, \overline{\alpha_i}^{-1}=\overline{v}$ since $\overline{\alpha_i}$ commutes with $\overline{\sigma}$ and $\overline{\sigma_2}$, and $\sigma\cdot \overline{u}$ and $\sigma\cdot \overline{v}$ are obtained as in the proof of the second part of \repr{gorinlin}(\ref{it:gorin3}). This completes the proof of \reth{lcdsbmsn}. 
\end{proof}

\section{The commutator subgroup of $B_m(\St\setminus\brak{x_1,x_2})$, $m\geq 2$}\label{sec:mgeq2}

Let $m\geq 2$. As we saw in \reth{lcdsbmsn}, the lower central series of $B_m(\St\setminus\brak{x_1,\ldots, x_n}))$, $n\geq 2$, is constant from the commutator subgroup onwards if $m\geq 3$. In this section, we study the case $n=2$ in more detail. The special case $m=n=2$ will also be analysed later in \resec{lcgsb22}, and the case $m\geq 3$ and $n=2$ will also be discussed in \resec{affineatil}.

From \rerems{annulus}, we know that $B_m(\St\setminus\brak{x_1, x_2}))$ is the $m$-string braid group of the annulus, and so is isomorphic to the Artin group of type~$B_m$. Presentations of these groups were obtained in~\cite{Lam,Ma}, as well as in~\cite{KP} (we will come back to this presentation in \repr{kent}). Annulus braid groups were also studied in~\cite{Cr,PR}.

Let $m\geq 2$. From \repr{iso}, it follows from part~(\ref{it:isob}) that $B_m(\St\setminus\brak{x_1,x_2})\cong B_m(\dt\setminus \brak{x_2})$, and from part~(\ref{it:isoc}) that 
\begin{equation*}
B_m(\dt\setminus \brak{x_2}, \brak{x_3,\ldots, x_{m+1}, x_{m+2}}) \cong B_{m,1}(\dt).
\end{equation*}
Hence $B_m(\St\setminus\brak{x_1,x_2})\cong B_{m,1}(\dt)$. But from part~(\ref{it:isod}), 
\begin{align*}
B_{m,1}(\dt)&\cong \pi_1(\dt\setminus\brak{x_3,x_4,\ldots, x_{m+2}},x_2)\rtimes B_m(\dt)\\
&\cong \F[m](A_{2,3}, \ldots, A_{2,m+2}) \rtimes B_m(\dt),
\end{align*}
where $B_m(\dt)$ is taken to be generated by $\sigma_3, \ldots, \sigma_{m+1}$, and the action $\phi$ of the $\sigma_i$, $3\leq i\leq m+1$, on the $A_{2,j}$, $3\leq j\leq m+2$ is that given by the Artin representation:
\begin{equation}\label{eq:sigij}
\sigma_i A_{2,j} \sigma_i^{-1}= 
\begin{cases}
A_{2,j+1} & \text{if $j=i$}\\
A_{2,j}^{-1} A_{2,j-1}A_{2,j} & \text{if $j=i+1$}\\
A_{2,j} & \text{otherwise.}
\end{cases}
\end{equation} 
From this, we may deduce that:
\begin{equation}\label{eq:sigiji}
\sigma_i^{-1} A_{2,j} \sigma_i= 
\begin{cases}
A_{2,j-1} & \text{if $j=i+1$}\\
A_{2,j} A_{2,j+1}A_{2,j}^{-1} & \text{if $j=i$}\\
A_{2,j} & \text{otherwise.}
\end{cases}
\end{equation} 

\begin{varthm}[\repr{aug}]
Let $m\geq 2$. Then:
\begin{enumerate}[(a)]
\item\label{it:bmit1} $B_m(\St\setminus\brak{x_1,x_2})\cong \F[m] \rtimes B_m(\dt)$, where the action $\phi$ is as given in \req{sigij}.
\item\label{it:bmit2} $\Gamma_2(B_m(\St\setminus\brak{x_1,x_2})) \cong \ker{\rho} \rtimes \Gamma_2(B_m(\dt))$, where
\begin{equation*}
\map{\rho}{\F[m](A_{2,3}, \ldots, A_{2,m+2})}[\Z]
\end{equation*}
is the augmentation homomorphism, and the action is that induced by $\phi$.
\end{enumerate}
\end{varthm}

\begin{proof}
Part~(\ref{it:bmit1}) was proved above, and in any case is a restatement of the results of \repr{iso}. So let us prove part~(\ref{it:bmit2}). Set $\F[m]= \F[m](A_{2,3}, \ldots, A_{2,m+2})$, and let $L$ be the subgroup of $\F[m]$ generated by $\Gamma_2(\F[m])$ and the normal subgroup generated by the elements of the form $\phi(g)(h)\cdot h^{-1}$, where $g\in B_m(\dt)$ and $h\in\F[m]$. By \repr{gammasemi}, it suffices to prove that $L=\ker{\rho}$.

First we show that $L\subseteq \ker{\rho}$. Since $\rho$ factors through Abelianisation, we have clearly that $\Gamma_2(\F[m]) \subseteq \ker{\rho}$. Further, since $\ker{\rho}$ is normal in $\F[m]$, it suffices to prove that $\phi(g)(h)\cdot h^{-1}\in \ker{\rho}$, where $g\in B_m(\dt)$ and $h\in\F[m]$. This is equivalent to showing that $\rho(h)= \rho(\phi(g)(h))= \rho(ghg^{-1})$ and may be achieved by double induction as follows. If $g$ and $h$ are both of length~$1$, in other words if they are generators or inverses of generators of their respective groups then the result holds using equations~\reqref{sigij} and~\reqref{sigiji}. Secondly, if $g$ is of length~$1$ then the result follows for all $h$ by applying induction on the word length of $h$ (relative to the given basis of $\F[m]$) and the fact that 
\begin{equation*}
gh_1h_2g^{-1}= gh_1g^{-1}\cdot gh_2g^{-1}
\end{equation*}
for all $h_1,h_2\in \F[m]$. Finally the result holds for all $g$ and all $h$ by applying induction on the word length of $g$ (relative to the given generators of $B_m(\dt)$) and the relation 
\begin{equation*}
g_1g_2 h (g_1g_2)^{-1}=g_1h' g_1^{-1},
\end{equation*}
where $g_1,g_2\in B_m(\dt)$ and $h'=g_2h g_2^{-1}\in \F[m]$. This proves that $L\subseteq \ker{\rho}$.

To see that $\ker{\rho}\subseteq L$, we determine a basis of $\ker{\rho}$ with the help of the Reidemeister-Schreier rewriting process. Taking $X=\brak{A_{2,3},\ldots, A_{2,m+2}}$ as a basis of $\F[m]$ and $U=\brak{A_{2,3}^i}_{i\in\Z}$ to be a Schreier transversal, we see that a basis of $\ker{\rho}$ is given by the elements of the form $\brak{A_{2,3}^i A_{2,j}A_{2,3}^{-(i+1)}}_{i\in\Z,\; j\in\brak{4,\ldots, m+2}}$, or in other words, the conjugates of the $A_{2,j}A_{2,3}^{-1}$ by $A_{2,3}^i$. Since $L$ is normal in $\F[m]$, it suffices to prove that the $A_{2,j}A_{2,3}^{-1}$ belong to $L$. This is the case, since for all $3\leq j\leq m+1$, 
\begin{gather*}
\text{$\phi(\sigma_j)(A_{2,j})A_{2,j}^{-1}=A_{2,j+1}A_{2,j}^{-1}\in L$, and}\\
A_{2,j+1}A_{2,3}^{-1}= A_{2,j+1}A_{2,j}^{-1}\cdot A_{2,j}A_{2,j-1}^{-1}\cdots A_{2,4}A_{2,3}^{-1}\in L.
\end{gather*}
Thus $\ker{\rho}\subseteq L$, which completes the proof of part~(\ref{it:bmit2}) of the proposition.
\end{proof}

We now investigate further the case $m=3$. By \reth{gorinlin}, $\Gamma_2(B_3(\dt))$ is a free group $\F[2](u,v)$ of rank~$2$, where $u=\sigma_4 \sigma_3^{-1}$ and $v=\sigma_3 \sigma_4 \sigma_3^{-2}$. For $i\in\Z$, we set
\begin{align*}
\alpha_i &= A_{2,3}^i A_{2,4}A_{2,3}^{-(i+1)}= A_{2,3}^i \alpha_0 A_{2,3}^{-i} \quad \text{and}\\
 \beta_i &= A_{2,3}^i A_{2,5} A_{2,3}^{-(i+1)}= A_{2,3}^i \beta_0 A_{2,3}^{-i}.
\end{align*}
Using relations~\reqref{sigij} and~\reqref{sigiji}, one may check that 
\begin{align*}
uA_{2,3}u^{-1}&= A_{2,3} A_{2,5} A_{2,3}^{-1} & vA_{2,3}v^{-1}&= A_{2,4}A_{2,5} A_{2,4} A_{2,5}^{-1} A_{2,4}^{-1}\\
uA_{2,4}u^{-1}&= A_{2,3} & vA_{2,4}v^{-1}&= A_{2,4}A_{2,5} A_{2,4}^{-1}\\
uA_{2,5}u^{-1}&= A_{2,5}^{-1} A_{2,4}A_{2,5}&
vA_{2,5}v^{-1}&= A_{2,5}^{-1} A_{2,4}^{-1} A_{2,3} A_{2,4}A_{2,5},
\end{align*}
then that
\begin{equation}\label{eq:uau}
\left.\begin{aligned}
u\alpha_0u^{-1} =& u A_{2,4}A_{2,3}^{-1} u^{-1} = (A_{2,3}A_{2,5}A_{2,3}^{-2})^{-1}= \beta_1^{-1}\\
u\beta_0 u^{-1}=& u A_{2,5}A_{2,3}^{-1} u^{-1} = (A_{2,3}^{-1} A_{2,5})^{-1} (A_{2,3}^{-1} A_{2,4}) (A_{2,5} A_{2,3}^{-1})\cdot\\
&(A_{2,3}A_{2,5} A_{2,3}^{-2})^{-1}
= \beta_{-1}^{-1} \alpha_{-1} \beta_0 \beta_1^{-1}\\
u A_{2,3}^i u^{-1} =& r_i \cdot A_{2,3}^i,
\end{aligned}\right\}
\end{equation}
where $r_i= A_{2,3} A_{2,5}^i A_{2,3}^{-(i+1)} \in\ker{\rho}$, and finally that
\begin{equation}\label{eq:vav}
\left.\begin{aligned}
v\alpha_0 v^{-1}=& v A_{2,4}A_{2,3}^{-1} v^{-1} = (A_{2,4} A_{2,3}^{-1}) (A_{2,3}A_{2,5}A_{2,3}^{-2}) \cdot\\
& (A_{2,3}^2A_{2,5}A_{2,3}^{-3})(A_{2,3}^2A_{2,4}A_{2,3}^{-3})^{-1} (A_{2,3} A_{2,5}A_{2,3}^{-2})^{-1}\cdot\\ 
& (A_{2,4} A_{2,3}^{-1})^{-1}= \alpha_0 \beta_1 \beta_2 \alpha_2^{-1} \beta_1^{-1} \alpha_0^{-1}\\
v\beta_0 v^{-1}=& v A_{2,5}A_{2,3}^{-1} v^{-1} = (A_{2,3}^{-1} A_{2,5})^{-1} (A_{2,3}^{-2}A_{2,4} A_{2,3})^{-1} \cdot\\ & (A_{2,3}^{-1} A_{2,4})  (A_{2,5} A_{2,3}^{-1}) (A_{2,3}A_{2,4}A_{2,3}^{-2})
(A_{2,3}^2 A_{2,5}A_{2,3}^{-3}) \cdot\\ & (A_{2,3}^2 A_{2,4}A_{2,3}^{-3})^{-1} (A_{2,3} A_{2,5}A_{2,3}^{-2})^{-1} (A_{2,4} A_{2,3}^{-1})^{-1}\\
&= \beta_{-1}^{-1} \alpha_{-2}^{-1} \alpha_{-1} \beta_0 \alpha_1 \beta_2 \alpha_2^{-1} \beta_1^{-1} \alpha_0^{-1}\\
v A_{2,3}^i v^{-1} =& s_i A_{2,3}^{-i},
\end{aligned}\right\}
\end{equation}
where $s_i= A_{2,4}A_{2,5} A_{2,4}^i A_{2,5}^{-1} A_{2,4}^{-1}A_{2,3}^{-i} \in\ker{\rho}$. Up to expressing the $r_i$ and $s_i$ in terms of the $\alpha_i$ and $\beta_i$, we thus obtain a complete set of relations for $\ker{\rho}\rtimes \F[2](u,v)$:
\begin{equation}\label{eq:ualpha}
\left.\begin{aligned}
u\alpha_i u^{-1} &= r_i \beta_{i+1}^{-1} r_i^{-1}\\
u\beta_i u^{-1} &= r_i \beta_{i-1}^{-1} \alpha_{i-1} \beta_i \beta_{i+1}^{-1} r_i^{-1}\\
v\alpha_i v^{-1} &= s_i \alpha_i \beta_{i+1} \beta_{i+2} \alpha_{i+2}^{-1} \beta_{i+1}^{-1} \alpha_i^{-1} s_i^{-1}\\
v\beta_i v^{-1} &= s_i \beta_{i-1}^{-1} \alpha_{i-2}^{-1} \alpha_{i-1} \beta_i \alpha_{i+1} \beta_{i+2} \alpha_{i+2}^{-1} \beta_{i+1}^{-1} \alpha_i^{-1}  s_i^{-1}
\end{aligned}\right\}
\end{equation}
for all $i\in\Z$.

Setting $\widetilde{\alpha_i}$ (resp.\ $\widetilde{\beta_i}$) to be the Abelianisation of $\alpha_i$ (resp.\ $\beta_i$), and Abelianising equations~\reqref{ualpha}, we obtain:
\begin{equation}\label{eq:alpha}
\left.
\begin{gathered}
\widetilde{\alpha_i}=-\widetilde{\beta_{i+1}}\\
\widetilde{\alpha_i}=\widetilde{\alpha_{i+3}}\\
\widetilde{\alpha_i}+\widetilde{\alpha_{i+1}}+ \widetilde{\alpha_{i+2}}=0
\end{gathered}\right\}
\end{equation}
for all $i\in\Z$. By \repr{aug}(\ref{it:bmit2}),  $\left.(B_3(\St\setminus\brak{x_1,x_2}))^{(1)}\right/ (B_3(\St\setminus\brak{x_1,x_2}))^{(2)}$ is Abelian, generated by the $\widetilde{\alpha_i}$ and $\widetilde{\beta_i}$, $i\in \Z$, and the Abelianisations of $u$ and $v$, subject to these relations, and so is a free Abelian group with basis $\widetilde{\alpha_0}, \widetilde{\alpha_1}$ and the Abelianisations of $u$ and $v$. Hence:
\begin{varthm}[\repr{b3z4}]
\begin{equation*}
\left(B_3\left(\St\setminus\brak{x_1,x_2}\right)\right)^{(1)} \left/ \left(B_3\left(\St\setminus\brak{x_1,x_2}\right)\right)^{(2)}\right. \cong \Z^4.\tag*{\qed}
\end{equation*}
\end{varthm}

With the help of \repr{b3z4}, we may obtain the following:
\begin{varthm}[\repr{z4z2}]
\begin{equation*}
B_3\left(\St\setminus\brak{x_1,x_2}\right)\left/ \left(B_3\left(\St\setminus\brak{x_1,x_2}\right)\right)^{(2)}\right. \cong \Z^4 \rtimes \Z^2,
\end{equation*}
where $\Z^4$ has a basis $\brak{\widetilde{\alpha_0}, \widetilde{\beta_0}, \widetilde{u}, \widetilde{v}}$, $\Z^2$ has a basis $\brak{\sigma,\rho_1}$, and the action is given by:
\begin{align*}
\sigma&\cdot \widetilde{u}=\widetilde{v} & \sigma&\cdot \widetilde{v}=-\widetilde{u}+ \widetilde{v}\\
\sigma&\cdot \widetilde{\alpha_0}=\widetilde{\beta_0} &\sigma&\cdot \widetilde{\beta_0}= \widetilde{\beta_0} -\widetilde{\alpha_0}\\
\rho_1&\cdot \widetilde{\alpha_0}=\widetilde{\alpha_0} &  \rho_1&\cdot \widetilde{\beta_0}=\widetilde{\beta_0}\\
\rho_1&\cdot \widetilde{u}=  -\widetilde{\alpha_0} -\widetilde{u} +\widetilde{v}&
\rho_1&\cdot \widetilde{v}=-\widetilde{\beta_0}-\widetilde{u}.
\end{align*}
\end{varthm}

\begin{proof}
Consider the following short exact sequence:
\begin{equation*}
1\to (B_{3,2})^{(1)}/(B_{3,2})^{(2)} \to B_{3,2}/(B_{3,2})^{(2)} \stackrel{\widehat{\alpha}}{\to} \gpab[B_3] \to 1,
\end{equation*}
where $B_{3,2}=B_3(\St\setminus\brak{x_1,x_2})$. From \repr{b3z4}, $(B_{3,2})^{(1)}/(B_{3,2})^{(2)}$ is a free Abelian group of rank~$4$ with basis $\widetilde{\alpha_0}$, $\widetilde{\beta_0}$, $\widetilde{u}$ and $\widetilde{v}$, where $\alpha_0=A_{2,4}A_{2,3}^{-1}$, $\beta_0=A_{2,5}A_{2,3}^{-1}$ and $\widetilde{u}$, $\widetilde{v}$ are the respective Abelianisations of $u=\sigma_4\sigma_3^{-1}$ and $v =\sigma_3 \sigma_4\sigma_3^{-2}$. Further, $\gpab[(B_{3,2})]$ is a free Abelian group of rank~$2$, with basis $\sigma,\rho_1$, where $\widehat{\alpha}(\overline{\sigma_3})= \widehat{\alpha}(\overline{\sigma_4})= \sigma$, and $\widehat{\alpha}(\overline{A_{1,j}})=\rho_1$ for $j=1,2,3$. Then $\sigma\mapsto \overline{\sigma_3}$ and $\rho_1 \mapsto \overline{A_{1,5}}$ defines a section for $\widehat{\alpha}$. Let us now determine the associated action. We have already seen that $\sigma\cdot u=v$, and $\sigma\cdot v=u^{-1}v$, so $\sigma\cdot \widetilde{u}=\widetilde{v}$, and $\sigma\cdot \widetilde{v}=-\widetilde{u}+ \widetilde{v}$. Further, from \req{sigij}, we have $\sigma\cdot \alpha_0= \sigma_3 \alpha_0\sigma_3^{-1}= A_{2,4}^{-1}A_{2,3}=\alpha_{-1}^{-1}$, and $\sigma\cdot \beta_0= \sigma_3 \alpha_0\sigma_3^{-1}= A_{2,5}A_{2,4}^{-1} =\beta_0 \alpha_0^{-1}$, so by \req{alpha}, $\sigma\cdot \widetilde{\alpha_0}=\widetilde{\beta_0}$ and $\sigma\cdot \widetilde{\beta_0}= \widetilde{\beta_0} -\widetilde{\alpha_0}$. As for the action of $\rho_1$, we have $\rho_1\cdot \alpha_0=\alpha_0$, so $\rho_1\cdot \widetilde{\alpha_0}=\widetilde{\alpha_0}$, and $\rho_1\cdot \beta_0=A_{1,5} \cdot A_{2,5} A_{2,3}^{-1} \cdot A_{1,5}^{-1}$. But $A_{1,5}= \sigma_4^{-1} \sigma_3^{-2} \sigma_4^{-1} A_{2,5}^{-1}$,
hence
\begin{align*}
\rho_1\cdot \beta_0 &=\sigma_4^{-1} \sigma_3^{-2} \sigma_3^{-1}   A_{2,5}  \sigma_4 \sigma_3^2 \sigma_4 A_{2,3}^{-1}\\
&= A_{2,3}A_{2,4} A_{2,5}A_{2,4}^{-1} A_{2,3}^{-2}= \alpha_1 \beta_2 \alpha_2^{-1},
\end{align*}
and thus $\rho_1\cdot \widetilde{\beta_0}= \widetilde{\beta_0}$. Also,
\begin{align*}
\rho_1\cdot u& =\sigma_4^{-1} \sigma_3^{-2} \sigma_3^{-1} A_{2,5}^{-1} u  A_{2,5} \sigma_4 \sigma_3^2 \sigma_4\\ 
&= A_{2,3}A_{2,4} A_{2,5}^{-1}A_{2,4}^{-1} A_{2,3}^{-1} A_{2,4} \sigma_4^{-1} \sigma_3^{-3} \sigma_4 \sigma_3^2 \sigma_4.
\end{align*}
But 
\begin{equation*}
\sigma_4^{-1} \sigma_3^{-3} \sigma_4 \sigma_3^2 \sigma_4 = \sigma_3^{-1} \sigma_4 \cdot \sigma_3 \sigma_4 \sigma_3^{-2}= u^{-1}v,
\end{equation*}
so $\rho_1\cdot u= \alpha_1 \beta_1^{-1} \alpha_0^{-1} \alpha_{-1}u^{-1}v$, and $\rho_1\cdot \widetilde{u}=  -\widetilde{\alpha_0} -\widetilde{u} +\widetilde{v}$. Finally,
\begin{align*}
\rho_1\cdot v& = A_{1,5} \sigma_3 A_{1,5}^{-1}\cdot \rho_1\cdot u\cdot A_{1,5} \sigma_3^{-1} A_{1,5}^{-1}= \sigma_3 \cdot \alpha_1 \beta_1^{-1} \alpha_0^{-1} \alpha_{-1}u^{-1}v \cdot\sigma_3^{-1}\\
&= A_{2,3}A_{2,4} A_{2,5}^{-1}A_{2,4}^{-1} A_{2,3}^{-1} A_{2,4}^{-1} A_{2,3}A_{2,4}v^{-1}u^{-1}v\\
&= \alpha_1 \beta_1^{-1} \alpha_0^{-1} \alpha_{-2}^{-1} \alpha_{-1} v^{-1}u^{-1}v,
\end{align*}
and so $\rho_1\cdot \widetilde{v}=-\widetilde{\beta_0}-\widetilde{u}$, which proves the proposition.
\end{proof}

We now give an alternative proof of \repr{b3z4}. Although it is long, we believe the method to be of interest. We will also make use of some results proved in \resec{gorinlin} to prove (\repr{f5f2}) that $(B_3(\St\setminus \brak{x_1,x_2}))^{(1)}$ is a semi-direct product of an infinite-rank subgroup of $\F[5](z_1,\ldots, z_5)$ by $\F[2](u,v)$.

Given $B_3(\St\setminus \brak{x_1,x_2})$, from the generalised Fadell-Neuwirth short exact sequence (\req{seqbnm}), we obtain
\begin{equation*}
1 \to B_3(\St\setminus \brak{x_1,x_2}) \stackrel{\iota}{\to} B_{3,1}(\St\setminus \brak{x_1}) \to B_1(\St\setminus \brak{x_1})\to 1.
\end{equation*}
Clearly $\iota$ is an isomorphism, and composing by the inclusion $B_{3,1}(\St\setminus \brak{x_1})\lhra B_4(\St\setminus \brak{x_1})$, we obtain an injective homomorphism 
\begin{equation*}
f \colon\thinspace B_3(\St\setminus \brak{x_1,x_2}) \to B_4(\St\setminus \brak{x_1}).
\end{equation*}
Further, we have the following commutative diagram of short exact sequences:
\begin{equation}\label{eq:p3st} 
\begin{xy}*!C\xybox{%
\xymatrix@C=0.7cm{%
1 \ar[r] & P_3(\St\setminus \brak{x_1,x_2}) \ar[r] \ar[d]^{\cong}_{\phi\bigr\rvert_{P_3(\St\setminus \brak{x_1,x_2})}}  & B_3(\St\setminus \brak{x_1,x_2}) \ar[r]^-{\pi} \ar[d]_{\phi} & \sn[3] \ar[r] \ar[d] & 1\\
1 \ar[r] & P_4(\dt) \ar[r] & B_4(\dt) \ar[r]^-{\pi} & \sn[4] \ar[r] & 1,}}
\end{xy}
\end{equation}
where $\pi$ is the homomorphism which to a braid associates its permutation, and $\phi$ is the composition of $f$ and the isomorphism $B_4(\St\setminus \brak{x_1})\cong B_4(\dt)$ given by \repr{iso}(\ref{it:isob}). The right-hand vertical homomorphism is the natural inclusion of $\sn[3]$ in $\sn[4]$. So $\phi$ is also injective, and $\phi(\sigma_i)= \sigma_i$ for $i=1,2$. The fact that the left-hand homomorphism $\phi\bigr\rvert_{P_3(\St\setminus \brak{x_1,x_2})}$ is an isomorphism follows from \repr{iso}(\ref{it:isoa}). From this, we obtain the following commutative diagram: 
\begin{equation*}
\xymatrix{%
(B_3(\St\setminus \brak{x_1,x_2}))^{(1)} \ar[r]^-{\pi} \ar[d]_{\phi} & \an[3] \ar[d]\\
(B_4(\dt))^{(1)} \ar[r]^-{\pi} & \an[4] ,}
\end{equation*}
$\an[n]$ being the alternating subgroup of $\sn[n]$. By abuse of
notation, we use the same symbols for the restriction homomorphisms.
If $H= (B_4(\dt))^{(1)}\cap \pi^{-1}(\an[3])$ then
$\phi(B_3(\St\setminus \brak{x_1,x_2}))^{(1)} \subseteq H$ by
commutativity of the diagram. Since
$\pi\bigr\rvert_{(B_4(\dt))^{(1)}}$ is surjective onto $A_4$, it
follows that $[(B_4(\dt))^{(1)}: H]=4$ (indeed, if $\map g{G_1}[G_2]$
is a surjective group homomorphism, and $H_2$ is a subgroup of $G_2$
then $H_1=g^{-1}(H_2)$ is a subgroup of $G_1$, and $G_1/H_1\to
G_2/H_2$, $xH_1\mapsto g(x_1)H_2$ defines a bijection, so $[G_1:
H_1]=[G_2: H_2]$). 

Recall from \reth{gorinlin}(\ref{it:gorinlin4}) that $(B_4(\dt))^{(1)}\cong \F[2](a,b)\rtimes \F[2](u,v)$, where $a=\sigma_3\sigma_1^{-1}$, $b= \sigma_2\sigma_3 \sigma_1^{-1} \sigma_2^{-1}$, $u=\sigma_2 \sigma_1^{-1}$, $v=\sigma_1 \sigma_2 \sigma_1^{-2}$, and the action is given by \req{gorinlin4}. The corresponding elements $u=\sigma_2 \sigma_1^{-1}$, $v=\sigma_1 \sigma_2 \sigma_1^{-2}$ of $B_3(\St\setminus \brak{x_1,x_2})$ (we use the same symbols to denote these elements) in fact belong to $(B_3(\St\setminus \brak{x_1,x_2}))^{(1)}$, generate a free group of rank~$2$ (by injectivity of $\phi$), and we have
\begin{equation}\label{eq:f2uvh}
\F[2](u,v)\subseteq \phi((B_3(\St\setminus
\brak{x_1,x_2}))^{(1)})\subseteq H. 
\end{equation}
Further, writing the elements of $(B_4(\dt))^{(1)}$ in the form $\F[2](a,b)\rtimes \F[2](u,v)$, if $(w,z)\in H$ then for all $z'\in \F[2](u,v)$, $z^{-1}z'\in H$ by \req{f2uvh}, so $(w,z)(e,z^{-1}z')= (w,z')\in H$, and thus $\brak{w}\times \F[2](u,v)\subseteq H$. Hence $H$ is of the form $H_1 \rtimes \F[2](u,v)$, where $H_1$ is an index~$4$ subgroup of $\F[2](a,b)$. Together with the identity, the elements $\pi(a)=(12)(34)$, $\pi(b)=(13)(24)$ and $\pi(ab)=(14)(23)$ form a set of coset representatives for $\an[4]/\an[3]$, so $e,a,b$ and $ab$ form a set of coset representatives for $\F[2](a,b)/H_1$. If $w=w(a,b)\in \F[2](a,b)$ then $\pi(w)=w(\pi(a), \pi(b))$, and since $\pi(a)$ and $\pi(b)$ generate a group isomorphic to $\Z_2 \times \Z_2$, we see that $\pi(w)\in A_3$ if and only if the exponent sums in $w$ of $a$ and $b$ are both even. In other words, $H_1= \ker{\F[2](a,b)\to \Z_2\times \Z_2}$, where $a\mapsto (1,0)$ and $b\mapsto (0,1)$, is nothing other than the free subgroup $N$ of $\F[2](a,b)$ of rank~$5$ described on page~\pageref{page:free5} of \resec{gorinlin}, possessing a basis $z_1=a^2$, $z_2=b^2$, $z_3=(ab)^2$, $z_4=ba^2b^{-1}$ and $z_5=ab^2a^{-1}$. Thus $H\cong \F[5](z_1,\ldots, z_5)\rtimes \F[2](u,v)$, where the action is given by equations~\reqref{acuxi} and~\reqref{acvxi}.

Hence we have a commutative diagram of the form:
\begin{equation}\label{eq:b3b4} 
\begin{xy}*!C\xybox{%
\entrymodifiers={!! <0pt, .8ex>+}
\xymatrix{%
& & & & \Z \ar[d]^{\psi}\\
\an[3] \ar@{^{(}->}[dd] & &\ar[ll]_-{\pi} \ar[dl]^{\phi} (B_{3,2})^{(1)} \ar@{^{(}->}[r] \ar[dd]^{\phi\bigr\rvert_{(B_{3,2})^{(1)}}}  & B_{3,2} \ar[r] \ar@{^{(}->}[dd]^{\phi} & {\underbrace{\gpab[({B_{3,2}})]}_{\cong\Z\times \Z}} \ar[dd]^{\widetilde{\phi}}\\
& H \ar[ul]_{\pi\bigr\rvert_H} \ar@{^{(}->}[dr] & & &\\
\an[4]  & &\ar[ll]_-{\pi} (B_4(\dt))^{(1)} \ar@{^{(}->}[r] & B_4(\dt) \ar[r] & {\underbrace{\gpab[{B_4(\dt)}]}_{\cong\Z}},}}
\end{xy}
\end{equation}
where $B_{3,2}= B_3(\St\setminus \brak{x_1,x_2})$, and 
\begin{equation*}
\map{\widetilde{\phi}}{\gpab[B_3(\St\setminus \brak{x_1,x_2})]}[{\gpab[B_4(\dt)]}]
\end{equation*}
is the homomorphism induced on the Abelianisations. Since $\gpab[{B_4(\dt)}]$ is infinite cyclic, generated by an element $\overline{\sigma}$, say, $\gpab[{B_3(\St\setminus \brak{x_1,x_2})}]$ is a free Abelian group of rank~$2$ with basis comprised of $\sigma$ and $A$, we have that $\widetilde{\phi}(\sigma)= \overline{\sigma}$, and $\widetilde{\phi}(A)=\overline{\sigma}^2$. So $\ker{\widetilde{\phi}} \cong \Z$ is the subgroup generated by $(-2,1)$ relative to the basis $(\sigma,A)$. Let 
\begin{equation*}
\map{\psi}{\Z}[{\gpab[B_3(\St\setminus \brak{x_1,x_2})]}]
\end{equation*}
be defined by $\psi(k)=k(-2,1)$. Then the final column of the diagram is exact.

Now the idea is the following: given $x\in H$, we may associate an element of $\Z$ using diagram~\reqref{b3b4}. We shall show that this is a homomorphism, $\epsilon'$ say, which satisfies $\epsilon'(z_i)=1$ if $i=1,2,3$, $\epsilon'(z_i)=-1$ if $i=4,5$, and $\epsilon'(u)= \epsilon'(v)=0$. From this, in particular, we obtain $\epsilon'=\eta\circ \epsilon$, where $\map{\epsilon}{\F[5](z_1,\ldots, z_5)\rtimes \F[2](u,v)}[\Z\oplus\Z^2]$ is the homomorphism defined by \req{defeps}, and $\map{\eta}{\Z\oplus \Z^2}[\Z]$ is defined by $\eta((1,0,0))=1$ and $\eta((0,1,0))=\eta((0,0,1))=0$.

To define $\epsilon'$, we first choose the following correspondence between $B_3(\St\setminus \brak{x_1,x_2})$ and $B_4(\dt)$: the three strings correspond of $B_3(\St\setminus \brak{x_1,x_2})$ to the first three strings of $B_4(\dt)$;  the puncture $x_1$ to the fourth (vertical) string; and the puncture $x_2$ to the boundary of $\dt$. In this representation, if $1\leq i\leq 3$ and $1\leq j\leq 2$ then $C_{i,j+3}$ will denote the element of $B_3(\St\setminus \brak{x_1,x_2})$ represented by a loop based at point $i$ which encircles the $j\up{th}$ puncture. Suppose that $z\in B_4(\dt)$ is such that $\pi(z)\in \sn[3]$. Then we claim that there exists $y_z\in B_3(\St\setminus \brak{x_1,x_2})$ such that $\phi(y_z)=z$; by injectivity of $\phi$, such a $y_z$ is unique. To prove the claim, notice there exists $y'\in B_3(\St\setminus \brak{x_1,x_2})$ such that $\pi(z)= \pi(y')= \pi\circ \phi(y')$ by commutativity of diagram~\reqref{p3st}. Hence $\phi(y')^{-1}z\in \ker{\pi}$. But since the first vertical arrow of that diagram is bijective, there exists $y''\in P_3(\St\setminus \brak{x_1,x_2})$ such that $\phi(y'')= \phi(y')^{-1}z$. Hence $z=\phi(y_z)$, where $y_z=y'y''$, and the claim is proved. 

So let $x\in H$. Since $H\subseteq (B_4(\dt))^{(1)}$, $x$ is sent to $0$ in $\gpab[{B_4(\dt)}]$. By the claim, there exists a unique $y_x \in B_3(\St\setminus \brak{x_1,x_2})$ such that $\phi(y_x)=x$, so by commutativity of diagram~\reqref{b3b4}, $\widetilde{y_x}\in \ker{\widetilde{\phi}}$, where $\widetilde{y_x}$ denotes the Abelianisation of $y_x$. Thus $\widetilde{y_x}= k(-2,1)$ relative to the basis $(\sigma,A)$, where $k\in\Z$, and so $x\mapsto k$ defines a map $\epsilon'$ from $H$ to $\Z$, well defined since $y_x$ is unique, and a homomorphism because $\phi$ is. Further, $\psi\circ \epsilon'(x)=\widetilde{y_x}$. Let us now calculate $\epsilon'$ on the given generating set $\brak{z_1,\cdots, z_5,u,v}$ of $H$. Now $y_u=u$ and $y_v=v$, and since $\F[2](u,v)\subseteq \phi(B_3(\St\setminus \brak{x_1,x_2}))^{(1)}$ by \reqref{f2uvh}, it follows that $\widetilde{y_u}=\widetilde{y_v}=0$, and so $\epsilon'(u)= \epsilon'(v)=0$. Now consider 
\begin{equation*}
x_1=a^2= (\sigma_3\sigma_1^{-1})^2= \sigma_3^2\sigma_1^{-2}.
\end{equation*}
From the given correspondence between $B_3(\St\setminus \brak{x_1,x_2})$ and $B_4(\dt)$, $z_1$ may be written as $C_{34}\sigma_1^{-2}$, and so under Abelianisation is sent to $(-2,1)$. Hence $\epsilon'(z_1)=1$. Since $z_2$ is conjugation of $z_1$ by $\sigma_2$, we obtain similarly that $\epsilon'(z_2)=1$. As for $z_3$,
\begin{align*}
z_3 &= (ab)^2= (\sigma_3\sigma_1^{-1}\cdot \sigma_2 \sigma_3\sigma_1^{-1} \sigma_2^{-1})^2\\
& = \sigma_1^{-1}\sigma_2 \sigma_3 \sigma_2 \sigma_1^{-1}\sigma_2^{-1} \sigma_1^{-1} \sigma_3 \sigma_2 \sigma_3 \sigma_1^{-1}\sigma_2^{-1}\\
&= \sigma_1^{-1}\sigma_2 \sigma_1^{-1}\sigma_3^2 \sigma_2 \sigma_1^{-1}\sigma_2^{-1}= \sigma_1^{-1}\sigma_2 \sigma_1^{-1} C_{34} \sigma_2 \sigma_1^{-1}\sigma_2^{-1},
\end{align*}
so $\epsilon'(z_3)=1$,
\begin{align*}
z_4 &=ba^2b^{-1}=\sigma_2 \sigma_3\sigma_1^{-1} \sigma_2^{-1} \sigma_1^{-2}\sigma_3^2 \sigma_2 \sigma_1\sigma_3^{-1} \sigma_2^{-1}\\
&=  \sigma_2 \sigma_3 \sigma_2^{-2} \sigma_1^{-1} \sigma_2^{-1} \sigma_3^2 \sigma_2 \sigma_1\sigma_3^{-1} \sigma_2^{-1}\\
&= \sigma_3^{-2} \sigma_2 \sigma_3 \sigma_1^{-1} \sigma_3 \sigma_2^2 \sigma_3^{-1} \sigma_1\sigma_3^{-1} \sigma_2^{-1}= \sigma_3^{-2} \sigma_2 \sigma_1^{-1} \sigma_3^2 \sigma_2^2 \sigma_1\sigma_3^{-2} \sigma_2^{-1}\\
&= C_{34}^{-1} \sigma_2 \sigma_1^{-1} C_{34} \sigma_2^2 \sigma_1 C_{34}^{-1} \sigma_2^{-1},
\end{align*}
so $\epsilon'(z_4)=-1$, and
\begin{align*}
z_5 &=ab^2a^{-1}= \sigma_1^{-1}\sigma_3 \sigma_2 \sigma_1^{-2}\sigma_3^2 \sigma_2^{-1} \sigma_3^{-1}\sigma_1\\
&= \sigma_1^{-1}\sigma_3 \sigma_2 \sigma_1^{-2} \sigma_2^{-1} \sigma_3^{-1} \sigma_2^2 \sigma_1\\
&= \sigma_1^{-2} \sigma_2^{-1}\sigma_3^{-2} \sigma_2 \sigma_1\sigma_2^2 \sigma_1= \sigma_1^{-2} \sigma_2^{-1}C_{34}^{-1} \sigma_2 \sigma_1\sigma_2^2 \sigma_1,
\end{align*}
so $\epsilon'(z_5)=-1$, and thus $\epsilon'=\eta\circ \epsilon$ as claimed.

Let $y\in (B_3(\St\setminus \brak{x_1,x_2}))^{(1)}$, and let $x=\phi(y)$. We know that $x\in H$, and $y_x=y$. But $\widetilde{y}=0$, hence $\epsilon'(x)=0$, and $\phi((B_3(\St\setminus \brak{x_1,x_2}))^{(1)})\subseteq \ker{\epsilon'}$. Conversely, let $x\in \ker{\epsilon'}$. Then $\phi(y_x)=x$, and $\widetilde{y_x}=\psi \circ \epsilon'(x)=0$, so $y_x\in (B_3(\St\setminus \brak{x_1,x_2}))^{(1)}$. Hence $x\in \phi((B_3(\St\setminus \brak{x_1,x_2}))^{(1)})$, and so $\ker{\epsilon'}=  \phi((B_3(\St\setminus \brak{x_1,x_2}))^{(1)})$. But $\phi$ is injective, and thus $\ker{\epsilon'}\cong (B_3(\St\setminus \brak{x_1,x_2}))^{(1)}$. To determine $\ker{\epsilon'}$, we use the following short exact sequence and the Reidemeister-Schreier rewriting process:
\begin{equation*}
1 \to (B_3(\St\setminus \brak{x_1,x_2}))^{(1)}\stackrel{\phi}{\lhra} \F[5](z_1,\ldots, z_5)\rtimes \F[2](u,v) \stackrel{\epsilon'}{\to} \Z\to 1.
\end{equation*}
The calculations are similar to those given in \resec{gorinlin} for the kernel of the homomorphism $\rho$ defined by \req{kerrho}; the difference is that here $\F[2](u,v) \subseteq \ker{\epsilon'}$. For all $j\in\Z$, set
\begin{equation*}
\alpha_{i,j}=
\begin{cases}
z_1^j z_i z_1^{-(j+1)} & \text{if $i=2,3$}\\
z_1^j z_i z_1^{-(j-1)} & \text{if $i=4,5$.}
\end{cases}
\end{equation*}
These elements form a basis of $\ker{\rho}$ (see Table~\ref{tab:basis} on page~\pageref{tab:basis}). To obtain a generating set of $\ker{\epsilon'}$, we need to add the elements $r_j=z_1^j uz_1^{-j}$ and $s_j= z_1^j vz_1^{-j}$, where $j\in\Z$. The relators of $\ker{\epsilon'}$ are of the form $z_1^j \mathcal{R} z_1^{-j}$, where $j\in\Z$ and $\mathcal{R}$ runs over the set of relators given by equations~\reqref{acuxi} and~\reqref{acvxi}. For $i=1,\ldots, 5$, let us set $t_i=uz_iu^{-1}$ and $w_i=vz_iv^{-1}$. First let $i=1$. Then for all $j\in\Z$, we have the relator:
\begin{equation*}
z_1^j uz_1 u^{-1}t_1^{-1} z_1^{-j}= r_jr_{j+1}^{-1}(z_1^j t_1^{-1} z_1^{-(j+1)})^{-1}.
\end{equation*}
from which it follows that 
\begin{equation}\label{eq:rj1rj} 
r_{j+1}=(z_1^j t_1 z_1^{-(j+1)})^{-1} r_j.
\end{equation}
This allows us to delete all of the $r_i$, $i\in\Z \setminus\brak{0}$, from the list of generators. By induction, we obtain:
\begin{equation}\label{eq:rj}
r_j=\begin{cases}
(z_1^{j-1} t_1 z_1^{-j})^{-1} (z_1^{j-2} t_1 z_1^{-j+1})^{-1} \cdots (t_1 z_1^{-1})^{-1}u & \text{if $j>0$}\\
(z_1^j t_1 z_1^{-(j+1)}) (z_1^{j+1} t_1 z_1^{-(j+2)}) \cdots (z_1^{-1} t_1) u  & \text{if $j< 0$.}
\end{cases}
\end{equation}
In a similar way, we may delete all of the $s_j$ except $s_0=v$ from the list of generators, and we obtain:
\begin{equation}\label{eq:sj}
s_j=
\begin{cases}
(z_1^{j-1} w_1 z_1^{-j})^{-1} (z_1^{j-2} w_1 z_1^{-j+1})^{-1} \cdots (w_1 z_1^{-1})^{-1}v & \text{if $j>0$}\\
(z_1^j w_1 z_1^{-(j+1)}) (z_1^{j+1} w_1 z_1^{-(j+2)}) \cdots (z_1^{-1} w_1) v &\text{if $j< 0$.}
\end{cases}
\end{equation}
Notice that in equations~\reqref{rj} and~\reqref{sj}, each of the bracketed terms belongs to $\ker{\rho}$, and hence so do $r_j u^{-1}$ and $s_j v^{-1}$. So they may be expressed in terms of the $\alpha_{i,j}$. Now let $i=2,3$ and $j\in\Z$. Then we have a relator
\begin{equation*}
z_1^j uz_iu^{-1}t_i^{-1} z_1^{-j}= r_j \alpha_{i,j} r_{j+1}^{-1} z_1^{j+1}t_i^{-1}z_1^{-j},
\end{equation*}
which yields a relation of the form $r_j \alpha_{i,j} r_j^{-1}= z_1^j t_it_1^{-1}z_1^{-j}$ by \req{rj1rj}. Using \req{rj}, we see that the elements $u z_1^jz_iz_1^{-(j+1)} u^{-1}$ may be expressed solely in terms of the $\alpha_{i,j}$. Indeed, 
\begin{equation}\label{eq:u23}
u \alpha_{i,j} u^{-1}= (r_j u^{-1})^{-1} z_1^j t_it_1^{-1}z_1^{-j} (r_j u^{-1}).
\end{equation} 
Similarly,
\begin{equation}\label{eq:v23}
v \alpha_{i,j} v^{-1}= (s_j v^{-1})^{-1} z_1^j w_iw_1^{-1}z_1^{-j} (s_j v^{-1}).
\end{equation} 
Finally, let $i=4,5$ and $j\in\Z$. Then we obtain analogously:
\begin{equation}\label{eq:u45}
u \alpha_{i,j} u^{-1}= (r_j u^{-1})^{-1} z_1^j t_it_1 z_1^{-j} (r_j u^{-1}).
\end{equation} 
Similarly,
\begin{equation}\label{eq:v45}
v \alpha_{i,j} v^{-1}= (s_j v^{-1})^{-1} z_1^j w_iw_1 z_1^{-j} (s_j v^{-1}).
\end{equation} 
This gives a complete set of relations for $\ker{\epsilon'}$. We conclude that:
\begin{prop}\label{prop:f5f2}
$(B_3(\St\setminus \brak{x_1,x_2}))^{(1)}\cong L \rtimes \F[2](u,v)$, where $L$ is the subgroup of $\F[5](z_1,\ldots, z_5)$ of infinite rank freely generated by $\brak{\alpha_{i,j}}_{i\in\Z,\, j\in\brak{2,3,4,5}}$, the action being given by equations~\reqref{u23},~\reqref{v23},~\reqref{u45}, and~\reqref{v45}, taking into account equations~\reqref{rj} and~\reqref{sj}.\qed
\end{prop}

From this, we may now determine $(B_3(\St\setminus \brak{x_1,x_2}))^{(1)}/ (B_3(\St\setminus \brak{x_1,x_2}))^{(2)}$ by Abelianising the presentation of $(B_3(\St\setminus \brak{x_1,x_2}))^{(1)}$ given by \repr{f5f2}, and thus reprove \repr{b3z4}. First notice that for all $2\leq i\leq 5$ and all $j,k\in\Z$, $z_1^k \alpha_{i,j} z_1^{-k}=\alpha_{i,j+k}$. If $w\in (B_3(\St\setminus \brak{x_1,x_2}))^{(1)}$, let $\widetilde{w}\in (B_3(\St\setminus \brak{x_1,x_2}))^{(1)}/(B_3(\St\setminus \brak{x_1,x_2}))^{(2)}$ denote its Abelianisation. A simple calculation shows that:
\begin{align*}
\widetilde{t_2t_1^{-1}} &= \widetilde{\alpha_{2,0}} + \widetilde{\alpha_{5,0}} + \widetilde{\alpha_{3,-1}} - \widetilde{\alpha_{2,-1}}- \widetilde{\alpha_{4,0}}\\
\widetilde{t_3t_1^{-1}} &= \widetilde{\alpha_{2,1}}- \widetilde{\alpha_{3,1}}+ \widetilde{\alpha_{5,1}} +\widetilde{t_2t_1^{-1}}\\
\widetilde{w_2w_1^{-1}} &= - \widetilde{\alpha_{4,1}} + \widetilde{\alpha_{2,1}}- \widetilde{\alpha_{3,1}}+ \widetilde{\alpha_{5,1}}  + \widetilde{\alpha_{5,0}}- \widetilde{\alpha_{4,0}}+ \widetilde{\alpha_{2,0}}\\
\widetilde{w_3w_1^{-1}} &= \widetilde{\alpha_{2,1}}+ \widetilde{w_2w_1^{-1}}.
\end{align*}
Abelianising equations~\reqref{u23} and~\reqref{v23} for $i=2$ then $i=3$ yields:
\begin{align}
\widetilde{\alpha_{5,j}}- \widetilde{\alpha_{4,j}} &= \widetilde{\alpha_{2,j-1}}- \widetilde{\alpha_{3,j-1}}\notag\\
\widetilde{\alpha_{3,j}}- \widetilde{\alpha_{2,j}} &= \widetilde{\alpha_{2,j+1}}- \widetilde{\alpha_{3,j+1}}+ \widetilde{\alpha_{5,j+1}}\notag\\
0= &\widetilde{\alpha_{5,j}}- \widetilde{\alpha_{4,j}} +\widetilde{\alpha_{5,j+1}}- \widetilde{\alpha_{4,j+1}} + \widetilde{\alpha_{2,j+1}}- \widetilde{\alpha_{3,j+1}}\notag\\
\widetilde{\alpha_{2,j+1}}+\widetilde{\alpha_{2,j}} &= \widetilde{\alpha_{3,j}}\label{eq:equa4}
\end{align}
for all $j\in\Z$. Substituting \req{equa4} into the three other equations, we obtain:
\begin{align}
\widetilde{\alpha_{5,j}}- \widetilde{\alpha_{4,j}} &= -\widetilde{\alpha_{2,j}}\label{eq:equa1}\\
 \widetilde{\alpha_{3,j+1}}&= \widetilde{\alpha_{5,j+1}}\label{eq:equa2}\\
\widetilde{\alpha_{2,j+1}}-\widetilde{\alpha_{2,j}} &= \widetilde{\alpha_{4,j+1}}.\label{eq:equa3}
\end{align}
Similarly, if $i=4,5$,
\begin{align*}
\widetilde{t_4t_1} =& \widetilde{\alpha_{2,0}}+ \widetilde{\alpha_{5,0}}\\
\widetilde{t_5t_1} =& \widetilde{\alpha_{2,0}}+ \widetilde{\alpha_{2,1}} -\widetilde{\alpha_{3,1}}+ \widetilde{\alpha_{5,1}}+ \widetilde{\alpha_{5,0}}+ \widetilde{\alpha_{2,-1}}\\
\widetilde{w_4w_1} =& \widetilde{\alpha_{3,0}}- \widetilde{\alpha_{4,1}} +\widetilde{\alpha_{2,1}}- \widetilde{\alpha_{3,1}}+ \widetilde{\alpha_{5,1}}\\
\widetilde{w_5w_1} =& \widetilde{\alpha_{3,-1}}- \widetilde{\alpha_{4,0}} +\widetilde{\alpha_{2,0}}- \widetilde{\alpha_{4,1}}+ 2\widetilde{\alpha_{2,1}} -\widetilde{\alpha_{3,1}}+ \widetilde{\alpha_{5,1}} +\\
&\widetilde{\alpha_{5,0}}+ \widetilde{\alpha_{3,-2}}- \widetilde{\alpha_{2,-2}}   -\widetilde{\alpha_{4,-1}}.
\end{align*}

Abelianising equations~\reqref{u45} and~\reqref{v45} for $i=4$ then $i=5$ and applying the previous equations yields:
\begin{align}
\widetilde{\alpha_{4,j}} &= \widetilde{\alpha_{2,j}}+\widetilde{\alpha_{5,j}} \; \text{which is equivalent to \reqref{equa1}}\notag\\
0  &= \widetilde{\alpha_{2,j-1}}+\widetilde{\alpha_{2,j}}+ \widetilde{\alpha_{2,j+1}} \label{eq:equa5}\\
\widetilde{\alpha_{4,j}} &=\widetilde{\alpha_{3,j}}+  \widetilde{\alpha_{2,j}}\label{eq:equa6}\\
0  &= \widetilde{\alpha_{2,j-1}}+\widetilde{\alpha_{2,j}}+ \widetilde{\alpha_{2,j+1}}\; \text{which is the same as \reqref{equa5}.}\notag
\end{align}
By equations~\reqref{equa5},~\reqref{equa4},~\reqref{equa3} and~\reqref{equa2} we obtain the solution
\begin{align*}
\widetilde{\alpha_{2,3k}}&= \widetilde{\alpha_{2,0}}= -\widetilde{\alpha_{3,3k+1}} = -\widetilde{\alpha_{5,3k+1}}\\
\widetilde{\alpha_{2,3k+1}}&= \widetilde{\alpha_{2,1}}= -\widetilde{\alpha_{3,3k+2}} = -\widetilde{\alpha_{5,3k+2}}\\
\widetilde{\alpha_{2,3k+2}}&= -(\widetilde{\alpha_{2,0}}+ \widetilde{\alpha_{2,1}})= -\widetilde{\alpha_{3,3k}} = -\widetilde{\alpha_{5,3k}}\\
\widetilde{\alpha_{4,3k}}&= 2\widetilde{\alpha_{2,0}}+ \widetilde{\alpha_{2,1}}\\
\widetilde{\alpha_{4,3k+1}}&= -\widetilde{\alpha_{2,0}}+ \widetilde{\alpha_{2,1}}\\
\widetilde{\alpha_{4,3k+2}}&= -\widetilde{\alpha_{2,0}}-2\widetilde{\alpha_{2,1}}
\end{align*}
for all $k\in\Z$, which satisfies the two remaining equations~\reqref{equa1} and~\reqref{equa6}. We conclude that $(B_3(\St\setminus \brak{x_1,x_2}))^{(1)}/ (B_3(\St\setminus \brak{x_1,x_2}))^{(2)}$ is a free Abelian group with basis $\brak{\widetilde{\alpha_{2,0}}, \widetilde{\alpha_{2,1}}, \widetilde{u}, \widetilde{v}}$, and this reproves \repr{b3z4}.

We will come back to the special case $m=n=2$ in the following section.

\section{The lower central and derived series of $B_2(\St\setminus\brak{x_1,x_2})$}\label{sec:lcgsb22}

From~\resec{bmn1n} and \reth{lcdsbmsn}, the only outstanding case for
the lower central series of the punctured sphere is the $2$-string
braid group. As we shall see in this section, it is particularly
challenging. We concentrate here on the case of the two-punctured
sphere. The group $B_2(\St\setminus\brak{x_1,x_2})$ has many
fascinating properties, and as a result, we are able to describe its
lower central and derived series in terms of those of the free product
$\Z_2\ast \Z$. In particular, we prove \reco{b2b2}, \repr{b2b2resid}
and \reth{lcsb2}. As we indicated in the Preface, the techniques used
in this section have since been applied in~\cite{BGG} to study the
$2$-string braid group of the torus, and similar results were obtained
(cf.\ \reth{bggtor}).

We start by determining $\Gamma_2(P_2(\St\setminus\brak{x_1,x_2}))$.
The map $F_2(\St\setminus\brak{x_1,x_2})\to
F_1(\St\setminus\brak{x_1,x_2})$ is a fibration, and the fibre over a
point $x_3$ of the base is of the form
$F_1(\St\setminus\brak{x_1,x_2,x_3})$. As in \req{split}, this gives
rise to a short exact sequence of the form:
\begin{multline*}
1\to \pi_1(\St\setminus\brak{x_1,x_2,x_3},x_4)\to P_2(\St\setminus\brak{x_1,x_2},\brak{x_3,x_4}) \stackrel{p_{\ast}}{\to}\\ 
\pi_1(\St\setminus\brak{x_1,x_2},x_3)\to 1.
\end{multline*}
We use the following notation for the generators $\gamma_{i,j}$,
$1\leq i,j \leq 2$ of
$P_2(\St\setminus\brak{x_1,x_2},\brak{x_3,x_4})$: the two punctures
correspond to the points $x_1,x_2$; and the two basepoints correspond
to $x_3,x_4$. The generator $\gamma_{i,j}$ is equal to the generator
$A_{i,j+2}$ of \repr{presbetanm}, and corresponds to a loop based at
$x_{j+2}$ which encircles $x_i$ in the positive direction. 

Let $\sigma$ be the standard Artin generator of
$B_2(\St\setminus\brak{x_1,x_2})$ which geometrically exchanges the
points $x_3$ and $x_4$. Then a (non-minimal) generating set of
$P_2(\St\setminus\brak{x_1,x_2})$ is given by the union of the
$\gamma_{i,j}$ and $\sigma^2$, and a generating set of
$B_2(\St\setminus\brak{x_1,x_2})$ is given by the union of the
$\gamma_{i,j}$ and $\sigma$.

The kernel $\pi_1(\St\setminus\brak{x_1,x_2,x_3},x_4)$ of $p_{\ast}$ is the free group $\F[2](a,b)$ of rank~$2$ on $a$ and $b$, where $a=\gamma_{1,2}$ and $b=\gamma_{2,2}$. The image of $p_{\ast}$ is an infinite cyclic group; the homomorphism which sends (one of) its generators to the element $c=\gamma_{2,1}$ of $P_2(\St\setminus\brak{x_1,x_2})$ defines a section for $p_{\ast}$. Hence
\begin{equation}\label{eq:p2semi}
P_2(\St\setminus\brak{x_1,x_2})\cong \F[2](a,b)\rtimes_{\phi} \Z,
\end{equation} 
where we identify the second factor $\Z$ with $\ang{c}$. The action $\phi$ on the kernel is given as follows (this may be checked using the presentation of $P_m(\St\setminus\brak{x_1,\ldots,x_n})$ given in~\cite{GG4}):
\begin{equation*}
\left\{ 
\begin{aligned}
\phi(c)(a)=cac^{-1} &= a\\
\phi(c)(b)= cbc^{-1} &= aba^{-1},
\end{aligned}\right.
\end{equation*}
which in fact is just conjugation by $a$.

As well as containing $\Gamma_2(\F[2](a,b))$, by \repr{gammasemi},
$\Gamma_2(P_2(\St\setminus\brak{x_1,x_2}))$ will also contain
elements of the form $[c^j,w]$, where $w\in \F[2](a,b)$ and $j\in\Z$. But from the form of the action $\phi$, $[c^j,w]=[a^j,w]$. Hence 
\begin{equation}\label{eq:gam2p2}
\Gamma_2(P_2(\St\setminus\brak{x_1,x_2}))= \Gamma_2(\F[2](a,b)),
\end{equation}
and thus the derived series (with the exception of the first term) of $P_2(\St\setminus\brak{x_1,x_2})$ is that of $\F[2](a,b)$.

By \repr{presbetanm}, $B_2(\St\setminus\brak{x_1,x_2})$ is generated by the $\gamma_{i,j}$, $1\leq i,j\leq 2$, and $\sigma$, subject to the four relations:
\begin{equation}\label{eq:b2s2}
\left.
\begin{aligned}
& \gamma_{1,2}\gamma_{2,2}\sigma^2=1\\
& \text{$\gamma_{1,1}\gamma_{2,1}\sigma^2=1$}\\
& \text{$\sigma \gamma_{i,1}\sigma^{-1}= \gamma_{i,2}$ for $i=1,2$.\;}
\end{aligned}\right\}
\end{equation}
Thus:
\begin{equation}\label{eq:abc}
\left.
\begin{aligned}
c &=\gamma_{2,1}=\gamma_{1,1}^{-1}\sigma^{-2}\\
b &=\gamma_{2,2}=\sigma \gamma_{1,1}^{-1} \sigma^{-3}\\
a &=\gamma_{1,2}= \sigma \gamma_{1,1}\sigma^{-1}.
\end{aligned}\right\}
\end{equation}
In particular, $B_2(\St\setminus\brak{x_1,x_2})$ is generated by $\sigma$ and $\gamma_{1,1}$.

In what follows, we shall sometimes write simply $P_{2,2}$ for $P_2(\St\setminus\brak{x_1,x_2})$, and $B_{2,2}$ for $B_2(\St\setminus\brak{x_1,x_2})$. 

\begin{prop}\mbox{}\label{prop:p2b2}
\begin{enumerate}
\item The commutator subgroup $[P_{2,2},P_{2,2}]$ of $P_{2,2}$ is a normal subgroup of  $[P_{2,2},B_{2,2}]$. 
\item\label{it:it2} The commutator subgroup $[P_{2,2},B_{2,2}]$ of $P_{2,2}$ and $B_{2,2}$ is a normal subgroup of $P_{2,2}$ and $B_{2,2}$.
\item The quotient group $[P_{2,2},B_{2,2}]\bigl/[P_{2,2},P_{2,2}]\bigr.$ is isomorphic to $\Z$, and is generated by the coset of the element $[\sigma,b]=b^{-1}c$.
\end{enumerate}
\end{prop}

\begin{proof}\mbox{}
\begin{enumerate}
\item This is clear since $P_{2,2}\triangleleft B_{2,2}$.
\item The fact that $[P_{2,2},B_{2,2}]$ is a subgroup of $P_{2,2}$ follows from projection into the symmetric group $\sn[2]$. Since $P_{2,2}\triangleleft B_{2,2}$, we see that $[P_{2,2},B_{2,2}]\triangleleft B_{2,2}$, and so $[P_{2,2},B_{2,2}]\triangleleft P_{2,2}$.
\item Using \req{abc}, we see easily that $\sigma^{-2}=ab$, and thus:
\begin{equation}\label{eq:sigma}
\left.\begin{aligned}\mbox{}
[\sigma,a] &= \sigma^2 \sigma^{-2} c^{-1} \sigma^{-2} a^{-1}= c^{-1}aba^{-1}=bc^{-1}=b(b^{-1}c)^{-1} b^{-1}\\
[\sigma,b] &= \sigma^2 \gamma_{1,1}^{-1} \sigma^{-4} b^{-1}= b^{-1}a^{-1}ca=b^{-1}c\\
[\sigma,c] &= \sigma \gamma_{1,1}^{-1} \sigma^{-1} \sigma^2 \gamma_{1,1}= b(b^{-1}c)^{-1} b^{-1}.
\end{aligned}\right\}
\end{equation}

We know that $[P_{2,2},B_{2,2}]$ is the normal closure in $B_{2,2}$ of the set of elements of the form $[\rho_1,\rho_2]$ and their inverses, where $\rho_1\in \brak{a,b,c}$ is a generator of $P_{2,2}$, and $\rho_2\in \brak{\sigma, a,b,c}$ is a generator of $B_{2,2}$. If $\rho_2\in \brak{a,b,c}$ then $[\rho_1,\rho_2]\in [P_{2,2},P_{2,2}]$. So we just need to consider the cosets of the conjugates of elements of the form $[\rho_1,\sigma]$. Consider the following relation:
\begin{equation}\label{eq:gam1}
\rho [\rho_1,\sigma] \rho^{-1}=\bigl[\rho, [\rho_1,\sigma]\bigr] \bigl[\rho_1,\sigma\bigr].
\end{equation}
If $\rho\in P_{2,2}$, then since $[\rho_1,\sigma]\in P_{2,2}$ by~(\ref{it:it2}), it follows that 
\begin{equation*}
\rho [\rho_1,\sigma] \rho^{-1} \equiv [\rho_1,\sigma] \quad \text{modulo $[P_{2,2},P_{2,2}]$.}
\end{equation*}
So suppose that $\rho\in B_{2,2}\setminus P_{2,2}$. Then $w=\rho\sigma^{-1}\in P_{2,2}$. By \req{sigma}, we see that
\begin{align*}
\sigma [\sigma,a]\sigma^{-1} &= \sigma bc^{-1}\sigma^{-1}= b^{-1} c b\cdot b^{-1}= b^{-1} [\sigma,a]^{-1}b\\
\sigma [\sigma,b]\sigma^{-1} &= \sigma b^{-1}c\sigma^{-1}=  b^{-1} c^{-1} b \cdot b = b^{-1} [\sigma,b]^{-1}b,\quad\text{and}\\
\sigma [\sigma,c]\sigma^{-1} &= \sigma [\sigma,a]\sigma^{-1}= b^{-1} [\sigma,c]^{-1}b.
\end{align*}
In other words, for all $\rho_1\in \brak{a,b,c}$, we have 
\begin{equation*}
\sigma [\sigma,\rho_1]\sigma^{-1}= b^{-1} [\sigma,\rho_1]^{-1} b.
\end{equation*}
Hence 
\begin{align*}
\bigl[\rho, [\rho_1,\sigma]\bigr] &= w \sigma [\sigma,\rho_1]\sigma^{-1} w^{-1} [\rho_1,\sigma]^{-1}\\
&=
w b^{-1} [\sigma,\rho_1]^{-1}b w^{-1} [\rho_1,\sigma]^{-1} = \bigl[
w b^{-1}, [\sigma,\rho_1]^{-1}\bigr] \bigl[ \rho_1, \sigma\bigr]^{-2}.
\end{align*}
Thus by \req{gam1}, we obtain
\begin{align*}
\rho [\rho_1,\sigma] \rho^{-1} &=\bigl[
w b^{-1}, [\sigma,\rho_1]^{-1}\bigr] \bigl[\rho_1,\sigma \bigr]^{-1}\\
& \equiv \text{$[\rho_1,\sigma]^{-1}$ modulo $[P_{2,2},P_{2,2}]$, since $[\sigma,\rho_1]\in P_{2,2}$.}
\end{align*}
By \req{sigma}, 
\begin{equation*}
\text{$[\sigma,a]\equiv [\sigma,c]\equiv  [\sigma,b]^{-1}$ modulo $[P_{2,2},P_{2,2}]$,}
\end{equation*}
and since $[B_{2,2},P_{2,2}]=[P_{2,2},B_{2,2}]$, we conclude that the quotient  $[P_{2,2},B_{2,2}]\bigl/[P_{2,2},P_{2,2}]\bigr.$ is infinite cyclic, and generated by the coset of the element $[\sigma,b]= b^{-1}c$ (using equations~\reqref{p2semi} and~\reqref{gam2p2}, one may check that $b^{-1}c\notin [P_{2,2},P_{2,2}]$).\qedhere
\end{enumerate}
\end{proof}

\begin{rem}\label{rem:stall}
Let us give an alternative proof of \repr{p2b2} using Stallings' exact sequence~\reqref{stallings}. Since $[P_{2,2},P_{2,2}], [P_{2,2},B_{2,2}] \trianglelefteq P_{2,2}$ and $[P_{2,2},P_{2,2}]\subseteq [P_{2,2},B_{2,2}]$, we see that \begin{equation*}
[P_{2,2},B_{2,2}]\bigl/[P_{2,2},P_{2,2}]\bigr. \trianglelefteq P_{2,2}\bigl/[P_{2,2},P_{2,2}]\bigr..
\end{equation*}
We thus have the following diagram:
\begin{equation*}
\xymatrix@C=0.45cm{%
              &     &        & 1 \ar[d] & \\
1  \ar[r]    &  [P_{2,2},B_{2,2}]\bigl/[P_{2,2},P_{2,2}]\bigr.  \ar[r]  &  P_{2,2}\bigl/[P_{2,2},P_{2,2}]\bigr.  \ar[r] & P_{2,2}\bigl/[P_{2,2},B_{2,2}]\bigr. \ar[d]   \ar[r] & 1\\
 &  &  &  B_{2,2}\bigl/[B_{2,2},B_{2,2}]\bigr. \ar[d]  & \\
 & & & \sn[2] \ar[d] & \\
&                       &       & 1. &}
\end{equation*}
The vertical short exact sequence is that of Stallings applied to the
usual exact sequence $1\to P_{2,2}\to B_{2,2}\to \sn[2]\to 1$. By
\repr{abbm}, $B_{2,2}\bigl/[B_{2,2},B_{2,2}]\bigr.$ is a free Abelian
group of rank~$2$, with basis $\brak{\sigma, \gamma_{1,1}}$
(notationally, here we do not distinguish an element of $B_{2,2}$ and
its Abelianisation). The kernel $P_{2,2}\bigl/[P_{2,2},B_{2,2}]\bigr.$
of the projection $B_{2,2}\bigl/[B_{2,2},B_{2,2}]\bigr.\to \sn[2]$
certainly contains the free subgroup of rank~$2$ with basis
$\brak{\sigma^2, \gamma_{1,1}}$, and in fact is equal to this subgroup
(for otherwise it would contain an element of the form $\sigma^p
\gamma_{1,1}^q$, where $p,q\in\Z$ and $p$ is odd, and thus would
contain $\sigma$, which is clearly not in the kernel). Since
$P_{2,2}\bigl/[P_{2,2},P_{2,2}]\bigr.$ is isomorphic to $\Z^3$ (by
\req{p2semi}), we see that the kernel
$[P_{2,2},B_{2,2}]\bigl/[P_{2,2},P_{2,2}]\bigr.$ of the horizontal
exact sequence is isomorphic to $\Z$. Further,
$P_{2,2}\bigl/[P_{2,2},P_{2,2}]\bigr.$ is freely generated by $a,b$
and $c$. From the relation $\sigma \gamma_{2,1}\sigma^{-1}=
\gamma_{2,2}$, we see that $b=[\sigma,c]\cdot c$, and so $b$ and $c$
project to the same element in $P_{2,2}\bigl/[P_{2,2},B_{2,2}]\bigr.$.
Hence (the coset of) $bc^{-1}$ is a non-trivial element of
$[P_{2,2},B_{2,2}]\bigl/[P_{2,2},P_{2,2}]\bigr.$, which yields the
result.
\end{rem}

We thus obtain a short exact sequence of the form:
\begin{equation*}
1\to [P_{2,2},P_{2,2}]\to [P_{2,2},B_{2,2}]\to \Z\to 1,
\end{equation*}
for which the homomorphism $\map {s}{\Z}[[P_{2,2},B_{2,2}]]$ defined by $s(1)= b^{-1}c$ defines a splitting. Since $[P_{2,2},P_{2,2}]$ is the normal closure in $P_{2,2}$ of the set of elements of the form $[\rho_1,\rho_2]$ and their inverses, where $\rho_1,\rho_2\in\brak{a,b,c}$, and the action in $\F[2](a,b)$ of conjugation by $c$ is just conjugation by $a$, we see that $[P_{2,2},P_{2,2}]$ is the normal closure in $\F[2](a,b)$ of the element $[a,b]$, and that:
\begin{equation*}
(b^{-1}c) w  (b^{-1}c)^{-1}= (b^{-1}a) w (b^{-1}a)^{-1}\quad\text{for all $w\in \F[2](a,b)$}.
\end{equation*}
Hence the action of $b^{-1}c$ on $[P_{2,2},P_{2,2}]$ is that of conjugation by $b^{-1}a$, and so by the above short exact sequence,
\begin{align*}
[P_{2,2},B_{2,2}] &\cong [P_{2,2},P_{2,2}] \rtimes_{\psi} \Z\\
&\cong \Gamma_2(\F[2](a,b)) \rtimes_{\psi} \Z\quad\text{by \req{gam2p2}},
\end{align*}
where the action $\psi$ of $\Z$ on $\Gamma_2(\F[2](a,b))$ is given by conjugation by $b^{-1}a$.

\begin{prop}\label{prop:b2b2}
$[P_{2,2},B_{2,2}]=[B_{2,2},B_{2,2}]$.
\end{prop}

From this, it follows immediately that:
\begin{varthm}[\reco{b2b2}]
$\Gamma_2(B_2(\St\setminus\brak{x_1,x_2})) \cong \Gamma_2(\F[2](a,b)) \rtimes_{\psi} \Z$. \qed
\end{varthm}

\begin{proof}[Proof of \repr{b2b2}] 
Consider the following commutative diagram of short exact sequences (obtained by taking the first two vertical sequences, and the second and third horizontal sequences, and then completing to the whole diagram):
\begin{equation*}
\xymatrix@C=0.4cm{%
              &  1 \ar[d]   &    1 \ar[d]    & 1 \ar[d] & \\
1  \ar[r]     &  [P_{2,2},B_{2,2}] \ar[d]   \ar[r]  &  [B_{2,2},B_{2,2}] \ar[d]  \ar[r] & [B_{2,2},B_{2,2}]\bigl/[P_{2,2},B_{2,2}]\bigr. \ar[d]   \ar[r] & 1\\
1 \ar[r] & P_{2,2} \ar[r] \ar[d] & B_{2,2} \ar[r] \ar[d] & \Z_2 \ar[r] \ar@{=}[d] & 1\\
1 \ar[r] & P_{2,2}\bigl/[P_{2,2},B_{2,2}]\bigr. \ar[r]  \ar[d] & B_{2,2}\bigl/[B_{2,2},B_{2,2}]\bigr. \ar[r] \ar[d] & \Z_2 \ar[r] \ar[d] & 1.\\
&  1                     &    1   & 1 &}
\end{equation*}
As in \rerem{stall}, the third row is Stallings' exact sequence~\reqref{stallings} applied to the second row. By exactness of the third vertical sequence, it follows that $[B_{2,2},B_{2,2}]=[P_{2,2},B_{2,2}]$.
\end{proof}

We may obtain an alternative description of $\Gamma_2(B_2(\St\setminus\brak{x_1,x_2}))$ as a free group of infinite rank. To see this, notice from part~(\ref{it:isob}) of \repr{iso} that $B_m(\St\setminus\brak{x_1,x_2})\cong B_m(\dt\setminus \brak{x_2})$, and from part~(\ref{it:isoc}) that 
\begin{equation*}
B_m(\dt\setminus \brak{x_2}, \brak{x_1,x_3,\ldots, x_{m+1}})\cong B_{m,1}(\dt).
\end{equation*}
Hence $B_2(\St\setminus\brak{x_1,x_2})\cong B_{2,1}(\dt)$. But from part~(\ref{it:isod}), 
\begin{align}
B_{2,1}(\dt)&\cong \pi_1(\dt\setminus\brak{x_3,x_4},x_2)\rtimes B_2(\dt)\nonumber\\
&\cong \F[2](\gamma_{2,1}, \gamma_{2,2}) \rtimes_{\phi} \ang{\sigma},\label{eq:b2semi}
\end{align}
where the action, obtained from equations~\reqref{b2s2},~\reqref{abc} and~\reqref{sigma}, is given by:
\begin{equation}\label{eq:actsig}
\left.
\begin{aligned}
\phi(\sigma)(\gamma_{2,1})&= \gamma_{2,2}\\
\phi(\sigma)(\gamma_{2,2})&= \gamma_{2,2}^{-1}\gamma_{2,1} \gamma_{2,2}.
\end{aligned}
\right\}
\end{equation}
So if $w=w(\gamma_{2,1}, \gamma_{2,2})\in \F[2](\gamma_{2,1}, \gamma_{2,2})$ then 
\begin{equation}\label{eq:sigw}
\phi(\sigma)(w)=\gamma_{2,2}^{-1}w(\gamma_{2,2}, \gamma_{2,1}) \gamma_{2,2},
\end{equation}
in other words, the action consists of exchanging $\gamma_{2,1}$ and $\gamma_{2,2}$, then conjugating by $\gamma_{2,2}^{-1}$. Let $N$ denote the normal closure in $\F[2](\gamma_{2,1}, \gamma_{2,2})$ of the elements of the form $\phi(\sigma^j)(w)\cdot w^{-1}$, where $j\in\Z$ and $w\in \F[2](\gamma_{2,1}, \gamma_{2,2})$, and let
$L$ be the subgroup of $\F[2](\gamma_{2,1}, \gamma_{2,2})$ generated by $\Gamma_2(\F[2](\gamma_{2,1}, \gamma_{2,2}))$ and $N$. Then it follows from \repr{gammasemi} and \req{b2semi} that $\Gamma_2(\gamma_{2,1}(\dt))\cong L$. 

\begin{prop}\label{prop:kerpsi}\mbox{}
\begin{enumerate}
\item $L$ is the kernel of the homomorphism $\map{\psi}{\F[2](\gamma_{2,1}, \gamma_{2,2})}[\Z]$, where $\psi$ is augmentation.
\item $L$ is a free group of infinite rank with basis $\brak{z_i}_{i\in\Z}$, where $z_i= \gamma_{2,1}^i \gamma_{2,2} \gamma_{2,1}^{-(i+1)}$ for all $i\in\Z$.
\end{enumerate}
\end{prop}

Since $L\cong  \Gamma_2(B_{2,1}(\dt)) \cong B_2(\St\setminus\brak{x_1,x_2})$, we obtain immediately:
\begin{cor}\label{cor:gam2b2s2}
$\Gamma_2(B_2(\St\setminus\brak{x_1,x_2}))$ is a free group of infinite rank with basis $\brak{z_i}_{i\in\Z}$, where $z_i= \gamma_{2,1}^i \gamma_{2,2} \gamma_{2,1}^{-(i+1)}$ for all $i\in\Z$.\qed
\end{cor}

\begin{proof}[Proof of \repr{kerpsi}]
First observe that $\psi$ factors through the Abelianisation of $\F[2](\gamma_{2,1}, \gamma_{2,2})$, and so $\Gamma_2(\F[2](\gamma_{2,1}, \gamma_{2,2}))\subseteq \ker{\psi}$. Secondly, from \req{actsig}, $\sigma$ commutes with $cb=\gamma_{2,1}\gamma_{2,2}$, and it follows from \req{sigw} that
\begin{multline*}
\sigma^m w(\gamma_{2,1}, \gamma_{2,2}) \sigma^{-m}=\\
\begin{cases}
(\gamma_{2,1}\gamma_{2,2})^{-m/2} w(\gamma_{2,1}, \gamma_{2,2}) (\gamma_{2,1}\gamma_{2,2})^{m/2}&\text{if $m$ is even}\\
(\gamma_{2,1}\gamma_{2,2})^{-(m-1)/2} \gamma_{2,2}^{-1}w(\gamma_{2,2}, \gamma_{2,1}) \gamma_{2,2} (\gamma_{2,1}\gamma_{2,2})^{(m-1)/2}&\text{if $m$ is odd.}
\end{cases}
\end{multline*}
So for all $j\in\Z$, $\psi(\phi(\sigma^j)(w)w^{-1})= \psi([\sigma^j,w])=0$. Since the same is true for products and conjugates in $\F[2](\gamma_{2,1}, \gamma_{2,2})$, we see that $N\subseteq \ker{\psi}$, and thus $L\subseteq \ker{\psi}$. Now let us show that $\ker{\psi}\subseteq L$. To see this, we first apply the Reidemeister-Schreier rewriting process in order to obtain a basis of $\ker{\psi}$ (which is a free group since it is a subgroup of $\F[2](\gamma_{2,1}, \gamma_{2,2})$). Taking $\brak{\gamma_{2,1}, \gamma_{2,2}}$ as the basis of $\F[2](\gamma_{2,1}, \gamma_{2,2})$ and $\brak{\gamma_{2,1}^i}_{i\in\Z}$ as a Schreier transversal, the process yields $\brak{\gamma_{2,1}^i \gamma_{2,2} \gamma_{2,1}^{-(i+1)}}_{i\in\Z}$ as a basis. But for all $i\in\Z$,
\begin{equation*}
\gamma_{2,1}^i \gamma_{2,2} \gamma_{2,1}^{-(i+1)}=\gamma_{2,1}^i \gamma_{2,2} \gamma_{2,1}^{-1}\gamma_{2,1}^{-i}= \gamma_{2,1}^i \phi(\sigma)(\gamma_{2,1})\gamma_{2,1}^{-1} \gamma_{2,1}^{-i},
\end{equation*}
which belongs to $L$ by definition. This proves that $\ker{\psi}=L$, and that $L$ is a free group of infinite rank with the given basis as required.
\end{proof}

\begin{rem}
Since $\Gamma_2(\F[2](\gamma_{2,1}, \gamma_{2,2}))$ is the normal closure of the commutator $[\gamma_{2,1}, \gamma_{2,2}]$ in $\F[2](\gamma_{2,1}, \gamma_{2,2})$, and
\begin{align*}
[\gamma_{2,1}, \gamma_{2,2}]&=(\gamma_{2,1}\cdot \gamma_{2,2}\gamma_{2,1}^{-1}\cdot \gamma_{2,1}^{-1})(\gamma_{2,2}\gamma_{2,1}^{-1})^{-1}\\
&=
(\gamma_{2,1}\cdot \phi(\sigma)(\gamma_{2,1})\gamma_{2,1}^{-1}\cdot \gamma_{2,1}^{-1})(\phi(\sigma)(\gamma_{2,1})\gamma_{2,1}^{-1})^{-1},
\end{align*}
it follows that $\Gamma_2(\F[2](\gamma_{2,1}, \gamma_{2,2}))$ is contained in $N$, and so $L=N=\ker{\psi}$.
\end{rem}

\begin{rems}\label{rem:b2s2}
In fact, the group $B_2(\St\setminus\brak{x_1,x_2})$ is of particular interest since it may be interpreted in several different ways.
\begin{enumerate}[(a)]
\item As well as being isomorphic to $B_{2,1}(\dt)$, it is also isomorphic to the $2$-string braid group of the annulus.
\item One may reduce the presentation given by \req{b2s2} to the following:
\begin{equation}
B_2(\St\setminus\brak{x_1,x_2})=\setang{\sigma, \gamma_{2,2}}{(\sigma \gamma_{2,2})^2= (\gamma_{2,2} \sigma)^2},
\end{equation}
which is nothing other than the Artin group of type $B_2$~\cite{Cr,T}.
\item\label{it:inter3} The above presentation shows that $B_2(\St\setminus\brak{x_1,x_2})$ is a one-relator group. Interpreting it as the $2$-string braid group of the annulus, it follows from~\cite{PR} that it has infinite cyclic centre generated by the full twist of $B_3(\dt)$, which written in terms of our generators, is of the form $(\sigma \gamma_{2,2})^2$. Further, the relation may be written as $[\sigma, (\sigma \gamma_{2,2})^2]=1$. In particular, $B_2(\St\setminus\brak{x_1,x_2})$ is a one-relator group with non-trivial centre.
\item\label{it:inter4} Setting $D=\sigma \gamma_{2,2}$, from above, we obtain the presentation 
\begin{equation}\label{eq:presbs}
\setang{\sigma, D}{[\sigma, D^2]=1}.
\end{equation}
So $B_2(\St\setminus\brak{x_1,x_2})$ is isomorphic to the Baumslag-Solitar group $\operatorname{BS}(2,2)$~\cite{BS}.
\item Following~\cite{FG}, using the presentation~\reqref{presbs}, consider the homomorphism of $B_2(\St\setminus\brak{x_1,x_2})$ onto $\Z[D]=\ang{D}$ given by taking the exponent sum of $D$. It follows from the Reidemeister-Schreier rewriting process that the kernel is a free group $\F[2](\sigma, D\sigma D^{-1})$ of rank two, and thus that
\begin{equation*}
B_2(\St\setminus\brak{x_1,x_2})\cong \F[2](\sigma, D\sigma D^{-1})\rtimes \Z[D],
\end{equation*}
where the action is given by $D\cdot (\sigma)=D\sigma D^{-1}$, and $D\cdot (D\sigma D^{-1})= \sigma$. In other words, the action exchanges the two basis elements of the kernel (and not just up to conjugation as in \req{sigw}), and so is an involution. From this, it follows that $B_2(\St\setminus\brak{x_1,x_2})$ is an HNN-extension of the free group $\F[2](\sigma, D\sigma D^{-1})$ with stable letter $D$.
\item Still following~\cite{FG} and using the presentation~\reqref{presbs}, consider the homomorphism of $B_2(\St\setminus\brak{x_1,x_2})$ onto $\Z[\sigma]=\ang{\sigma}$ given by taking the exponent sum of $\sigma$. Applying the Reidemeister-Schreier rewriting process, one sees that that the kernel is generated by an infinite number of generators $x_i=\sigma_i D\sigma_i^{-1}$, $i\in \Z$, subject to the relations $x_i^2=x_0^2=D^2$ for all $i\in \Z$.
\end{enumerate}
\end{rems}

Applying~\cite{KMc,McC} to \rerems{b2s2}(\ref{it:inter3}) above, we see immediately that:
\begin{varthm}[\repr{b2b2resid}]
$B_2(\St\setminus\brak{x_1,x_2})$ is residually nilpotent and residually a finite $2$-group.\qed
\end{varthm}

Using the algorithm given in~\cite{CFL}, one may determine the quotient groups of the lower central series of $B_2(\St\setminus\brak{x_1,x_2})$. But these quotients may also be obtained explicitly using the results of~\cite{Ga,La}:
\begin{varthm}[\reth{lcsb2}]
For all $i\geq 2$, $\Gamma_i(B_2(\St\setminus\brak{x_1,x_2})) \cong \Gamma_i(\Z_2\ast \Z)$, and:
\begin{align*}
\Gamma_i(B_2(\St\setminus\!\brak{x_1,x_2}))/ \Gamma_{i+1}(B_2(\St\setminus\!\brak{x_1,x_2}))&\cong \Gamma_i(\Z_2\ast \Z)/ \Gamma_{i+1}(\Z_2\ast \Z)\\
& \cong \underbrace{\Z_2\oplus \cdots \oplus \Z_2}_{\text{$R_i$ times}},
\end{align*}
where
\begin{equation*}
R_i=\sum_{j=0}^{i-2}\; \left( \sum_{\substack{k\mid i-j\\ k>1}}\; \mu\left( \frac{i-j}{k}\right) \frac{k\alpha_k}{i-j} \right),
\end{equation*}
$\mu$ is the M\"obius function, and
\begin{equation*}
\alpha_k=\frac{1}{k} \left( \tr{\begin{pmatrix}
0 & -1\\
-1 & 1
\end{pmatrix}^k}-1\right).
\end{equation*}
\end{varthm}

\begin{rems}\mbox{}
\begin{enumerate}[(a)]
\item One may check that $\displaystyle R_{i+1}=R_i+\sum_{\substack{k\mid i+1\\ k>1}}\; \mu\left( \frac{i+1}{k}\right) \frac{k\alpha_k}{i+1}$, and that $\displaystyle \tr{\begin{pmatrix}
0 & -1\\
-1 & 1
\end{pmatrix}^k}=\left( \frac{1-\sqrt{5}}{2}\right)^k+ \left( \frac{1+\sqrt{5}}{2}\right)^k$.
\item By induction, one obtains $\begin{pmatrix}
0 & -1\\
-1 & 1
\end{pmatrix}^k=\begin{pmatrix}
f_{k-1} & -f_k\\
-f_k & f_{k+1}
\end{pmatrix}$, where $(f_k)_{k\geq 0}$ is the classical Fibonacci sequence defined by $f_0=0$, $f_1=1$, and $f_{k+2}=f_{k+1}+f_k$ for all $k\geq 0$. 
\item A simple calculation shows that $R_2=1$, $R_3=2$, $R_4=3$, $R_5=5$ and $R_6=7$.
\end{enumerate}
\end{rems}

The following lemma and corollary will be used in the proof of \reth{lcsb2}.

\begin{lem}\label{lem:ordern}
Let $G$ be a finitely-generated group. 
\begin{enumerate}[(a)]
\item\label{it:orderna} Suppose that there exists $i\geq 2$ such that $\Gamma_i(G)/\Gamma_{i+1}(G)$ is a torsion group. Then for all $j\geq i$, $\Gamma_j(G)/\Gamma_{j+1}(G)$ is a torsion group.
\item\label{it:ordernb} Suppose that there exists $i\geq 2$ and $n\in \N$ such that $x^n=1$ for all $x\in \Gamma_i(G)/\Gamma_{i+1}(G)$. Then for all $j\geq i$, $y^n=1$ for all $y\in \Gamma_j(G)/\Gamma_{j+1}(G)$.
\end{enumerate}
\end{lem}

\begin{proof}[Proof of \relem{ordern}]
Let $X$ be a finite set of generators of $G$. From~\cite{MKS}, we recall that for all $i\geq 2$, $\Gamma_i(G)/\Gamma_{i+1}(G)$ is a finitely-generated Abelian group, generated by the cosets of the simple $i$-fold commutators of elements of $X$. We prove part~(\ref{it:orderna}) by induction on $j$: suppose that $\Gamma_j(G)/\Gamma_{j+1}(G)$ is a torsion group for some $j\geq 2$. Now let $y\in \Gamma_{j+1}(G)/\Gamma_{j+2}(G)$. Then there exist simple $j$-fold commutators $x_1,\ldots,x_k\in \Gamma_j(G)$, $z_1,\ldots, z_k\in G$ and $\epsilon_1, \ldots, \epsilon_k\in \brak{\pm 1}$ such that $y$ is equal to the $\Gamma_{j+2}(G)$-coset of $[x_1,z_1]^{\epsilon_1}\cdots [x_k,z_k]^{\epsilon_k}$. By hypothesis, there exist $m_1, \ldots, m_k\in\N$ such that $x_i^{m_i}\in \Gamma_{j+1}(G)$ for $i=1,\ldots,k$. Set $m=\ppcm{m_1,\ldots}{m_k}$. Then modulo  $\Gamma_{j+2}(G)$, 
\begin{equation*}
y^m \equiv ([x_1,z_1]^{\epsilon_1}\cdots [x_k,z_k]^{\epsilon_k})^m \equiv
[x_1^m,z_1]^{\epsilon_1}\cdots [x_k^m,z_k]^{\epsilon_k}\equiv 1,
\end{equation*}
since each of the commutators $[x_i^m,z_i]^{\epsilon_i}$ belongs to $\Gamma_{j+2}(G)$. This proves part~(\ref{it:orderna}). Part~(\ref{it:ordernb}) follows similarly, taking $m_1= \cdots =m_k=n$ in the above proof.
\end{proof}

\begin{cor}\label{cor:order2}
The lower central series quotients of $\Z_2\ast \Z$ are isomorphic to the direct sum of a finite number of copies of $\Z_2$.
\end{cor}

\begin{proof}[Proof of \reco{order2}]
Let $x,y$ generate $\Z_2$ and $\Z$ respectively. Then $\Gamma_2(\Z_2\ast \Z)/ \Gamma_3(\Z_2\ast \Z)$ is a cyclic group generated by the coset of $[x,y]$. But modulo $\Gamma_3(\Z_2\ast \Z)$, $[x,y]^2\equiv [x^2,y]\equiv 1$, and the result follows from \relem{ordern} and using the fact that the lower central series quotients of $\Z_2\ast \Z$ are finitely-generated Abelian groups.
\end{proof}

\begin{proof}[Proof of \reth{lcsb2}]
Consider the presentation~\reqref{presbs} of the group $B_2(\St\setminus\brak{x_1,x_2})$. Let $\Z_2\ast\Z=\setang{\overline{D}, \overline{\sigma}}{\overline{D}^2=1}$. Since the centre of $B_2(\St\setminus\brak{x_1,x_2})$ is generated by $D^2$, we obtain the following central extension:
\begin{equation*}
1\to \ang{D^2} \to B_2(\St\setminus\brak{x_1,x_2}) \stackrel{\psi}{\to} \Z_2\ast\Z \to 1,
\end{equation*}
where $\psi(D)=\overline{D}$ and $\psi(\sigma)=\overline{\sigma}$.
Since $\psi$ is surjective, for $i\geq 2$, it induces a surjection $\map{\psi_i}{\Gamma_i(B_2(\St\setminus\brak{x_1,x_2}))}[\Gamma_i(\Z_2\ast \Z)]$. Using the fact that $\gpab[(B_2(\St\setminus\brak{x_1,x_2}))]= \ang{D,\sigma}\cong\Z^2$, this gives rise to the following commutative diagram of short exact sequences:
\begin{equation*}
\xymatrix@C=0.7cm{%
1 \ar[r] & \Gamma_2(B_2(\St\setminus\brak{x_1,x_2})) \ar[r] \ar[d]_{\psi_2} & B_2(\St\setminus\brak{x_1,x_2}) \ar^(0.63){\text{Ab}}[r] \ar[d]_{\psi} & \Z\oplus\Z \ar[r] \ar[d] & 1\\
1 \ar[r] & \Gamma_2(\Z_2\ast \Z) \ar[r]  & \Z_2\ast\Z \ar^{\text{Ab}}[r]  & \Z_2\oplus\Z \ar[r]  & 1,}
\end{equation*}
where $\text{Ab}$ denotes Abelianisation. Now $\psi_2$ is injective, since if $x\in\ker{\psi_2}$ then $x\in\ker{\psi}$, so there exists $k\in\Z$ such that $x=D^{2k}$. But since $x\in \Gamma_2(B_2(\St\setminus\brak{x_1,x_2}))$, its Abelianisation is trivial, so $k=0$. Hence $\psi_2$ is an isomorphism. But for $i\geq 2$, since $\psi_{i+1}$ is the restriction of $\psi_2$ to $\Gamma_{i+1}(B_2(\St\setminus\brak{x_1,x_2}))$ onto $\Gamma_{i+1}(\Z_2\ast \Z)$, it follows that $\psi_i$ is an isomorphism for all $i\geq 2$, and that 
\begin{equation*}
\Gamma_i(B_2(\St\setminus\brak{x_1,x_2}))/ \Gamma_{i+1}(B_2(\St\setminus\brak{x_1,x_2}))\cong \Gamma_i(\Z_2\ast \Z)/ \Gamma_{i+1}(\Z_2\ast \Z).
\end{equation*}
This proves the first part of the theorem. 

We now calculate the successive lower central series quotients $\Gamma_i(\Z_2\ast \Z)/\Gamma_{i+1}(\Z_2\ast \Z)$. This may be done by applying the results of~\cite{Ga,La}; we follow those of~\cite{Ga}.
From \reco{order2}, for each $i\geq 2$, $\Gamma_i(\Z_2\ast \Z)/\Gamma_{i+1}(\Z_2\ast \Z)$ is the direct sum of a finite number, denoted by $R_i$ in~\cite{Ga}, of copies of $\Z_2$. 

To determine $R_i$, one may first check that in Theorem~2.2
of~\cite{Ga}, $U_{\infty}(x)=0$ and $R_k^{\infty}=0$ for all $k\geq 2$
($R_k^{\infty}$ represents the rank of the free abelian  factor of
$\Gamma_k(\Z_2\ast \Z)/\Gamma_{k+1}(\Z_2\ast \Z)$). Secondly,
referring to the notation of Section~3 of that paper, we see that
$y=x$, $z=\frac{x}{1-x}$, $U(x)=\frac{x^2}{1-x}$, and
\begin{align*}
\frac{d}{dx}\left( \ln(1-U(x)) \right)&=\frac{x(x-2)}{(x-1)(x^2+x-1)}\\
&=\frac{-1}{x-1}+ \frac{1}{x-\lambda_+}+ \frac{1}{x-\lambda_-},
\end{align*}
where $\lambda_{\pm}=\frac{-1\pm\sqrt{5}}{2}$ are the roots of $x^2+x-1$. So from equation~(3.22) of~\cite{Ga}, we observe that for $k\geq 2$,
\begin{equation*}
\alpha_k=\frac{1}{k} \left( \tr{M^k}-1\right),
\end{equation*}
where $M=\begin{pmatrix}
0 & -1\\
-1 & 1
\end{pmatrix}$, and $\tr{M^k}=(-1)^k\left( \lambda_+^k+\lambda_-^k\right)$. The second part of the theorem then follows from Theorem~3.4 of~\cite{Ga}.
\end{proof}

We may thus describe the derived series of
$B_2(\St\setminus\brak{x_1,x_2})$ in terms of that of the free group
of rank~$2$:
\begin{cor}\label{cor:dsb22}
For all $i\in\N$,
\begin{equation*}
(B_2(\St\setminus\brak{x_1,x_2}))^{(i)}\cong\pi((\Z\ast \Z)^{(i)}),
\end{equation*}
where $\map{\pi}{\Z\ast \Z}[\Z_2 \ast \Z]$ is the homomorphism obtained by taking the first factor modulo~$2$.
\end{cor}

\begin{proof}
Let $G_1,G_2$ be two groups. If $\map{\pi}{G_1}[G_2]$ is a surjective homomorphism, then the restriction $\map{\pi\vert_{(G_1)^{(1)}}}{(G_1)^{(1)}}[(G_2)^{(1)}]$, and by induction on $i$, so is the restriction $\map{\pi\vert_{(G_1)^{(i)}}}{(G_1)^{(n)}}[(G_2)^{(i)}]$. Taking $G_1=\Z\ast\Z$ and $G_2=\Z_2\ast \Z$, it follows that
\begin{equation*}
(\Z_2\ast\Z)^{(i)} =\pi((\Z\ast\Z)^{(i)}).
\end{equation*}
But 
\begin{equation*}
(B_2(\St\setminus\brak{x_1,x_2}))^{(i)}\cong(\Z_2\ast\Z)^{(i)}
\end{equation*}
by \reth{lcsb2} which proves the corollary.
\end{proof}

We now determine explicitly $\Gamma_3(B_2(\St\setminus\brak{x_1,x_2}))$. 
\begin{prop}\label{prop:gam3b2}
Let $\map{\rho_2}{\Gamma_2(B_2(\St\setminus\brak{x_1,x_2}))}[\Z_2]$ be the homomorphism defined by $\rho_2(z_n)=1$ for all $n\in \Z$, where $\brak{z_n}_{n\in\Z}$ is the basis given by \reco{gam2b2s2}. Then $\Gamma_3(B_2(\St\setminus\brak{x_1,x_2}))= \ker{\rho_2}$. In particular, $\Gamma_3(B_2(\St\setminus\brak{x_1,x_2}))$ is a free group of infinite rank with a basis given by $\brak{z_n z_0^{-1}}_{n\in\Z\setminus\brak{0}}\bigcup \brak{z_m^2}_{m\in\Z}$, and \begin{equation*}
\Gamma_2(B_2(\St\setminus\brak{x_1,x_2}))/ \Gamma_3(B_2(\St\setminus\brak{x_1,x_2}))\cong \Z_2.
\end{equation*}
\end{prop}

\begin{rem}
Since $R_2=1$, this agrees with the result of \reth{lcsb2} in the case $i=2$.  
\end{rem}

\begin{proof}
We start by calculating the action under conjugation of the generators $\gamma_{2,1}, \gamma_{2,2}$ and $\sigma$ of $B_2(\St\setminus\brak{x_1,x_2})$ on the generators $z_n$ of $\Gamma_2(B_2(\St\setminus\brak{x_1,x_2}))$. Clearly $\gamma_{2,1} z_n \gamma_{2,1}^{-1}= z_{n+1}$ and $\gamma_{2,2}z_n \gamma_{2,2}^{-1}= z_0 z_{n+1} z_0^{-1}$. Further, it follows from \req{sigw} that 
\begin{equation*}
\sigma z_n \sigma^{-1}= \gamma_{2,2}^{n-1} \gamma_{2,1} \gamma_{2,2}^{-n},
\end{equation*}
which rewriting in terms of the $z_i$ yields:
\begin{equation*}
\sigma z_n \sigma^{-1}= 
\begin{cases}
z_0 z_1\cdots z_{n-2} z_{n-1}^{-1} z_{n-2}^{-1}\cdots z_1^{-1}z_0^{-1} & \text{if $n> 0$}\\
z_{-1}^{-1}\cdots z_{-\lvert n\rvert}^{-1} z_{-(\lvert n\rvert+1)}^{-1} z_{-\lvert n\rvert}\cdots z_{-1}& \text{if $n\leq 0$.}
\end{cases}
\end{equation*}
Let us apply the Reidemeister-Schreier rewriting process to the basis
$\brak{z_n}_{n\in\Z}$ of $\Gamma_2(B_2(\St\setminus\brak{x_1,x_2}))$,
taking the Schreier transversal $\brak{1,z_0}$ for $\rho_2$. This
yields a basis $\brak{z_n z_0^{-1}}_{n\in\Z\setminus\brak{0}}\bigcup
\brak{z_0 z_m}_{m\in\Z}$ of $\ker{\rho_2}$, or equivalently a basis
$\brak{z_n z_0^{-1}}_{n\in\Z\setminus\brak{0}}\bigcup
\brak{z_m^2}_{m\in\Z}$. Since 
\begin{equation*}
z_n z_0^{-1}=
\begin{cases}
(z_n z_{n-1}^{-1})(z_{n-1} z_{n-2}^{-1})\cdots (z_1 z_0^{-1})& \text{for all $n>0$}\\
(z_{n+1} z_n^{-1})^{-1} (z_{n+2} z_{n+1}^{-1})^{-1} \cdots (z_0 z_{-1}^{-1})^{-1}&\text{for all $n<0$,}
\end{cases}
\end{equation*}
and $z_{i+1} z_i^{-1}=[a,z_i]\in \Gamma_3(B_2(\St\setminus\brak{x_1,x_2}))$ for all $i\in \Z$, we see that $z_n z_0^{-1} \in \Gamma_3(B_2(\St\setminus\brak{x_1,x_2}))$ for all $n\neq 0$. Finally, if $m\in\Z$ then $z_mz_0=(z_mz_0^{-1})z_0^2$. But $[\sigma,z_1]= z_0^{-1}z_1^{-1}$, so $z_0^2= (z_1 z_0^{-1})^{-1} [\sigma,z_0]^{-1} \in \Gamma_3(B_2(\St\setminus\brak{x_1,x_2}))$. Thus $\ker{\rho_2} \subseteq\Gamma_3(B_2(\St\setminus\brak{x_1,x_2}))$.

To prove the converse, observe first that
$\Gamma_3(B_2(\St\setminus\brak{x_1,x_2}))$ is the normal closure in
$B_2(\St\setminus\brak{x_1,x_2})$ of the commutators $[\gamma_{2,1},
z_n]$, $[\gamma_{2,2}, z_n]$ and $[\sigma, z_n]$, where $n\in\Z$. It
follows easily from the above expressions that these elements belong
to $\ker{\rho_2}$. Further, conjugation by each of $\gamma_{2,1}$,
$\gamma_{2,2}$ and $\sigma$ induce automorphisms of
$\Gamma_2(B_2(\St\setminus\brak{x_1,x_2}))$, each of which leaves
$\ker{\rho_2}$ invariant, and so induces an automorphism of $\Z_2$,
which is in fact the identity in all three cases. Hence for all
$n\in\Z$, all conjugates of $[\gamma_{2,1}, z_n]$, $[\gamma_{2,2},
z_n]$ and $[\sigma, z_n]$ by elements of
$B_2(\St\setminus\brak{x_1,x_2})$ belong to $\ker{\rho_2}$, and so
$\Gamma_3(B_2(\St\setminus\brak{x_1,x_2}))\subseteq \ker{\rho_2}$. We
conclude that $\Gamma_3(B_2(\St\setminus\brak{x_1,x_2}))\subseteq
\ker{\rho_2}$, and $\Gamma_2(B_2(\St\setminus\brak{x_1,x_2}))/
\Gamma_3(B_2(\St\setminus\brak{x_1,x_2}))\cong \Z_2$.
\end{proof}

\section{The commutator subgroup of $B_m(\St\setminus\brak{x_1,x_2})$, $m\geq 3$}\label{sec:affineatil}

As we already observed in \rerems{annulus}, $B_m(\St\setminus\brak{x_1,x_2})$ may be identified with the $m$-string braid group of the annulus. The case $m=2$ having already been studied in \resec{lcgsb22}, let us now suppose that $m\geq 3$. In this case, we know from \reth{lcdsbmsn} that the lower central series is constant from the commutator subgroup onwards. The following presentation of $B_m(\St\setminus\brak{x_1,x_2})$ was obtained by Kent and Peifer:
\begin{prop}[\cite{KP}]\label{prop:kent}
If $m\geq 3$ then $B_m(\St\setminus\brak{x_1,x_2})$ admits a presentation of the following form:
\begin{enumerate}
\item[\underline{\textbf{generators:}}] $\sigma_0, \sigma_1,\ldots, \sigma_{m-1}$ and $\tau$.
\item[\underline{\textbf{relations:}}] 
\begin{gather}
\text{$\si{i}\si{j}=\si{j}\si{i}$ if $\lvert i-j\rvert\neq 1,m-1$ and $0\leq i,j\leq m-1$}\label{eq:annu1}\\
\text{$\si{i}\si{i+1}\si{i}=\si{i+1}\si{i}\si{i+1}$ for $0\leq i\leq m-1$, and}\label{eq:annu2}\\
\text{$\tau^{-1}\sigma_i\tau=\sigma_{i+1}$ for $0\leq i\leq m-1$.}\label{eq:annu3}
\end{gather}
The indices should be taken modulo $m$.
\end{enumerate}
\end{prop}
The $m$ points should be thought of as being arranged around the centre of the annulus. The generator $\sigma_0$ corresponds to a positive half-twist between the $m\up{th}$ and $1\up{st}$ point, while $\tau$ is represented geometrically by a rigid rotation of the annulus about the centre by an angle $2\pi/n$. It follows from this presentation that:
\begin{cor}[\cite{KP}]\label{cor:kent}
If $m\geq 3$ then $B_m(\St\setminus\brak{x_1,x_2})$ is isomorphic to the semi-direct product of the affine Artin group $\widetilde{A}_{m-1}$ (generated by $\sigma_0, \sigma_1,\ldots, \sigma_{m-1}$, and subject to relations~\reqref{annu1} and~\reqref{annu2}) by the infinite cyclic group generated by $\tau$, the action being that of conjugation given by relation~\reqref{annu3}.
\end{cor}

Then we have the following result:
\begin{prop}\label{prop:gamma2annu}\mbox{}
\begin{enumerate}[(a)]
\item If $m\geq 3$ then $\Gamma_2(B_m(\St\setminus\brak{x_1,x_2}))$ is generated by the elements $p_k=\sigma_1^k\sigma_2 \sigma_1^{-(k+1)}$, $r_k=\sigma_1^k\sigma_0 \sigma_1^{-(k+1)}$, for all $k\in\Z$, and $q_i=\sigma_i\sigma_1^{-1}$ for $3\leq i\leq m-1$.
\item\label{it:m3gam2} If $m=3$, then $\Gamma_2(B_3(\St\setminus\brak{x_1,x_2}))$ is defined by the following relations:
\begin{gather}
p_{k+1} p_{k+2}^{-1} p_k^{-1}=1\label{eq:m3pp}\\
r_{k+1} r_{k+2}^{-1} r_k^{-1}=1\label{eq:m3rr}\\
r_k p_{k+1} r_{k+2} p_{k+2}^{-1} r_{k+1}^{-1} p_k^{-1}=1,\label{eq:m3rp}
\end{gather}
where $k\in\Z$.
\item\label{it:m4gam2} If $m\geq 4$ then $\Gamma_2(B_m(\St\setminus\brak{x_1,x_2}))$ is defined by the following relations:
\begin{gather}
p_{k+1} p_{k+2}^{-1} p_k^{-1}=1\label{eq:m4pp}\\
r_{k+1} r_{k+2}^{-1} r_k^{-1}=1\label{eq:m4rr}\\
p_k q_3 p_{k+2} q_3^{-1} p_{k+1}^{-1} q_3^{-1}=1\label{eq:m4pq3}\\
\text{$p_k q_i p_{k+1}^{-1} q_i^{-1}=1$ for all $4\leq i\leq m-1$} \label{eq:m4pqi}\\
\text{$q_iq_jq_i^{-1}q_j^{-1}=1$ for all $3\leq i<j-1\leq m-2$} \label{eq:m4qiqj}\\
\text{$q_iq_{i+1}q_i=q_{i+1}q_iq_{i+1}$ for $3\leq i\leq m-2$} \label{eq:m4qi1}\\
r_k p_{k+1} r_{k+1}^{-1} p_k^{-1}=1 \label{eq:m4rp}\\
\text{$r_k q_i r_{k+1}^{-1} q_i^{-1}=1$ for all $3\leq i\leq m-2$} \label{eq:m4rq}\\
r_k q_{m-1} r_{k+2} q_{m-1}^{-1} r_{k+1}^{-1} q_{m-1}^{-1}=1,\label{eq:m4rm}
\end{gather}
where $k\in\Z$.
\end{enumerate}
\end{prop}

We may thus deduce the first derived series quotient of the group $\Gamma_2(B_m(\St\setminus\brak{x_1,x_2}))$:
\begin{varthm}[\reco{dsannu}]
Let $m\geq 3$. Then
\begin{equation*}
\left(B_m \left(\St\setminus\brak{x_1,x_2}\right)\right)^{(1)}\left/ \left(B_m\left(\St\setminus\brak{x_1,x_2}\right)\right)^{(2)}\right. \cong
\begin{cases}
\Z^4 & \text{if $m=3$}\\
\Z^2 & \text{if $m=4$}\\
\Z & \text{if $m\geq 5$.}
\end{cases}
\end{equation*}
\end{varthm}

\begin{proof}[Proof of \repr{gamma2annu}]
We start by applying \repr{gammasemi} to the result of \reco{kent}, namely that 
\begin{equation*}
B_m(\St\setminus\brak{x_1,x_2}) \cong \widetilde{A}_{m-1} \rtimes \ang{\tau}. 
\end{equation*}
If $w=\sigma_{i_1}^{\epsilon_{i_1}} \cdots \sigma_{i_k}^{\epsilon_{i_k}}\in \widetilde{A}_{m-1}$, it follows from the action, given by \req{annu3}, that for all $l\in\Z$,
\begin{equation*}
\tau^{-l}w\tau^l\cdot w^{-1}= \sigma_{i_1+l}^{\epsilon_{i_1}} \cdots \sigma_{i_k+l}^{\epsilon_{i_k}} \cdot \sigma_{i_k}^{-\epsilon_{i_k}}\cdots \sigma_{i_1}^{-\epsilon_{i_1}},
\end{equation*}
where the indices should be taken modulo~$m$. Hence $\tau^{-l}w\tau^l\cdot w^{-1}\in \Gamma_2(\widetilde{A}_{m-1})$, and it follows from \repr{gammasemi} that 
\begin{equation*}
\Gamma_2(B_m(\St\setminus\brak{x_1,x_2}))\cong \Gamma_2(\widetilde{A}_{m-1}).
\end{equation*}

A presentation of $\Gamma_2(\widetilde{A}_{m-1})$ may be obtained by observing that $\gpab[(\widetilde{A}_{m-1})] \cong \Z$, and by applying the Reidemeister-Schreier rewriting process to the generating set $\brak{\sigma_0,\sigma_1,\ldots, \sigma_{m-1}}$ of $\widetilde{A}_{m-1}$ and the Schreier transversal $\brak{\sigma_1^k}_{k\in\Z}$. The generators and relations not containing $\sigma_0$ define a group isomorphic to $B_m(\dt)$, and using~\cite{GL}, we obtain all of the generators and relations of \repr{gamma2annu} not containing $r_k$. The generator $\sigma_0$ of $\widetilde{A}_{m-1}$ gives rise to generators $r_k=\sigma_1^k \sigma_0 \sigma_1^{-(k+1)}$ of $\Gamma_2(\widetilde{A}_{m-1})$, where $k\in\Z$. The relation~\reqref{annu2} with $j=0$ and $i=1$ yields relations of the form $r_{k+1} r_{k+2}^{-1} r_k^{-1}=1$, $k\in\Z$, in $\Gamma_2(\widetilde{A}_{m-1})$. If $m=3$ then we obtain relations~\reqref{m3rp} in $\Gamma_2(\widetilde{A}_2)$ from relation \reqref{annu2} with $j=0$ and $i=2$, and so we deduce the presentation given in part~(\ref{it:m3gam2}). If $m\geq 4$, taking $j=0$ in relations~\reqref{annu1} with $i=2$ (resp.\ $3\leq i\leq m-2$) yields relations~\reqref{m4rp} (resp.~\reqref{m4rq}) in $\Gamma_2(\widetilde{A}_{m-1})$. Finally we obtain relations~\reqref{m4rm} in $\Gamma_2(\widetilde{A}_{m-1})$ by taking $j=0$ and $i=m-1$ in relation~\reqref{annu2}, and this gives the presentation of part~(\ref{it:m4gam2}).
\end{proof}

\begin{proof}[Proof of \reco{dsannu}]
It suffices to Abelianise the presentations of \repr{gamma2annu}, in
other words, we add the commutation relations of all of the generators
to the given presentations. First let $m=3$. Equation~\reqref{m3rp}
becomes trivial using equations~\reqref{m3pp} and~\reqref{m3rr}.
Further, it follows from equations~\reqref{m3pp} (resp.~\reqref{m3rr})
that all of the $p_k$ (resp.\ $r_k$) may be expressed uniquely in
terms of $p_0$ and $p_1$ (resp.\ $r_0$ and $r_1$), and hence
$\gpab[(B_3(\St\setminus\brak{x_1,x_2}))]$ is a free Abelian group of rank~$4$ with basis $\brak{p_0,p_1,r_0,r_1}$

Let $m=4$. By \req{m4pp}, it follows from \req{m4pq3} that $q_3$
Abelianises to the trivial element, and then \req{m4rr} implies that
\req{m4rm} becomes trivial. By \req{m4rp}, $p_k=r_k$ for all $k\in\Z$.
As above, all of the $p_k$ (resp.\ $r_k$) may be expressed uniquely in
terms of $p_0$ and $p_1$ (resp.\ $r_0$ and $r_1$), and thus
$\gpab[(B_4(\St\setminus\brak{x_1,x_2}))]$ is a free Abelian group of rank~$4$ with basis $\brak{p_0,p_1}$.

Finally, if $m\geq 5$, by \req{m4qi1} we obtain additionally that all of the $q_i$ Abelianise to the trivial element. By \req{m4pqi} (resp.~\reqref{m4rq}), $p_k=p_{k+1}$ (resp.\ $r_k=r_{k+1}$). Thus $\gpab[(B_m(\St\setminus\brak{x_1,x_2}))]$ is a infinite cyclic group generated by $p_0$.
\end{proof}

\section{The series of
$B_m(\St\setminus\brak{x_1,x_2,x_3})$}\label{sec:bmn3}

The situation seems to be more difficult in the case of the braid group of the $3$-punctured sphere. As we remark below, if $m\geq 2$, $B_m(\St\setminus\brak{x_1,x_2,x_3})$ is isomorphic to the affine Artin group of type $\widetilde{C}_m$ for which little seems to be known~\cite{All,ChP}. We have not even been able to describe the commutator subgroup. We may however obtain some partial results, notably in \repr{b2n3} the fact that the successive lower central series quotients of $B_2(\St\setminus\brak{x_1,x_2,x_3})$ are direct sums of $\Z_2$, which generalises part of \reth{lcsb2}.

We begin by considering the case $m=2$.
\begin{prop}[\cite{BG}]\label{prop:presbg1}
The following constitutes a presentation of the group $B_2(\St\setminus\brak{x_1,x_2,x_3})$:
\begin{enumerate}
\item[\underline{\textbf{generators:}}] $\sigma$, $\rho_1$ and $\rho_2$.
\item[\underline{\textbf{relations:}}] 
\begin{gather}
(\sigma\rho_1)^2 = (\rho_1\sigma)^2\label{eq:sigrho1}\\
(\sigma\rho_2)^2 = (\rho_2\sigma)^2\label{eq:sigrho2}\\
\rho_1\rho_2=\rho_2\rho_1.\notag
\end{gather}
\end{enumerate}
\end{prop}
Geometrically, $B_2(\St\setminus\brak{x_1,x_2,x_3})$ may be considered as the $2$-string braid group of the twice-punctured disc, which in turn may be considered as a subgroup of $B_4(\dt)$ whose first and fourth strings are vertical. Then with the usual notation, $\rho_1=A_{1,2}$, $\rho_2=A_{3,4}$, and $\sigma$ is the positive half-twist of the second and third strings.

Let $G_1$ be the group generated by $\sigma$ and $\rho_1$ subject to
the relation~\reqref{sigrho1}, and let $G_2$ be the group
generated by $\sigma$ and $\rho_2$ subject to the
relation~\reqref{sigrho2}.  It follows from the above proposition and \rerem{b2s2}(\ref{it:inter4})
that $B_2(\St\setminus\brak{x_1,x_2,x_3})$ may be considered as the
amalgamated product $G_1 \ast_{\ang{\sigma}} G_2$ of two copies of the
Baumslag-Solitar group $\operatorname{BS}(2,2)$, subject to the
additional relation $[\rho_1,\rho_2]=1$. We wonder if it would be possible to obtain determine the commutator subgroup via this amalgamated product.

The following gives a generalisation to $B_2(\St\setminus\brak{x_1,x_2,x_3})$ of part of \reth{lcsb2}.
\begin{prop}\label{prop:b2n3}
For all $i\geq 2$, the lower central series quotient $\Gamma_i(B_2(\St\setminus\brak{x_1,x_2,x_3}))\left/ \Gamma_{i+1}(B_2(\St\setminus\brak{x_1,x_2,x_3}))\right.$ is isomorphic to the direct sum of a finite number of copies of $\Z_2$.
\end{prop}

\begin{proof}
As in the proof of \relem{ordern}, since $B_2(\St\setminus\brak{x_1,x_2,x_3})$ is finitely generated, it follows that the lower central quotient $\Gamma_i(B_2(\St\setminus\brak{x_1,x_2,x_3}))\left/ \Gamma_{i+1}(B_2(\St\setminus\brak{x_1,x_2,x_3}))\right.$ is a finitely-generated Abelian group. By part~(\ref{it:ordernb}) of \relem{ordern}, it suffices to prove the result in the case $i=2$, which we do using the presentation of \repr{presbg1}. We know that $\Gamma_2(B_2(\St\setminus\brak{x_1,x_2,x_3}))\left/ \Gamma_3(B_2(\St\setminus\brak{x_1,x_2,x_3}))\right.$ is generated by the $\Gamma_3$-cosets of the commutators of the form $[x,y]$, where $x,y\in \brak{\sigma,\rho_1, \rho_2}$, and thus of the commutators $[\sigma,\rho_i]$ for $i=1,2$. But $[\sigma,\rho_i]=[\rho_i^{-1} \sigma^{-1}]^{-1}$ by relations~\reqref{sigrho1} and~\reqref{sigrho2}. So modulo $\Gamma_3$, $[\sigma,\rho_i]$ is congruent to $[\sigma,\rho_i]^{-1}$, in other words, $[\sigma,\rho_i]^2$ is trivial modulo $\Gamma_3$, which proves the result. 
\end{proof}

As was pointed out in~\cite{All,BG},
$B_2(\St\setminus\brak{x_1,x_2,x_3})$ is isomorphic to the affine
Artin braid group $\widetilde{C}_2$. More generally, for $m\geq 2$,
$B_m(\St\setminus\brak{x_1,x_2,x_3})$ is isomorphic to
$\widetilde{C}_{m}$ and by~\cite{BG} has a presentation of the form:
\begin{enumerate}
\item[\underline{\textbf{generators:}}] $\rho_1,\rho_m$ and $\sigma_i$, $1\leq i\leq m-1$.
\item[\underline{\textbf{relations:}}] 
\begin{gather*}
\text{$\sigma_i\sigma_j=\sigma_j\sigma_i$ if $\lvert i-j\rvert\geq 2$ and $1\leq i,j\leq m-1$}\\
\text{$\sigma_i\sigma_{i+1}\sigma_i= \sigma_{i+1}\sigma_i\sigma_{i+1}$ for all $1\leq i\leq m-2$}\\
\rho_1 \comm \rho_m\\
\text{$\rho_1\comm \sigma_i$ for all $2\leq i\leq m-1$}\\
\text{$\rho_m\comm \sigma_i$ for all $1\leq i\leq m-2$}\\
(\sigma_1\rho_1)^2= (\rho_1\sigma_1)^2\\
(\sigma_{m-1}\rho_m)^2= (\rho_m\sigma_{m-1})^2.
\end{gather*}
\end{enumerate}

The following result yields information about the derived series quotients of $B_m(\St\setminus\brak{x_1,x_2,x_3})$.
\begin{prop}\label{prop:bmn3i}
Let $m\geq 2$. Then $B_m(\St\setminus\brak{x_1,x_2,x_3})$ is a semi-direct product of a group $K_0$ by $B_m(\dt)$. In particular, for all $i\geq 1$ $(B_m(\St\setminus\brak{x_1,x_2,x_3}))^{(i)}$ is a semi-direct of a group $K_i$ by $(B_m(\dt))^{(i)}$.
\end{prop}

\begin{proof}
Consider the homomorphism of $B_m(\St\setminus\brak{x_1,x_2,x_3})$ to $B_m(\dt)$ which sends $\rho_1$ and $\rho_m$ onto the trivial element. From the above presentation, it is clearly surjective, and it admits an obvious section. So if $K_0$ denotes the kernel then $B_m(\St\setminus\brak{x_1,x_2,x_3}) \cong K_0\rtimes B_m(\dt)$. The second part is obtained by induction on $i$, using \repr{gammasemi}. 
\end{proof}

\chapter{Presentations for $\Gamma_2(B_n(\St))$, $n\geq 4$}\label{chap:prescom}

In this chapter, we give various presentations of $\Gamma_2(B_n(\St))$, $n\geq 4$. In \resec{gen}, we begin by giving a general presentation obtained using the Reidemeister-Schreier rewriting process. In \resec{presg2bn4}, we consider the case $n=4$, and derive the presentation given in \reth{lcsbn}(\ref{it:lcs4}). In \resec{presg2bn5}, we restate the presentation given by \repr{fullpres} for the case $n=5$, and for $n\geq 6$, we refine the presentation to obtain \repr{g2b6}.

\section{A general presentation of $\Gamma_2(B_n(\St))$ for $n\geq 4$}\label{sec:gen}

\begin{prop}\label{prop:fullpres}
Let $n\geq 4$. The following constitutes a presentation of the group $\Gamma_2(B_n(\St))$:
\begin{enumerate}
\item[\underline{\textbf{generators:}}] 
\begin{gather*}
w=\si[2n-2]1\\
u_1=\si{2}\sii{1}, u_2=\si{1}\si{2}\sii[2]{1}, \ldots, u_{2n-2}=\si[2n-3]{1} \si{2}\sii[(2n-2)]{1}\\ v_1=\si{3}\sii{1}, \ldots, v_{n-3}=\si{n-1}\sii{1}.
\end{gather*}
\item[\underline{\textbf{relations:}}] 
\begin{gather}
\text{$v_iv_j=v_jv_i$ if $\lvert i-j\rvert\geq 2$ and $1\leq i,j\leq n-3$}\label{eq:eq1}\\
\text{$v_i v_{i+1}v_i=v_{i+1}v_iv_{i+1}$ for all $1\leq i\leq n-4$}\label{eq:eq2}\\
w \comm v_i\label{eq:eq3}\\
\text{$u_iv_j u_{i+1}^{-1}v_j^{-1}=1$ for $j\geq 2$ and $i=1,\ldots ,2n-3$}\label{eq:eq4}\\
\text{$u_{2n-2}v_j w u_1^{-1}w^{-1}v_j^{-1}=1$ for $2\leq j\leq n-3$}\label{eq:eq5}\\
\text{$u_iv_1 u_{i+2}v_1^{-1}u_{i+1}^{-1}v_1^{-1}=1$ for $i=1,\ldots, 2n-4$}\label{eq:eq6}\\
u_{2n-3}v_1 wu_1 w^{-1}v_1^{-1}u_{2n-2}^{-1}v_1^{-1}=1\label{eq:eq7}\\
u_{2n-2}v_1 wu_2 v_1^{-1}u_1^{-1}w^{-1}v_1^{-1}=1\label{eq:eq8}\\
\text{$u_{i+1}u_{i+2}^{-1}u_i^{-1}=1$ for all $i=1,\ldots, 2n-4$}\label{eq:eq9}\\
u_{2n-2}wu_1^{-1}w^{-1}u_{2n-3}^{-1}=1\label{eq:eq10}\\
wu_1u_2^{-1}w^{-1}u_{2n-2}^{-1}=1\label{eq:eq11}\\
u_2(v_1\cdots v_{n-4}v_{n-3}^2 v_{n-4}\cdots v_1) u_{2n-3}w=1\label{eq:eq12}\\
u_3(v_1\cdots v_{n-4}v_{n-3}^2 v_{n-4}\cdots v_1) u_{2n-2}w=1\label{eq:eq13}\\
\text{$u_i(v_1\cdots v_{n-4}v_{n-3}^2 v_{n-4}\cdots v_1) wu_{i-3}=1$ for $i=4,\ldots, 2n-2$}\label{eq:eq14}\\
u_1(v_1\cdots v_{n-4}v_{n-3}^2 v_{n-4}\cdots v_1) u_{2n-4}w=1.\label{eq:eq15}
\end{gather}
\end{enumerate}
\end{prop}

In what follows, we shall denote by equation~$(m_i)$ the equation~$(m)$ of the above system for the parameter value $i$.

\begin{proof}
Taking the standard presentation~\reqref{presnbns} of $B_n(\St)$, and the set $\brak{1, \sigma_1, \sigma_1^2,\ldots, \sigma_1^{2n-3}}$ as a Schreier tranversal, we apply the Reide\-meister-Schreier rewriting process to the following short exact sequence:
\begin{equation*}
\xymatrix{%
1\ar[r] & \Gamma_2(B_n(\St)) \ar[r] & B_n(\St) \ar[r] & \bnab{n} \ar[r] & 1.}
\end{equation*}
As generators of $\Gamma_2(B_n(\St))$, we obtain $w=\sigma_1^{2n-2}$, $\sigma_1^j \sigma_i \sigma_1^{-(j+1)}$ and $\sigma_1^{2n-3} \sigma_i$, where $2\leq i\leq n-1$ and $0\leq j\leq 2n-4$. We replace the latter by $\sigma_1^{2n-3} \sigma_i\cdot w^{-1}= \sigma_1^{2n-3} \sigma_i\sigma_1^{-(2n-2)}$. Now turning to the relations, if $i\geq 3$ then for $j=0,\ldots,2n-4$, the relator $\sigma_1 \sigma_i\sigma_1^{-1} \sigma_i^{-1}$ of $B_n(\St)$ gives rise to relators 
\begin{equation*}
\sigma_1^{j} \sigma_1 \sigma_i\sigma_1^{-1} \sigma_i^{-1} \sigma_1^{-j}= \sigma_1^{j+1} \sigma_i\sigma_1^{-(j+2)} \cdot \sigma_1^{j+1}\sigma_i^{-1} \sigma_1^{-j}
\end{equation*}
of $\Gamma_2(B_n(\St))$, so 
\begin{equation*}
\sigma_1^{j+1} \sigma_i\sigma_1^{-(j+2)}= \sigma_1^{j}\sigma_i \sigma_1^{-(j+1)}= \sigma_i \sigma_1^{-1}=v_{i-2}.
\end{equation*}
If $j=2n-3$ then we have a relator of the form
\begin{equation*}
\sigma_1^{2n-3} \sigma_1 \sigma_i\sigma_1^{-1} \sigma_i^{-1} \sigma_1^{-(2n-3)}= \sigma_1^{2n-2} \cdot \sigma_i\sigma_1^{-1} \cdot \sigma_1^{2n-2} \sigma_1^{-(2n-2)}\sigma_i^{-1} \sigma_1^{-(2n-3)},
\end{equation*}
and thus $v_{i-2}$ commutes with $w$, which gives relation~\reqref{eq3}. If $1\leq i,j\leq n-3$ and $\lvert i-j\rvert\geq 2$ then the relator $\sigma_{i+2} \sigma_{j+2}\sigma_{i+2}^{-1} \sigma_{j+2}^{-1}$ gives rise to the single relator $v_iv_jv_i^{-1}v_j^{-1}$, while if $1\leq i\leq n-4$, the relator $\sigma_{i+2} \sigma_{i+3}\sigma_{i+2} \sigma_{i+3}^{-1}\sigma_{i+2}^{-1} \sigma_{i+3}^{-1}$ yields the single relator
$v_i v_{i+1}v_i=v_{i+1}v_iv_{i+1}$, thus we obtain equations~\reqref{eq1} and~\reqref{eq2}.

Now for $i=1,\ldots, u_{2n-2}$, let $u_i=\sigma_1^{i-1} \sigma_2 \sigma_1^{-i}$. From the relator $\sigma_1^{j-1}\sigma_2 \sigma_1\sigma_2 \sigma_1^{-1}\sigma_2^{-1} \sigma_1^{-1} \sigma_1^{-(j-1)}$, we obtain the relators $u_ju_{j+2}u_{j+1}^{-1}$ if $j=1,\ldots, 2n-4$, $u_{2n-3}wu_1 w^{-1} u_{2n-2}^{-1}$ if $j=2n-3$, and $u_{2n-2}wu_2u_1^{-1}w^{-1}$ if $j=2n-2$, which gives respectively equations~\reqref{eq9},~\reqref{eq10} and~\reqref{eq11}.

If $2\leq i\leq n-3$ then the relator $\sigma_1^{j-1}\sigma_{i+2} \sigma_2 \sigma_{i+2}^{-1}\sigma_2^{-1} \sigma_1^{-(j-1)}$ yields  relators $v_i u_{j+1}v_i^{-1}u_j^{-1}$ if $j=1,\ldots, 2n-3$ and $v_iw u_1w^{-1} v_i^{-1}u_{2n-2}^{-1}$ if $j=2n-2$, and so we recover equations~\reqref{eq4} and~\reqref{eq5}.

From the relator $\sigma_1^{j-1}\sigma_3 \sigma_2 \sigma_3
\sigma_2^{-1} \sigma_3^{-1} \sigma_2^{-1} \sigma_1^{-(j-1)}$, we
obtain the relators $v_1 u_{j+1}v_1 u_{j+2}^{-1} v_1^{-1}u_j^{-1}$ if
$j=1,\ldots, 2n-4$, $v_1 u_{2n-2}v_1 w u_1^{-1} w^{-1}
v_1^{-1}u_{2n-3}^{-1}$ if $j=2n-3$, and  $v_1 w u_1v_1 u_2^{-1} w^{-1}
v_1^{-1}u_{2n-2}^{-1}$ if $j=2n-2$, which gives
equations~\reqref{eq6},~\reqref{eq7} and~\reqref{eq8}.

Finally,
\begin{multline*}
\sigma_1 \sigma_2\cdots \sigma_{n-2} \sigma_{n-1}^2 \sigma_{n-2} \cdots \sigma_2 \sigma_1= \sigma_1\sigma_2 \sigma_1^{-2}\cdot \sigma_1^2\sigma_3 \sigma_1^{-3} \cdots\\  \cdots\sigma_1^{n-3}\sigma_{n-2} \sigma_1^{-(n-2)}\cdot \sigma_1^{n-2}\sigma_{n-1} \sigma_1^{-(n-1)}\cdot \sigma_1^{n-1}\sigma_{n-1} \sigma_1^{-n}\cdot\\ \sigma_1^n \sigma_{n-2} \sigma_1^{-(n+1)}\cdot \sigma_1^{2n-4} \sigma_2 \sigma_1^{-(2n-3)}\cdot \sigma_1^{2n-2},
\end{multline*}
and conjugating by $\sigma_1^{j-1}$, we obtain relators 
\begin{equation*}
\begin{cases}
u_2(v_1\cdots v_{n-4}v_{n-3}^2 v_{n-4}\cdots v_1) u_{2n-3}w & \text{if $j=1$}\\
u_3(v_1\cdots v_{n-4}v_{n-3}^2 v_{n-4}\cdots v_1) u_{2n-2}w & \text{if $j=2$}\\
u_{j+1}(v_1\cdots v_{n-4}v_{n-3}^2 v_{n-4}\cdots v_1) wu_{j-2} & \text{if $j=3,\ldots, 2n-3$}\\
wu_1(v_1\cdots v_{n-4}v_{n-3}^2 v_{n-4}\cdots v_1) u_{2n-4}& \text{if $j=2n-2$}.
\end{cases}
\end{equation*}
This yields the remaining equations~\reqref{eq12},~\reqref{eq13},~\reqref{eq14} and~\reqref{eq15}.
\end{proof}

We now simplify somewhat the presentation of $\Gamma_2(B_n(\St))$ given by \repr{fullpres}. From equations~\reqref{eq4} and~\reqref{eq9}, for $i=1,2$ we obtain the following equations:
\begin{align}
u_1v_j &=v_j u_2 \quad\text{for all}\quad j\geq 2\tag{\unreqref[1]{eq4}}\\
u_2v_j &=v_j u_1^{-1}u_2 \quad\text{for all}\quad j\geq 2\tag{\unreqref[2]{eq4}}.
\end{align}
This allows us to eliminate \req{eq5} as follows. For all $j\geq 2$, we have:
\begin{align*}
v_j w u_1w^{-1}v_j^{-1}&= w  v_ju_1^{-1}v_j^{-1}w^{-1} \quad\text{by \req{eq3}}\\
&= w  v_j u_2 v_j^{-1} u_2^{-1}w^{-1}\quad\text{by \req[2]{eq4}}\\
&= w  u_1 u_2^{-1}w^{-1}\quad\text{by \req[1]{eq4}}\\
&= u_{2n-2} \quad\text{by \req{eq11}},
\end{align*}
and this is equivalent to \req{eq5}, which we thus delete from the list of relations.

Suppose that for some $2\leq i\leq 2n-4$, we have equations~\reqref[i-1]{eq4} and~\reqref[i]{eq4}. We now show that they imply \reqref[i+1]{eq4}. For all $j\geq 2$, we have:
\begin{align*}
u_{i+1}v_j u_{i+2}^{-1}v_j^{-1} &= u_{i-1}^{-1}u_i v_j u_{i+1}^{-1}u_i v_j^{-1} \quad\text{by \req{eq9}}\\
&= u_{i-1}^{-1}v_j u_i v_j^{-1} \quad\text{by \req[i]{eq4}}\\
&= 1 \quad\text{by \req[i-1]{eq4}},
\end{align*}
which yields \req[i+1]{eq4}. So we may successively delete equations~\reqref[2n-3]{eq4}, \reqref[2n-4]{eq4}, \ldots, \reqref[3]{eq4} from the list of relations.

We now show that we may delete all but one of the surface relations~\reqref{eq12}--\reqref{eq15}. First suppose that we have \req{eq12}. Now
\begin{align*}
u_{2n-2}wu_3 &=u_{2n-3}wu_1u_3 \quad\text{by \req{eq10}}\\
&= u_{2n-3}wu_2 \quad\text{by \req[1]{eq9}}.
\end{align*}
This implies \req{eq13} which we delete from the list of relations.

Now suppose that we have \req[i+1]{eq14} for some $5\leq i\leq 2n-2$. Let us write $A=v_1\cdots v_{n-4}v_{n-3}^2 v_{n-4}\cdots v_1$. Then $wu_{i-2}u_{i+1}=A^{-1}$. So
\begin{align*}
wu_{i-3}u_i &= wu_{i-2}u_{i+1}\quad\text{by \req[1]{eq9}}\\
&= A^{-1} \quad\text{by above}.
\end{align*}
This yields \req[i]{eq14}, and so we may delete successively equations~\reqref[4]{eq14}, \ldots, \reqref[2n-3]{eq14}. 

Now suppose that we have \reqref{eq15}, so $Au_{2n-4}wu_1=1$. Then 
\begin{align*}
Awu_{2n-5}u_{2n-2} &=Aw u_{2n-4} wu_1w^{-1}\quad\text{by equations~\reqref{eq9} and~\reqref{eq10}}\\
&= w(Au_{2n-4}wu_1)w^{-1}\quad\text{by \req{eq3}}\\
&= 1 \quad\text{by above}.
\end{align*}
This implies \req[2n-2]{eq14} which we delete from the list of relations.

Finally, suppose that we have \req{eq12}. Then
\begin{align*}
Au_{2n-4}wu_1&= Au_{2n-3}u_{2n-2}^{-1}wu_1\quad\text{by \req[2n-4]{eq9}}\\
&= Au_{2n-3}wu_2 \quad\text{by \req{eq11}}\\
& =1  \quad\text{by above}.
\end{align*}
This yields \req{eq15} which we delete from the list. It thus follows that we may delete all but one of the surface relations; let us keep \req{eq12}.

Summing up, we may thus delete relations~\reqref{eq5},~\reqref[i]{eq4} for $i=3,\ldots, 2n-3$ and~\reqref{eq13}--\reqref{eq15} from the presentation of $\Gamma_2(B_n(\St))$ given by \repr{fullpres}.

\section{The derived subgroup of $B_4(\St)$}\label{sec:presg2bn4}

The aim of this section is to use \repr{fullpres} to derive the
presentation of $\Gamma_2(B_4(\St))$ given in
\reth{lcsbn}(\ref{it:lcs4}), from which we were able to see that
$\Gamma_2(B_4(\St))\cong \F[2]\rtimes \mathcal{Q}_8$.

We first remark that in this case, the relations~\reqref{eq1},~\reqref{eq2},~\reqref{eq4} and~\reqref{eq5} do not exist. Further, from relations~\reqref{eq9}, we may obtain the following:
\begin{align*}
u_2 &= u_3u_4^{-1}&
u_1 &= u_3u_4^{-1}u_3^{-1}\\
u_5 &= u_3^{-1}u_4&
u_6 &= u_4^{-1}u_3^{-1}u_4,
\end{align*}
which we take to be definitions of $u_1,u_2,u_5$ and $u_6$, so we delete equation~\reqref{eq9} from the list of relations. From \req{eq12}, we see that
\begin{equation*}
w=u_4^{-1}u_3 v_1^{-2} u_4u_3^{-1}.
\end{equation*}
We conclude that $\Gamma_2(B_4(\St))$ is generated by $u_3,u_4$ and $v_1$.

Let us return momentarily to the situation of the previous section. Before deleting all but one of the surface relations, we shall derive some other useful relations. 

Consider the surface relations~\reqref{eq12}--\reqref{eq15}. From relations~\req{eq12} and~\reqref[5]{eq14} (resp.\ \reqref{eq15}  and~\reqref[4]{eq14}), it follows that $u_5$ (resp.\ $u_4$) commutes with $v_1^2$. But these two equations are equivalent to the relations 
\begin{align}
& u_3\comm v_1^2,\quad\text{and}\label{eq:eq16}\\
& u_4\comm v_1^2.\label{eq:eq17}
\end{align}

Further, equations~\reqref{eq15} and~\reqref{eq17} imply \req[4]{eq14}, and equations~\reqref{eq12},~\reqref{eq16} and~\reqref{eq17} imply \req[5]{eq14}, so we replace equations~\reqref[4]{eq14} and \reqref[5]{eq14} by equations~\reqref{eq16} and~\reqref{eq17}.

As in \resec{gen}, we can then delete equations~\reqref[6]{eq14} and \reqref{eq15} from the list of relations, which becomes: \reqref{eq3}, \reqref{eq6}, \reqref{eq7}, \reqref{eq8}, \reqref{eq10} and \reqref{eq11}. We now analyse these relations in further detail.

From equation~\reqref{eq3} and the definition of $w$, we see that $v_1\comm u_4^{-1}u_3u_4u_3^{-1}$. Up to conjugacy, \req[1]{eq6} may be written as follows:
\begin{align*}
1 &= u_3v_1^{-1}u_3^{-1} u_3 u_4u_3^{-1} v_1^{-1}u_3 u_4^{-1}u_3^{-1} v_1\\
&= u_3v_1^{-1}u_3^{-1} u_4 u_4^{-1} u_3u_4 u_3^{-1} v_1^{-1} u_3 u_4^{-1}u_3^{-1} v_1= u_3v_1^{-1}u_3^{-1} u_4 v_1^{-1}u_4^{-1}v_1,
\end{align*}
and hence we may replace \req[1]{eq6} by:
\begin{equation}
u_3v_1u_3^{-1}= u_4 v_1^{-1}u_4^{-1}v_1.\label{eq:eq61p}
\end{equation}
Up to conjugacy, \req[2]{eq6} may be written:
\begin{equation}
u_3^{-1}v_1u_3= u_4^{-1}v_1u_4 v_1^{-1}.\label{eq:eq62p}
\end{equation}
By equations~\reqref{eq61p} and~\reqref{eq16}, the left-hand side of \req[3]{eq6} may be written:
\begin{equation*}
u_3v_1u_3^{-1}u_4 v_1^{-1}u_4^{-1}v_1^{-1}= u_3v_1^2 u_3^{-1}v_1^{-2}=1,
\end{equation*}
so relation~\reqref[3]{eq6} is automatically satisfied, and we thus delete it from the list.

Using the fact that $v_1\comm u_4^{-1}u_3u_4u_3^{-1}$, \req[4]{eq6} may be written:
\begin{align*}
1 &= u_4v_1 u_4^{-1}u_3^{-1}u_4 v_1^{-1}u_4^{-1} u_3 v_1^{-1}= u_4v_1 u_3^{-1}u_3 u_4^{-1}u_3^{-1}u_4 v_1^{-1}u_4^{-1} u_3 v_1^{-1}\\
&= u_4v_1 u_3^{-1}v_1^{-1} u_3 u_4^{-1} v_1^{-1},
\end{align*}
and from this, we obtain \req{eq62p}, using the fact that $v_1^2$ commutes with $u_4$. So we delete \req[4]{eq6} from the list.

We now consider \req{eq3}. Using equations~\reqref{eq61p} and~\reqref{eq62p}, we obtain:
\begin{align*}
1&= u_4^{-1}u_3 u_4u_3^{-1}v_1 u_3u_4^{-1} u_3^{-1}u_4 v_1^{-1} = u_4^{-1}u_3 v_1u_4 v_1^{-1}u_4^{-1} u_3^{-1}u_4 v_1^{-1}\\
&= u_4^{-1}u_3 v_1u_4 v_1^{-1}u_4^{-1}  v_1u_3^{-1}v_1^{-1} u_4= u_4^{-1}u_3 v_1u_3 v_1u_3^{-2}v_1^{-1} u_4,
\end{align*}
which up to conjugacy, and using the fact that $v_1^2$ commutes with $u_3$ yields:
\begin{equation}
u_3^{-2} v_1^{-1}u_3v_1^{-1}u_3 v_1^{-1}\cdot v_1^{4}=1.\label{eq:eq3p}
\end{equation}
We replace \req{eq3} by this relation.

From equations~\reqref{eq16} and~\reqref{eq17}, the left-hand side of equations~\reqref{eq10} and~\reqref{eq11} collapse, and so we delete them from the list.

After immediate cancellations, \req{eq8} becomes:
\begin{align*}
1&= u_4^{-1} u_3^{-1}u_4 v_1 u_4^{-1} u_3 v_1^{-1} u_4 v_1^{-1}\\
& = u_4^{-1} u_3^{-1}u_4 v_1 u_4^{-1}v_1^{-1} v_1 u_3 v_1^{-1} u_4 v_1 u_4^{-1} u_4 v_1^{-2}\\
&= u_4^{-1} u_3^{-1} u_3 v_1u_3^{-1}v_1 u_3 u_3 v_1^{-1}u_3^{-1}u_4 v_1^{-2},
\end{align*}
which up to conjugacy and inversion yields \req{eq3p}. So we delete \req{eq8} from the list.

After immediate cancellations, the left-hand side of \req{eq7} becomes:
\begin{align*}
u_3^{-1}u_4 v_1 u_4^{-1} u_3 u_4^{-1} u_3^{-1}u_4 v_1^{-1} u_4^{-1} &= u_3 u_4 v_1^{-1} = u_3^{-1}u_4 v_1 u_4^{-1} v_1^{-1} u_3 v_1^{-1}\\
&= u_3^{-1} u_3 v_1 u_3^{-1}u_3 v_1^{-1}=1,
\end{align*}
using the fact that $v_1\comm u_4^{-1}u_3u_4u_3^{-1}$, and applying  equations~\reqref{eq61p} and~\reqref{eq17}. So we delete \req{eq7} from the list.

We are thus left with relations~\reqref{eq16},~\reqref{eq17},~\reqref{eq61p},~\reqref{eq62p} and~\reqref{eq3p}. We now multiply together equations~\reqref{eq61p} and~\reqref{eq62p}. The product of the left-hand sides, by \req{eq3p}, is given by:
\begin{equation*}
u_3v_1u_3^{-2} v_1u_3=v_1,
\end{equation*}
while by equations~\reqref{eq16},~\reqref{eq17},~\reqref{eq61p},~\reqref{eq62p} and~\reqref{eq3p}, the product of the right-hand sides is given by:
\begin{align*}
u_4 v_1^{-1}  u_4^{-1} v_1 u_4^{-1} v_1 u_4 v_1^{-1}&= u_4 v_1^{-1}  u_4^{-1} v_1 u_4^{-1} v_1^{-1} u_4 v_1\\
&= v_1^{-1} v_1 u_4 v_1^{-1}  u_4^{-1} u_3^{-1} v_1^{-1} u_3 v_1\\
&= v_1^{-1} v_1^{-1} u_4 v_1  u_4^{-1} u_3^{-1} v_1^{-1} u_3 v_1\\
&= v_1^{-1} u_3 v_1^{-1} u_3^{-2} v_1^{-1} u_3 v_1=v_1^{-3}.
\end{align*}
From these two equations, we conclude that:
\begin{equation}\label{eq:v14}
v_1^4=1,
\end{equation}
and so \req{eq3p} becomes:
\begin{equation}\label{eq:eq3pp}
u_3^{-2} v_1^{-1}u_3v_1^{-1}u_3 v_1^{-1} =1.
\end{equation}

The list of relations now becomes: \reqref{eq16}, \reqref{eq17}, \reqref{v14}, \reqref{eq3pp}, \reqref{eq61p} and \reqref{eq62p}. We may rewrite the corresponding presentation as follows:

\begin{prop}\label{prop:presg2b4} The following constitutes a
presentation of the group $\Gamma_2(B_4(\St))$:
\begin{enumerate}
\item[\underline{\textbf{generators:}}] $g_1,g_2,g_3$, where in terms of the usual generators of $B_4(\St)$, $g_1=u_3=\sigma_1^2\sigma_2 \sigma_1^{-3}$, $g_2=u_4=\sigma_1^3\sigma_2 \sigma_1^{-4}$ and $g_3=v_1=\sigma_3 \sigma_1^{-1}$.
\item[\underline{\textbf{relations:}}] 
\begin{gather*}
g_3^4=1\\
g_3^2 \comm g_1\\
g_3^2 \comm g_2\\
g_3 \comm g_2g_1\\
g_2^{-1} g_1^{-1} g_3^{-1} g_1 g_2 g_3^{-1}=1\\
g_1^{-2} g_3^{-1} g_1 g_3^{-1} g_1 g_3^{-1}=1.
\end{gather*}
\end{enumerate}
\end{prop}

\begin{proof}
Rewriting $u_3,u_4$ and $v_1$ in terms of the $g_i$, we obtain directly the first three and the last of the given relations.  As for the fourth and fifth relations, we obtain respectively:
\begin{align*}
g_3 g_2g_1 g_3^{-1} g_1^{-1} g_2^{-1} &= v_1 u_4 u_3^{-1} v_1^{-1} u_3^{-1} u_4^{-1}\\
&= u_4(u_4^{-1} v_1u_4 v_1^{-1}v_1 u_3^{-1} v_1^{-1} u_3^{-1}) u_4^{-1}\\
&= u_4u_3(u_3^{-2} v_1^{-1} u_3v_1^{-1} u_3v_1^{-1} v_1^{4})u_3^{-1} u_4^{-1}=1
\end{align*}
by equations~\reqref{eq62p},~\reqref{v14} and~\reqref{eq3pp}, and
\begin{align*}
g_2^{-1} g_1^{-1} g_3^{-1} g_1 g_2 g_3^{-1}&=  u_4^{-1} u_3^{-1} v_1^{-1} u_3 u_4 v_1^{-1}\\
&= v_1u_4^{-1} (u_4v_1^{-1}u_4^{-1}v_1 v_1^{-1} u_3^{-1} v_1^{-1} u_3) u_4 v_1^{-1}\\
&= v_1u_4^{-1}u_3 (v_1 u_3^{-1} v_1^{-1} u_3^{-1} v_1^{-1} u_3^2) u_3^{-1} u_4 v_1^{-1}=1
\end{align*}
by equations~\reqref{eq61p},~\reqref{v14} and~\reqref{eq3pp}. Thus the presentation we derived with generators $u_3,u_4$ and $v_1$ implies the system given by \repr{presg2b4}. Conversely, given this system, we have
\begin{equation*}
u_3v_1u_3^{-1} = u_3^{-1} v_1^{-1} u_3 v_1^{-1}= u_4v_1^{-1} u_4^{-1} v_1,
\end{equation*}
which is \req{eq61p}, and
\begin{align*}
u_3^{-1} v_1u_3 &= u_3^{-1} v_1^{-1} u_3 v_1^2= u_3v_1 u_3^{-1}v_1 v_1^2 = u_4^{-1}v_1u_4 v_1^3= u_4^{-1}v_1u_4 v_1^{-1},
\end{align*}
which is \req{eq62p}. Hence the system given by \repr{presg2b4} is equivalent to our presentation with generators $u_3,u_4$ and $v_1$, and so in particular is a presentation of $\Gamma_2(B_4(\St))$.
\end{proof}

\section{The derived subgroup of $B_5(\St)$}\label{sec:presg2bn5}

For the case $n=5$, we obtain the following presentation directly from \repr{fullpres}:
\begin{prop}\label{prop:g2b5}
The following constitutes a presentation of the group $\Gamma_2(B_5(\St))$:
\begin{enumerate}
\item[\underline{\textbf{generators:}}] 
\begin{gather*}
w=\si[8]1\\
u_1=\si{2}\sii{1}, u_2=\si{1}\si{2}\sii[2]{1}, \ldots, u_8=\si[7]{1} \si{2}\sii[8]{1}\\
v_1=\si{3}\sii{1}, v_2=\si{4}\sii{1}.
\end{gather*}
\item[\underline{\textbf{relations:}}] 
\begin{gather*}
v_1v_2v_1=v_2v_1v_2\\
\text{$w \comm v_i$ for $i=1,2$}\\
u_1v_2 =v_2 u_2\\
u_2v_2 =v_2 u_1^{-1}u_2\\
\text{$u_iv_1 u_{i+2}v_1^{-1}u_{i+1}^{-1}v_1^{-1}=1$ for $i=1,\ldots, 6$}\\
u_7v_1 wu_1 w^{-1}v_1^{-1}u_8^{-1}v_1^{-1}=1\\
u_8v_1 wu_2 v_1^{-1}u_1^{-1}w^{-1}v_1^{-1}=1\\
\text{$u_{i+1}u_{i+2}^{-1}u_i^{-1}=1$ for $i=1,\ldots, 6$}\\
u_8wu_1^{-1}w^{-1}u_7^{-1}=1\\
wu_1u_2^{-1}w^{-1}u_8^{-1}=1\\
u_2(v_1v_2^2 v_1) u_7w=1.\tag*{\mbox{\qed}}
\end{gather*}
\end{enumerate}
\end{prop}

\section{The derived subgroup of $B_n(\St)$ for $n\geq 6$}\label{sec:presg2bn6}

We now suppose that $n\geq 6$. Then the generator $v_3$ exists.

Suppose that \req[i]{eq6} holds for some $1\leq i\leq 2n-5$. Let us take $j\geq 3$. We eliminate \req[i+1]{eq6} as follows: applying successively equations~\reqref{eq4} and~\reqref{eq1}, we obtain:
\begin{align*}
u_{i+1}v_1 u_{i+3}v_1^{-1}u_{i+2}^{-1}v_1^{-1} &= v_j^{-1}u_i v_j  v_1 v_j^{-1}u_{i+2} v_j v_1^{-1}
v_j^{-1}u_{i+1}^{-1} v_j v_1^{-1}\\
&= v_j^{-1}(u_i  v_1 u_{i+2}  v_1^{-1}
u_{i+1}^{-1}  v_1^{-1})v_j\\
&=1 \;\text{by \req[i]{eq6}}.
\end{align*}
It thus follows that we may delete successively equations~\reqref[2n-4]{eq6}, \ldots, \reqref[2]{eq6} from the list of relations.

Suppose that \req{eq7} holds. Applying the idea of the previous paragraph, we eliminate \req{eq8}:
\begin{align*}
u_{2n-2}v_1 wu_2 v_1^{-1}u_1^{-1}w^{-1}v_1^{-1} =& v_j^{-1}u_{2n-3} v_j v_1 w v_j^{-1}u_1 v_j v_1^{-1} w^{-1} v_j^{-1}\cdot\\
& u_{2n-2}^{-1} v_j v_1^{-1}\\
=&  v_j^{-1}(u_{2n-3}  v_1 w u_1 v_1^{-1} w^{-1} u_{2n-2}^{-1}  v_1^{-1})v_j=1.
\end{align*}

Let us suppose that \req[1]{eq6} holds. Then so does \req[2n-4]{eq6}. We eliminate \req{eq7} as follows.
\begin{align*}
u_{2n-3}v_1 wu_1 w^{-1}v_1^{-1}u_{2n-2}^{-1}v_1^{-1}=&  v_j^{-1}u_{2n-4} v_j v_1 v_j^{-1}u_{2n-2} v_j v_1^{-1}  v_j^{-1}\cdot\\
&u_{2n-3}^{-1} v_j v_1^{-1}\\
=& v_j^{-1}(u_{2n-4}  v_1 u_{2n-2} v_1^{-1} u_{2n-3}^{-1}  v_1^{-1})v_j=1.
\end{align*}

\begin{prop}
Let $n\geq 6$. The following constitutes a presentation of the group $\Gamma_2(B_n(\St))$:
\begin{enumerate}
\item[\underline{\textbf{generators:}}] 
\begin{gather*}
w=\si[2n-2]1\\
u_1=\si{2}\sii{1}, u_2=\si{1}\si{2}\sii[2]{1}, \ldots, u_{2n-2}=\si[2n-3]{1} \si{2}\sii[(2n-2)]{1}\\ v_1=\si{3}\sii{1}, \ldots, v_{n-3}=\si{n-1}\sii{1}.
\end{gather*}
\item[\underline{\textbf{relations:}}] 
\begin{gather}
\text{$v_iv_j=v_jv_i$ if $\lvert i-j\rvert\geq 2$}\label{eq:eqa}\\
\text{$v_i v_{i+1}v_i=v_{i+1}v_iv_{i+1}$ for all $1\leq i<j\leq n-4$}\label{eq:eqb}\\
w \comm v_i\label{eq:eqc}\\
\text{$u_1v_j =v_j u_2$, where $j\geq 2$}\label{eq:eqd}\\
\text{$u_2v_j =v_j u_1^{-1}u_2$, where $j\geq 2$}\label{eq:eqe}\\
u_1v_1 u_1^{-1}u_2 v_1^{-1}u_2^{-1}v_1^{-1}=1\label{eq:eqf}\\
\text{$u_{i+1}u_{i+2}^{-1}u_i^{-1}=1$ for all $i=1,\ldots, 2n-4$}\label{eq:eqg}\\
u_{2n-2}wu_1^{-1}w^{-1}u_{2n-3}^{-1}=1\label{eq:eqh}\\
wu_1u_2^{-1}w^{-1}u_{2n-2}^{-1}=1\label{eq:eqi}\\
u_2(v_1\cdots v_{n-4}v_{n-3}^2 v_{n-4}\cdots v_1) u_{2n-3}w=1.\label{eq:eqj}
\end{gather}\qed
\end{enumerate}
\end{prop}

This presentation may be refined further. Set 
\begin{equation*}
\text{$A=v_1\cdots v_{n-4}v_{n-3}^2 v_{n-4}\cdots v_1$ and $y=u_2^{-1} u_1u_2 u_1^{-1}$.}
\end{equation*} 
Applying equations~\reqref{eqh} and~\reqref{eqi} to \req{eqj}, we have:
\begin{equation*}
1=u_2A u_{2n-3}w=u_2A u_{2n-2}w u_1^{-1}=u_2 Aw u_1u_2^{-1} u_1^{-1},
\end{equation*}
so
\begin{equation*}
w=A^{-1} y.
\end{equation*}
Since $A$ commutes with $w$ by \req{eqc}, we see that $A$ commutes with $y$. Equations~\reqref{eqh} and~\reqref{eqi} are then equivalent to:
\begin{align}
u_{2n-3} &= A^{-1} u_2^{-1} y^{-1}A\label{eq:eqhp}\\
u_{2n-2} &= A^{-1} u_2^{-1} u_1y^{-1} A.\label{eq:eqip}
\end{align}

Let $i\geq 2$. One may check using relations~\reqref{eqa} and~\reqref{eqb} that $A$ commutes with $v_i$. Relation~\reqref{eqc} is then equivalent to $v_i$ commutes with $y$. But this is implied by equations~\reqref{eqd} and~\reqref{eqe}. Indeed, from these two relations we see that $v_iu_2v_i^{-1}=u_1$ and $v_iu_1v_i^{-1}=u_1u_2^{-1}$, and then one may check directly that $v_i$ commutes with $u_2^{-1} u_1u_2 u_1^{-1}$. This implies that we may delete equations~\reqref[i]{eqc} for $2\leq i\leq n-3$.

From \req{eqg}, we may calculate $u_3,\ldots u_{2n-4},u_{2n-3}$ and $u_{2n-2}$ in terms of $u_1$ and $u_2$. Since all but the last two of these elements do not appear anywhere in the rest of the presentation, we may delete relations~\reqref[i]{eqg} for $i=1,\ldots, 2n-4$, provided that we keep (as definitions) the expressions for $u_{2n-3}$ and $u_{2n-2}$ in terms of $u_1$ and $u_2$. Let us calculate the general term $u_i$ in terms of $u_1$ and $u_2$.

For $i\in\N$, we define $v_i$ as follows:
\begin{equation*}
v_i= 
\begin{cases}
u_1u_2^{-1} & \text{if $i\equiv 0\bmod 6$}\\
u_1 & \text{if $i\equiv 1\bmod 6$}\\
u_2 & \text{if $i\equiv 2\bmod 6$}\\
u_1^{-1}u_2 & \text{if $i\equiv 3\bmod 6$}\\
u_1^{-1} & \text{if $i\equiv 4\bmod 6$}\\
u_2^{-1} & \text{if $i\equiv 5\bmod 6$.}
\end{cases}
\end{equation*}

\begin{lem}
Let $i\in\N$, and let $k\geq 0$ and $0\leq l\leq 5$ be such that $i=6k+l+1$. Then:
\begin{equation*}
u_i= 
\begin{cases}
y^k v_i y^{-k} & \text{if $l= 0,1,2$}\\
y^k u_2^{-1}u_1 v_i u_1^{-1}u_2 y^{-k} & \text{if $l=3,4,5$.}
\end{cases}
\end{equation*}
\end{lem}

\begin{proof}
The proof is by induction on $i$, one considers the six possible cases depending on the value of $i\bmod 6$.
\end{proof}

We can then determine equations~\reqref{eqhp} and~\reqref{eqip} in the three possible cases. We let $k\geq 0$ and $0\leq l\leq 5$ be such that $2n-2=6k+l+1$.
\begin{enumerate}
\item\label{it:pres6a} \underline{$2n-2\equiv 0\bmod 6$ ($l=5$):}
\begin{align*}
y^k u_2^{-1}  y^{-k} &= A^{-1} u_2^{-1} A\\
y^k u_2^{-1} u_1 y^{-k} &= A^{-1} u_2^{-1} u_1 A.
\end{align*}
Hence:
\begin{equation*}
u_1,u_2\comm Ay^k.
\end{equation*}
\item\label{it:pres6b} \underline{$2n-2\equiv 2\bmod 6$ ($l=1$):}
\begin{align*}
y^k u_1 y^{-k} &= A^{-1} u_2^{-1} y^{-1}A\\
y^k u_2 y^{-k} &= A^{-1} u_2^{-1} u_1 y^{-1}A.
\end{align*}
\item\label{it:pres6c} \underline{$2n-2\equiv 4\bmod 6$ ($l=3$):}
\begin{align*}
y^k u_1^{-1} y^{-k} &= A^{-1} u_1^{-1} u_2 A\\
y^k u_2 y^{-k} &= A^{-1} u_2^{-1}u_1^{-1} u_2 A.
\end{align*}
\end{enumerate}

\begin{prop}\label{prop:g2b6}
Let $n\geq 6$. The following constitutes a presentation of the group $\Gamma_2(B_n(\St))$:
\begin{enumerate}
\item[\underline{\textbf{generators:}}] 
\begin{gather*}
u_1=\si{2}\sii{1}, u_2=\si{1}\si{2}\sii[2]{1}\\ 
v_1=\si{3}\sii{1}, \ldots, v_{n-3}=\si{n-1}\sii{1}.
\end{gather*}
\item[\underline{\textbf{relations:}}] 
\begin{gather*}
\text{$v_iv_j=v_jv_i$ if $\lvert i-j\rvert\geq 2$}\\
\text{$v_i v_{i+1}v_i=v_{i+1}v_iv_{i+1}$ for all $1\leq i<j\leq n-4$}\\
y \comm v_1\\
v_j u_2v_j^{-1}= u_1, \quad\text{where $j\geq 2$}\\
v_j u_1v_j^{-1}= u_1u_2^{-1}, \quad\text{where $j\geq 2$}\\
u_1v_1 u_1^{-1}u_2 v_1^{-1}u_2^{-1}v_1^{-1}=1,
\end{gather*}
plus the two corresponding relations from~(\ref{it:pres6a}),~(\ref{it:pres6b}) and~(\ref{it:pres6c})  of the previous paragraph, where 
\begin{equation*}
\text{$y=u_2^{-1} u_1u_2 u_1^{-1}$ and $A=v_1\cdots v_{n-4}v_{n-3}^2 v_{n-4}\cdots v_1$.}\tag*{\mbox{\qed}}
\end{equation*}
\end{enumerate}
\end{prop}

\begin{rem}
From this presentation, one could also delete, for example, the generator $u_2$.
\end{rem}

\backmatter

\bibliographystyle{amsalpha}

\begin{thebibliography}{[McCa]}
{\small

\bibitem[AW]{AW} W.~A.~Adkins and S.~Weintraub,  Algebra, an approach
via module theory, Graduate Texts in Mathematics, \textbf{136},
Springer-Verlag, New York, 1992.

\bibitem[All]{All} D.~Allcock, Braid pictures for Artin groups, \emph{Trans.\ Amer.\ Math.\ Soc.} \textbf{354} (2002), 3455--3474.

\bibitem[A1]{A1} E.~Artin, Theorie der Z\"opfe, \emph{Abh.\ Math.\
Sem.\ Univ.\ Hamburg} \textbf{4} (1925), 47--72.

\bibitem[A2]{A2} E.~Artin, Theory of braids, \emph{Ann.\ Math.}
\textbf{48} (1947), 101--126.

\bibitem[A3]{A3} E.~Artin, Braids and permutations, \emph{Ann.\
Math.} \textbf{48} (1947), 643--649.

\bibitem[BS]{BS} G.~Baumslag and D.~Solitar, Some two-generator one-relator non-Hopfian groups, \emph{Bull.\ Amer.\ Math.\ Soc.} \textbf{68} (1962),199--201.

\bibitem[BGG]{BGG} P.~Bellingeri, S.~Gervais and J.~Guaschi, Lower central series for surface braid groups, arXiv preprint math.GT/0512155.

\bibitem[BG]{BG} P.~Bellingeri and E.~Godelle, Questions on surface
braid groups, Preprint LMNO~2005-15 University of Caen, arXiv preprint math.GR/05036578.

\bibitem[Big]{Big} S.~Bigelow, Braid groups are linear, \emph{J.\ Amer.\ Math.\ Soc.} \textbf{14} (2001), 471--486

\bibitem[Bi1]{Bi1} J.~S.~Birman, On braid groups, \emph{Comm.\ Pure
and Appl.\ Math.} \textbf{22} (1969), 41--72.

\bibitem[Bi2]{Bi2} J.~S.~Birman, Braids, links and mapping class
groups, \emph{Ann.\ Math.\ Stud.} \textbf{82}, Princeton University
Press, 1974.

\bibitem[Bi3]{Bi3} J.~S.~Birman, Mapping class groups of surfaces, in
Braids (Santa Cruz, CA, 1986), 13--43, \emph{Contemp.\ Math.}
\textbf{78}, Amer.\ Math.\ Soc., Providence, RI, 1988. 

\bibitem[Bri]{Bri} E.~Brieskorn, \emph{Sur les groupes de tresses 
(d'apr\`es V.~I.~Arnol'd)}, S\'eminaire Bourbaki, 24\`eme ann\'ee 
(1971/1972), Exp.\ No.\ 401, Lecture Notes in Mathematics \textbf{317}, Springer, Berlin, 1973, 21--44. 

\bibitem[Bro]{Bro} K.~S.~Brown, Cohomology of groups, Graduate Texts in Mathematics, \textbf{87}, Springer-Verlag, New York, 1982.

\bibitem[BZ]{BZ} G.~Burde and H.~Zieschang,
Knots, Second edition, de Gruyter Studies in Mathematics, \textbf{5}.
Walter de Gruyter \& Co., Berlin, 2003.


\bibitem[ChP]{ChP} R.~Charney and D.~Peifer, The $K(\pi,1)$-conjecture
for the affine braid groups, \emph{Comment.\ Math.\ Helv.} \textbf{78}
(2003), 584--600.

\bibitem[CFL]{CFL} K.~T.~Chen, R.~H.~Fox and R.~C.~Lyndon, \emph{Ann.\
Math.} \textbf{68} (1958), 81--95.

\bibitem[Ch]{Ch} W.-L.~Chow, On the algebraical braid group. 
\emph{Ann.\ Math.} \textbf{49} (1948) 654--658.

\bibitem[CG]{CG} F.~R.~Cohen and S.~Gitler, On loop
spaces of configuration spaces, \emph{Trans.\ Amer.\
Math.\ Soc.} \textbf{354} (2002), 1705--1748.

\bibitem[Cr]{Cr} J.~Crisp, Injective maps between Artin groups, in Geometric group theory down under, 1996, 119--137, Eds.\ J.~Cossey, C.~F.~Miller~III, W.~D.~Neumann, M.~Shapiro, de Gruyter, 1999.

\bibitem[CrP]{CP} J.~Crisp and L.~Paris, Artin groups of type~B and~D, \emph{Adv.\ Geom.} \textbf{5} (2005), 607--636.

\bibitem[Fa]{Fa} E.~Fadell, Homotopy groups of configuration spaces
and the string problem of Dirac, \emph{Duke Math.\ Journal}
\textbf{29} (1962), 231--242.

\bibitem[FH]{FH} E.~Fadell and S.~Y.~Husseini,  Geometry and topology
of configuration spaces, Springer Monographs in Mathematics.
Springer-Verlag, Berlin, 2001.

\bibitem[FaN]{FaN} E.~Fadell and L.~Neuwirth, Configuration spaces, \emph{Math.\ Scandinavica} \textbf{10} (1962), 111--118.

\bibitem[FVB]{FvB} E.~Fadell and  J.~Van~Buskirk,
The braid groups of $\mathbb{E}^2$ and $\St$, \emph{Duke
Math.\ Journal} \textbf{29} (1962), 243--257.

\bibitem[FR1]{FR1} M.~Falk and R.~Randell,  The lower central series
of a fiber-type arrangement, \emph{Invent.\ Math.} \textbf{82} (1985),
77--88.

\bibitem[FR2]{FR2} M.~Falk and R.~Randell,  The lower central series
of generalized pure braid groups, in Geometry and topology (Athens,
Ga., 1985), 103--108, Lecture Notes in Pure and Appl.\ Math.\
\textbf{105}, Dekker, New York, 1987. 

\bibitem[FR3]{FR3} M.~Falk and R.~Randell, Pure braid groups and
products of free groups, Braids (Santa Cruz, CA, 1986), 217--228,
\emph{Contemp.\ Math.} \textbf{78}, Amer.\ Math.\ Soc., Providence,
RI, 1988.

\bibitem[FG]{FG} A.~Fel'shtyn and D.~L.~Gon\c{c}alves, Twisted conjugacy classes of automorphisms of Baumslag-Solitar groups, arXiv preprint math.GR/0405590.

\bibitem[FoN]{FoN} R.~H.~Fox and L.~Neuwirth, The braid groups,
\emph{Math.\ Scandinavica} \textbf{10} (1962), 119--126.

\bibitem[Ga]{Ga} A.~M.~Gaglione, Factor groups of the lower central series for special free products, \emph{J.~Algebra} \textbf{37} (1975), 172--185.

\bibitem[GVB]{GVB} R.~Gillette and J.~Van Buskirk, The word
problem and consequences for the braid groups and mapping class
groups of the $2$-sphere, \emph{Trans.\ Amer.\ Math.\ Soc.}
\textbf{131} (1968), 277--296. 

\bibitem[GG1]{GG1} D.~L.~Gon\c{c}alves and J.~Guaschi, On the
structure of surface pure braid groups, \emph{J.~Pure Appl.\ Algebra} \textbf{182} (2003), 33--64 (due to a printer's error, this article was republished in its entirety with the reference \textbf{186} (2004) 187--218).

\bibitem[GG2]{GG2} D.~L.~Gon\c{c}alves and J.~Guaschi, The roots of
the full twist for surface braid groups, \emph{Math.\ Proc.\ Camb.\ Phil.\ Soc.} \textbf{137} (2004), 307--320.

\bibitem[GG3]{GG3} D.~L.~Gon\c{c}alves and J.~Guaschi, The braid groups of the projective plane, \emph{Algebraic and Geometric Topology} \textbf{4} (2004), 757--780.

\bibitem[GG4]{GG4} D.~L.~Gon\c{c}alves and J.~Guaschi, The braid group
$B_{{n,m}}(\St)$ and the generalised Fadell-Neuwirth short exact
sequence, \emph{J.~Knot Theory and its Ramifications} \textbf{14}
(2005), 375--403.

\bibitem[GG5]{GG5} D.~L.~Gon\c{c}alves and J.~Guaschi, The quaternion group as a subgroup of the sphere braid groups, arXiv preprint math.GT/0603377.

\bibitem[GMP]{GMP} J.~Gonz\'alez-Meneses and L.~Paris, Vassiliev invariants for braids on surfaces, \emph{Trans.\ Amer.\ Math.\ Soc.}  \textbf{356} (2004), 219--243

\bibitem[GL]{GL} E.~A.~Gorin and V~Ja.~Lin, Algebraic equations with continuous coefficients and some problems of the algebraic theory of braids, \emph{Math.~USSR Sbornik} \textbf{7} (1969), 569--596.

\bibitem[Gu]{Gu} J.~Guaschi, Nielsen theory, braids and fixed points
of surface homeomorphisms, \emph{Topology and its Applications}
\textbf{117} (2002), 199--230.

\bibitem[Hal]{H} M.~Hall, The theory of groups, Macmillan, New York, 1959.

\bibitem[Han]{Han} V.~L.~Hansen, Braids and Coverings: Selected
topics, \emph{London Math. Society Student Text}, \textbf{18},
Cambridge University Press, 1989.

\bibitem[HMR]{HMR} P.~Hilton, G.~Mislin and L.~Roitberg, Localization of nilpotent groups and spaces, North-Holland Mathematics Studies, No.\ 15, Notas de Matem\'atica, No.\ 55, North-Holland Publishing Co., Amsterdam, 1975.

\bibitem[J]{J} D.~L.~Johnson, Presentation of groups, LMS Lecture Notes \textbf{22} (1976), Cambridge University Press.

\bibitem[KP]{KP} R.~P.~Kent~IV and D.~Peifer, A geometric and algebraic description of annular braid groups, \emph{Int.~J.~Algebra and Computation} \textbf{12} (2002) 85--97

\bibitem[KMc]{KMc} G.~Kim and J.~McCarron, Some residually $p$-finite one-relator groups, \emph{J.~Algebra} \textbf{169} (1994), 817--826.

\bibitem[Ko]{Ko} T.~Kohno,  S\'erie de Poincar\'e-Koszul associ\'ee
aux groupes de tresses pures,  \emph{Invent.\ Math.} \textbf{82}, 57--75.

\bibitem[Lab]{La} J.~P.~Labute, The lower central series of the group $\setang{x,y}{x^p=1}$, \emph{Proc.\ Amer.\ Math.\ Soc.} \textbf{66} (1977), 197--201.

\bibitem[Lam]{Lam} S.~Lambropoulou,  Braid structures in knot complements, handlebodies and 3-manifolds, in Knots in Hellas~'98 (Delphi), 274--289, Ser.\ Knots Everything, \textbf{24}, World Sci.\ Publishing, River Edge, NJ, 2000.

\bibitem[MKS]{MKS} W.~Magnus, A.~Karrass and D.~Solitar, Combinatorial group theory, Second revised edition, Dover Publications Inc., New York, 1976.

\bibitem[Ma]{Ma} S.~Manfredini, Some subgroups of Artin's braid group, \emph{Topology and its Applications} \textbf{78} (1997), 123--142.

\bibitem[McCa]{McC} J.~McCarron, Residually nilpotent one-relator groups with nontrivial centre, \emph{Proc.\ Amer.\ Math.\ Soc.} \textbf{124} (1996), 197--201.

\bibitem[McCl]{McCl} J.~McCleary, A user's guide to spectral sequences, Second edition. Cambridge Studies in Advanced Mathematics, Vol.~58, Cambridge University Press, 2001.

\bibitem[Mu]{M} K.~Murasugi, Seifert fibre spaces and braid groups, \emph{Proc.\ London Math.\ Soc.} \textbf{44} (1982), 71--84.

\bibitem[MK]{MK} K.~Murasugi and B.~I.~Kurpita,  A study of braids,
Mathematics and its Applications \textbf{484}, Kluwer Academic
Publishers, Dordrecht, 1999.

\bibitem[PR]{PR} L.~Paris and D.~Rolfsen, Geometric subgroups of
surface braid groups \emph{Ann.\ Inst.\ Fourier} \textbf{49}
(1999), 417--472.

\bibitem[R]{R} D.~Rolfsen,  New developments in the theory of Artin's
braid groups, Proceedings of the Pacific Institute for the
Mathematical Sciences Workshop `Invariants of Three-Manifolds'
(Calgary, AB, 1999),  \emph{Topology and its Applications}
\textbf{127} (2003), 77--90.

\bibitem[Sc]{S} G.~P.~Scott, Braid groups and the group of homeomorphisms of a surface, \emph{Proc.\ Camb.\ Phil.\ Soc.} \textbf{68} (1970), 605--617.

\bibitem[St]{St} J.~Stallings, Homology and central series of groups,
\emph{J.~Algebra} \textbf{2} (1965), 170--181.

\bibitem[T]{T} J.~Tits, Normalisateurs de tores~I: groupes de Coxeter \'etendus, \emph{J.~Algebra} \textbf{4} (1966), 96--116.

\bibitem[VB]{vB} J.~Van~Buskirk, Braid groups of compact
$2$-manifolds with elements of finite order, \emph{Trans.\ Amer.\
Math.\ Soc.} \textbf{122} (1966), 81--97.

\bibitem[Z1]{Z1} O.~Zariski,  On the Poincar\'e group of rational
plane curves, \emph{Amer.\ J.\ Math.} \textbf{58} (1936), 607--619.

\bibitem[Z2]{Z2} O.~Zariski, The topological discriminant group of a
Riemann surface of genus $p$, \emph{Amer.\ J.\ Math.} \textbf{59}
(1937), 335--358.
}

\end{thebibliography}

\end{document}